\documentclass{amsart}

\usepackage{amsmath, amssymb, amsthm, amscd}
\usepackage{calc}
\usepackage{psfrag,graphicx,color}
\usepackage{url}
\usepackage{listings}
\usepackage{pdflscape}
\usepackage{ifpdf}

\ifpdf
\usepackage[notightpage]{pst-pdf}
\fi

\usepackage{verbatim}
\usepackage[all]{xy}
\newdir{ >}{{}*!/-10pt/@{>}}
\newdir^{ (}{{}*!/-5pt/@^{(}}
\newdir_{ (}{{}*!/-5pt/@_{(}}
\def\cgaps#1{}
\def\Cgaps#1{}
\allowdisplaybreaks

\def\undersetbrace#1\to#2{\underbrace{#2}_{#1}}
\def\oversetbrace#1\to#2{\overbrace{#2}^{#1}}
\def\AMSunderset#1\to#2{\underset{#1}{#2}}
\def\AMSoverset#1\to#2{\overset{#1}{#2}}

\def\norm#1{\left\|{#1}\right\|}
\def\gs{\:\begin{smallmatrix}< \\ > \end{smallmatrix}\:}
\swapnumbers

\newtheorem*{prop*}{Proposition}

\newtheorem*{thm*}{Theorem}

\newtheorem*{lem*}{Lemma}

\newtheorem*{cor*}{Corollary}
\numberwithin{equation}{subsection}

\def\ign#1{}             
\def\o{\circ}
\def\X{\mathfrak X}
\def\al{\alpha}
\def\be{\beta}

\def\ze{\zeta}

\def\th{\theta}
\def\io{\iota}
\def\ka{\kappa}

\def\ph{\varphi}

\def\om{\omega}
\def\Ga{\Gamma}
\def\De{\Delta}
\def\Th{\Theta}
\def\La{\Lambda}

\def\Ph{\Phi}

\def\Om{\Omega}
\def\i{^{-1}}
\def\x{\times}
\def\p{\partial}
\let\on=\operatorname
\def\L{\mathcal L}
\def\ol{\overline}
\def\grad{\on{grad}}%
\def\AMSonly#1{}
\def\bp{\on{\boxplus}}
\def\id{\on{id}}
\def\R{\mathbb{R}}
\def\Tr{\on{Tr}}
\def\vol{\on{vol}}
\def\Imm{\on{Imm}}
\def\Emb{\on{Emb}}
\def\g{\bar{g}}
\def\grad{\on{grad}}
\def\Vol{\on{Vol}}
\def\Diff{\on{Diff}}
\def\dist{{\on{dist}}}
\def\Id{\on{Id}}
\def\hor{\on{hor}}
\def\Nor{\on{Nor}}
\def\Fl{\on{Fl}}
\def\Len{\on{Len}}
\def\m{\mu} 
\def\v{v} 

\begin{document}
\title{Almost local metrics on shape space of hypersurfaces in $n$-space}
\author{Martin Bauer, Philipp Harms, Peter W. Michor}
\address{
Martin Bauer, Peter W. Michor:
Fakult\"at f\"ur Mathematik, Universit\"at Wien,
Nordbergstrasse 15, A-1090 Wien, Austria.
}
\email{Bauer.Martin@univie.ac.at}
\email{Peter.Michor@univie.ac.at}
\address{Philipp Harms: EdLabs, Harvard University, 44 Brattle Street, Cambridge, MA 02138, USA}
\email{pharms@edlabs.harvard.edu}

\thanks{All authors were supported by 
FWF Project 21030 and by NSF-FRG grant DMS-0456253} 
\date{\today}
\begin{abstract}
This paper extends parts of the results from [P.W.Michor and D. Mumford, \emph{Appl. Comput. Harmon. Anal.,} 23 (2007), pp. 74--113] for plane curves to the case of hypersurfaces in 
$\mathbb R^n$.
Let $M$ be a compact connected oriented $n-1$ dimensional manifold without 
boundary like the sphere or the torus. 
Then shape space is either the manifold of submanifolds of $\mathbb R^n$ of
type $M$, or the orbifold of immersions from $M$ to $\mathbb R^n$
modulo the group of diffeomorphisms of $M$.
We investigate almost local Riemannian metrics on shape space. These are induced by 
metrics of the following form on the space of immersions:
$$ G_f(h,k) = \int_{M} \Ph(\on{Vol}(f),\on{Tr}(L))\g(h, k) \on{vol}(f^*\g),$$
where $\g$ is the Euclidean metric on $\mathbb R^n$, $f^*\g$ is the induced metric on $M$,
$h,k\in C^\infty(M,\mathbb R^n)$ are tangent vectors at $f$ to the space of embeddings or 
immersions, 
where $\Ph:\mathbb R^2\to \mathbb R_{>0}$ is a suitable smooth function, 
$\on{Vol}(f) = \int_M\on{vol}(f^*\g)$ is the total hypersurface volume of $f(M)$, and the 
trace $\on{Tr}(L)$ of the Weingarten mapping is the mean curvature. 
For these metrics we compute the geodesic equations both on the space of immersions and on shape 
space, the conserved momenta arising from the obvious symmetries, and the sectional curvature.
For special choices of $\Ph$ we give complete formulas for the sectional curvature.
Numerical experiments illustrate the behavior of these metrics. 
\end{abstract}

\noindent{\small 
SIAM J. Imaging Sci. 5 (2012), pp. 244-310. 
}

\maketitle

 \section*{Contents}

\noindent
{\hphantom{9}1} {Introduction}
\\
{\hphantom{9}2} {Shape space and the Hamiltonian approach}
\\
{\hphantom{9}3} {Differential geometry of surfaces and notation}
\\
{\hphantom{9}4} {Variational formulas}
\\
{\hphantom{9}5} {The geodesic equation on $\operatorname {Imm}(M,\mathbb R^n)$}
\\
{\hphantom{9}6} {The geodesic equation on $B_i(M,\mathbb R^n)$}
\\
{\hphantom{9}7} {Sectional curvature on $B_i(M,\mathbb R^n)$}
\\
{\hphantom{9}8} {Geodesic distance on $B_i(M,\mathbb R^n)$}
\\
{\hphantom{9}9} {The set of concentric spheres}
\\
{10} {Special cases of almost local metrics}
\\
{11} {Numerical results}
\\
\hphantom{90} {References}
\\
\hphantom{90} {The AMPL model file}

\markboth{Almost local metrics on shape space of hypersurfaces in $n$-space }{M. Bauer, P. Harms, P. Michor}
\newpage
\section{INTRODUCTION}\label{sec:introduction}{\hfill}\par
\smallskip

Many procedures in science, engineering, and medicine produce data in the form of shapes of point clouds in 
$\mathbb R^n$. If one expects such a cloud to follow roughly a submanifold of a certain type in 
$\mathbb R^n$, then it is of utmost importance to describe the space of all possible submanifolds 
of this type (we call it a shape space hereafter) 
and equip it with a significant metric which is able to distinguish special features 
of the shapes. Almost local metrics are a contribution towards this aim. 

This paper benefited from discussions with David Mumford,  Hermann Schichl, who taught us 
about the use of AMPL, Johannes Wallner and Tilak Ratnanather.
Parts of this paper can be found in the Ph.D. thesis of Martin Bauer \cite{Bauer2010}.

\subsection{Reading suggestions}
A reader who wants to see results before immersing himself in the theoretical background
is recommended to pick up the necessary definitions in the introduction and to jump directly
to the last two sections containing special cases and numerical results.
In section \ref{hamiltonian} we build the fundaments for  shape analysis in a Riemannian setting. 
Section \ref{no} presents background material in differential geometry 
and can serve  as a reference for our notation.
Throughout this work we will use covariant derivatives of vector fields along immersions. This 
concept might not be known to all readers; thus we have decided to give a careful description 
in section \ref{no:co}. 
The main results of the work are in sections \ref{geodesic_equation_Bi}--\ref{numerics}.

In the following introduction we give a non-technical presentation of our approach.
Parts of it can also be found in the Ph.D. thesis of 
Philipp Harms \cite{Harms2010}. 

\subsection{The Riemannian setting}
Most of the metrics used today in data analysis and computer vision are of an ad-hoc 
and naive nature. One embeds shape space in some Hilbert space or Banach space and uses the 
distance therein. Shortest paths are then line segments, but they leave shape space quickly. 
For several reasons the Riemannian setting for shape analysis is a better 
solution. 
\begin{itemize}
\item It formalizes an \emph{intuitive notion of similarity} of shapes: Shapes that
differ only by a small \emph{deformation} are similar to each other. To compare shapes, 
we measure the length of a deformation. A deformation of a shape is a path in shape space. 
Remember that in a Riemannian manifold, 
the geodesic distance between two points is the infimum over the length of all 
paths connecting them. 
\item Riemannian metrics on shape space have been used successfully in 
\emph{computer vision} for a long time, often without any mention of the 
underlying metric. Gradient flows for shape smoothing are an example.  
An underlying metric is needed for the definition of a gradient.
Often, the metric  used implicitly is the $L^2$-metric
which has, however, turned out to be too weak.
\item The exponential map (if it exists) that is induced by a Riemannian metric permits us to
\emph{linearize shape space}: When shapes are represented as initial velocities of 
geodesics connecting them to a fixed reference shape, one effectively works
in the linear tangent space over the reference shape. Curvature will play an 
essential role in quantifying the deviation of curved shape space from its 
linearized approximation. 
\item The linearization of shape space by the exponential map
allows us to do \emph{statistics} on shape space. 
\end{itemize}
However a disadvantage of the Riemannian approach is that shapes can be compared with each other
only when there is a deformation between them, i.e., when they have the same topology. 

\subsection{Shape spaces}
In mathematics and computer vision, shapes have been represented in many ways.
Point clouds, meshes, level-sets, graphs of a function, currents, and measures are
but some of the possibilities. The notion of shapes underlying this work is that 
of immersed or embedded submanifolds of $\R^n$ of co-dimension one.
Any such submanifold  can be represented as a fixed immersion or embedding modulo 
reparametrizations. 

\begin{figure}[ht]
\centering
\includegraphics[width=.72\textwidth]{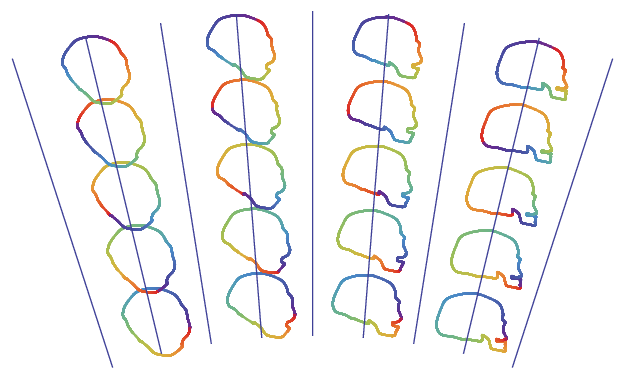}
\caption{Illustration of the space of immersions of $S^1$ in $\R^2$. }
\label{immersions}
\end{figure}

The space of all immersions is illustrated in figure~\ref{immersions}.
The colors are there to help the reader  visualize parametrizations. 
Immersions differing only by a reparametrization
are drawn along vertical lines. These lines are the orbits of 
the reparametrization group. Immersions in the same orbit correspond to 
the same shape in shape space. In other words, shape space is the
space of orbits of the reparametrization group acting on the space of immersions.

As mentioned
in the previous section, only shapes with the same topology  can be compared. Thus
we assume that all shapes (i.e., submanifolds) are diffeomorphic to the same compact 
connected oriented $n-1$ dimensional manifold $M$. 
We will deal only with smooth shapes. They form the core of actual shape space, 
which can be viewed as the Cauchy completion with respect to geodesic distance 
for one of the Riemannian metrics that we treat in this paper. 
See section~\ref{hamiltonian} for a formal definition 
of shape space.

\subsection{Riemannian metrics on shape space}
\begin{figure}[ht]
\centering
\includegraphics[width=.9\textwidth]{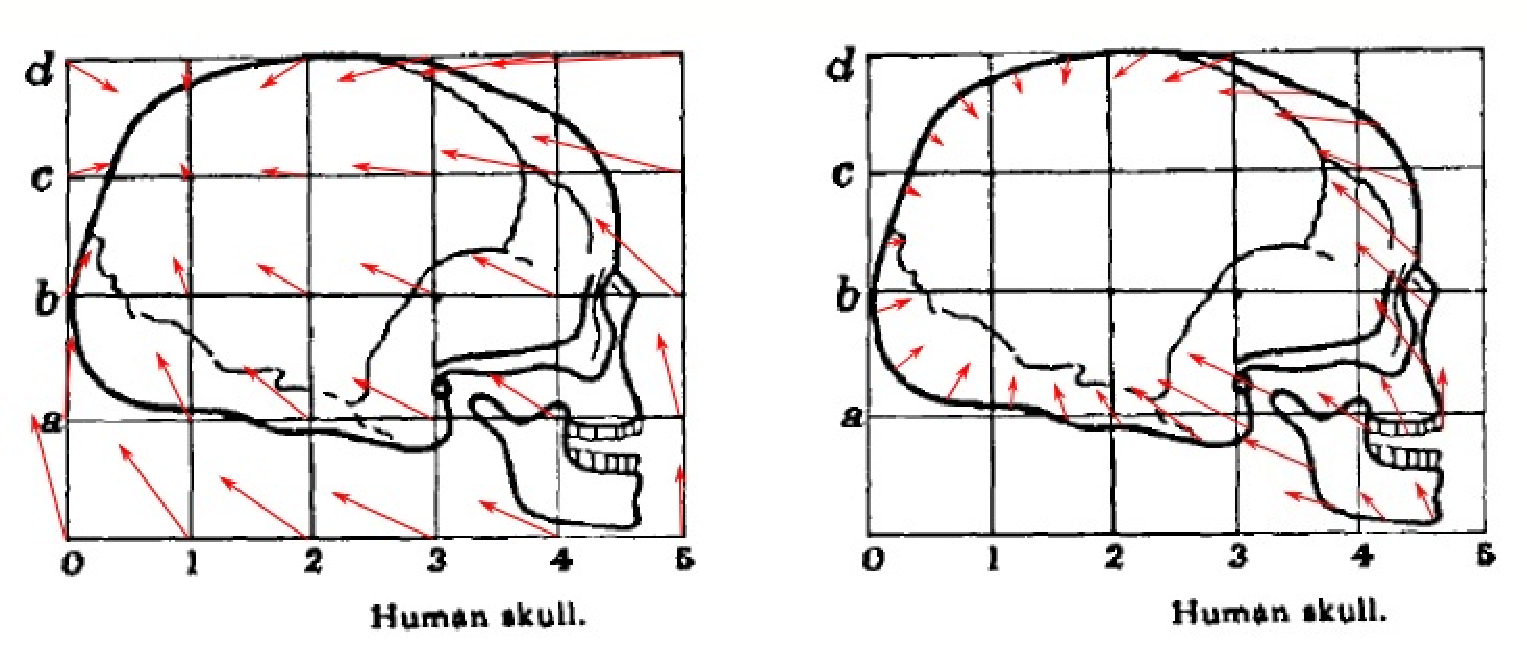}
\caption{Infinitesimal transformation of ambient space (left) as measured
by an outer metric and infinitesimal transformation of the shape itself (right)
as measured by an inner metric.}
\label{in:inout}
\end{figure}

Riemannian metrics measure \emph{infinitesimal deformations}. 
Riemannian metrics on shape space come in two flavors: 
\begin{itemize}
\item Outer metrics measure how much ambient space has to be deformed 
in order to yield the desired deformation of the shape.
An infinitesimal deformation of ambient space is a vector field on ambient 
space and could be pictured as a small arrow attached to every point in ambient space; 
see figure~\ref{in:inout}.\footnote{The 
graphic is an adaptation by the authors
of a graphic in \cite{Thompson1942}.}
\item Inner metrics measure deformations of the shape itself 
within a fixed ambient space.
An infinitesimal deformation of the shape itself is a
vector field along the shape. It could be pictured as a small arrow attached to every point
of the shape; see figure~\ref{in:inout}. 
\end{itemize}

The metrics treated in this work are inner metrics. 

\subsection{Where this paper comes from}\label{introduction:wherethispaper}
In \cite{Michor98}, Michor and Mumford  investigated inner metrics on the 
shape space of planar curves. 
The simplest such metric is the $L^2$-metric given by
\begin{align*}
G^0_f(h,k) &= \int_{S^1} \g\big( h(\th),k(\th) \big)\,|f'(\th)|\,d\th,
\end{align*}
where $f,h,k: S^1\to\R^2$ are  smooth functions. $f$ is the curve representing the shape, and 
$h,k$ are deformation vector fields of $f$. The Euclidean metric on $\R^2$
is denoted by $\g$. We use integration by arc length $ds=|f'(\th)|\,d\th$. 
This makes the metric invariant under reparametrizations of $f,h,k$. 
This invariance is needed
when factoring out reparametrizations; see section~\ref{hamiltonian}.
Since the metric has to be positive definite, it is natural to require 
that $f'(\th)\neq 0$ everywhere, i.e., $f$ is an immersion. 

Unfortunately it turned out that the $L^2$ metric induces vanishing geodesic distance 
on shape space \cite{Michor98}.
This means that any two shapes can be connected by an arbitrarily short path in shape space, 
when path length is measured with the $L^2$ metric.
Later in \cite{Michor102} it was found that the vanishing geodesic distance phenomenon for the 
$L^2$-metric occurs also in 
the more general shape space where $S^1$ is replaced by a compact 
manifold $M$ and Euclidean $\mathbb R^2$ is replaced by a
Riemannian manifold $N$. It also occurs on the full diffeomorphism group 
$\on{Diff}(N)$.\footnote{But not on 
the subgroup $\on{Diff}(N,\on{vol})$ of volume preserving diffeomorphisms, where the geodesic 
equation for the $L^2$-metric is the Euler equation of an incompressible fluid.} 
These results together imply 
vanishing geodesic distance on spaces of immersions.\footnote{This has not been stated 
in \cite{Michor98,Michor102}, but it follows easily. First choose a short horizontal path going from 
the immersion $f_0$ to the $\Diff(M)$-orbit of the immersion $f_1$. 
Then choose another short path in the $\Diff(M)$-orbit 
of $f_1$ connecting the endpoint of the previous path to $f_1$.}

The discovery of the degeneracy of the $L^2$ metric was the starting point of 
a quest for better Riemannian metrics. 
A possibility excluding the degenerate paths encountered in \cite{Michor98} 
is to penalize high length and/or curvature. This led to the
investigation of a class of metrics which were called almost local metrics; 
see  \cite{Michor107,Michor98}. A better name might be weighted $L^2$-metrics. 
They are of the form 
$$
G^\Ph_f(h,k) = \int_{S^1}\Ph(\ell(f),\ka_f(\th)) \g\big( h(\th),k(\th) \big)|f'(\th)|\,d\th , 
$$
where $\Ph:\mathbb R^2\to \mathbb R_{>0}$ is a suitable smooth function, 
$\ell(f)=\int_{S^1}|f'(\th)|\,d\th$ is the length of $f$, and $\ka_f$ is the curvature 
of $f$. If $\Ph=\Ph(\ell(f))$, then this is just a conformal change of the metric; 
it was proposed and investigated independently in \cite{Shah2008} and in 
\cite{YezziMennucci2004,YezziMennucci2004a,YezziMennucci2005}.
For $\Ph=1+A\ka_f^2$, where $A$ is a positive constant,
the metric was investigated in great detail in \cite{Michor98}.

\subsection{Almost local metrics, geodesics, and curvature}
In this paper we take up the investigation of almost local metrics from \cite{Michor107} 
and we generalize it to higher dimensions. 
For surfaces in $\R^3$ this leads to metrics of the form
$$ G^\Ph_f(h,k) = \int_{M} \Ph(\on{Area}(f),\mu)\,\g( h, k)\,d\on{Area}.$$
Here $M$ is a compact connected oriented two-dimensional manifold, $\g$ is the 
Euclidean metric on $\R^3$, $f:M\to \R^3$ is an immersion, 
$h,k: M\to \R^3$ are seen as deformation vector fields of the immersion, 
and $\Ph:\mathbb R^3\to \mathbb R_{>0}$ is again a suitable positive smooth function
depending on the area of the immersed surface $f(M)$ and on the mean curvature $\mu$. 
In two dimensions, the mean curvature $\mu=\Tr(L)$ and Gau{\ss} curvature $\ka=\det(L)$ 
are all the invariants of the Weingarten mapping $L$. 
However, in this paper we do not treat the metric which involves the Gau{\ss} curvature.
This is done in the paper \cite{Michor120}. 

We do not restrict ourselves to surfaces in $\R^3$.
Instead we treat almost local metrics on spaces of hyper-surfaces in $\R^n$. 
More specifically, these are metrics of the form
$$ G^\Ph_f(h,k) = \int_{M} \Ph(\on{Vol}(f),\on{Tr}(L))\,\g( h, k)\on{vol}(g).$$
Here $M$ is a compact connected oriented $n-1$ dimensional manifold, $f:M\to \R^n$ is an immersion, 
$h,k: M\to \R^n$ are  deformation vector fields, and
$\Ph:\mathbb R^2\to \mathbb R_{>0}$ is again a suitable positive smooth function.
We denote the pullback of the Euclidean metric $\g$ to $M$ via the immersion $f$ by $g=f^*\g$.
It is sometimes called the first fundamental form of the surface and is given in 
coordinates by $g_{ij}=\g(\p_i f,\p_j f)$.
The natural replacement of $dA$ is the $n-1$ dimensional volume density $\vol(g)$
induced by $g$. In coordinates $u^1,\ldots u^{n-1}$ on $M$ 
it is given by $\sqrt{\det(g_{ij})}\,|du^1\wedge \cdots\wedge du^{n-1}|$.
The total $n-1$ dimensional volume of $f(M)$ is denoted by $\on{Vol}(f) = \int_M\on{vol}(g)$. Finally,  
$\on{Tr}(L)$ is the mean curvature, which is the trace of the Weingarten mapping $L$.
See section~\ref{no}
for a rigorous definition of the objects that are used in the definition of the metric.

As mentioned above, metrics of this form will be called almost local as in \cite{Michor107}. 
A better name might be weighted $H^0$-metrics or weighted $L^2$-metrics. 
It is natural to consider Gau\ss{}-curvature 
weighted metrics as well. This is done in \cite{Michor120}. It might also be worth 
considering other curvature invariants.

In the theoretical part of this work $\Ph$ is general. The special choices of $\Ph$ 
investigated in the more practical two last sections are
$$\Ph=\Vol^k,\quad \Ph=e^{\Vol},\quad \Ph=1+A\Tr(L)^{2k}, \quad \Ph=\Vol^{\frac{1+n}{1-n}}+A \frac{\Tr(L)^2}{\Vol},$$
where $A>0$ and $k \in \mathbb N$ are constants. The last choice of $\Ph$ induces a scale-invariant metric. 

The $G^\Ph$-metrics are weak Riemannian metrics on the manifold 
of all immersions $M\to \mathbb R^n$, which is 
an open subset of the Fr\'echet space of all smooth mappings. 
But we are interested in the induced Riemannian metric on shape space, which is the quotient space 
under the identification of immersions  differing  only by a reparametrization; see \ref{sh:sub}. 
Shape space is difficult to
handle, but geodesics on shape space are images of so-called \emph{horizontal geodesics} 
on the space of immersions.
A geodesic on the manifold of immersions is called horizontal 
if its velocity vector is a horizontal tangent vector at each time. 
A tangent vector is called horizontal if it is $G^\Ph$-perpendicular to the reparametrization orbits.
The length of a horizontal (minimizing) geodesic defines the distance between its endpoints, 
which is what we are interested in. 

In general, geodesics are critical points of the \emph{energy} functional  
$$E(f)=\frac12\int_0^1 G^\Ph_{f(t)}(\p_t f,\p_t f)\,dt,$$ where 
$f$ is a smooth curve of immersions.
A curve $f(t)$ is a critical point of the 
\emph{horizontal energy} functional 
$$E^{\text{hor}}(f)=\frac12\int_0^1 G^\Ph_{f(t)}\big((\p_t f)^{\hor},(\p_t f)^{\hor}\big)\,dt$$
if and only if a suitable reparametrization of it is a horizontal 
geodesic. In the above formula $(\p_t f)^{\hor}$ 
denotes the horizontal part of the velocity $\p_t f$; see section \ref{sh:sub}.

Almost local metrics have the great advantage that
a tangent vector $h:M\to \mathbb R^n$ with footpoint an immersion $f:M\to \mathbb R^n$ 
is horizontal if and only if $h(x)$ is normal to 
$T_{f(x)}f(M)$ in $\mathbb R^n$ for all $x\in M$; see section \ref{geodesic_equation_Bi:horizontal}. 
Therefore the horizontal energy for almost local metrics is given by an easy and computable formula.
This makes the numerical approach in this paper possible; see section \ref{numerics}. 
But an analytical proof of the existence of critical points for the horizontal energy
(which can be viewed as an anisotropic plateau problem) is still lacking. 
The simple form of the horizontal bundle also opens  the way to computations of \emph{sectional 
curvature} on shape space; see section \ref{sectional_curvature}. We are interested in sectional curvature
because it will eventually be important for doing statistics on shape space and  for the computation of
conjugate points. Furthermore,
unbounded positive sectional curvature 
might be related to vanishing geodesic distance (see  \cite{Michor102}).

\subsection{Contributions of this work}
Almost local metrics are generalized to higher dimensions. 
Some estimates for geodesic distance on shape space are 
proven. They show that almost local metrics with suitable functions $\Ph$ 
overcome the degeneracy of the $L^2$ metric.
In addition, almost local metrics are compared to the Fr\'echet metric.

The geodesic equation and conserved quantities are calculated on shape space
and on the full space of immersions.
For this aim the Hamiltonian formalism developed in 
\cite{Michor107} is updated to the 
more general situation here. 
Furthermore, the Riemann curvature tensor is calculated on shape space.
It contains some negative, positive, and indefinite terms.
Explicit formulas for special choices of $\Ph$ are given.

For all these calculations, first and second derivatives of the metric, volume form, second fundamental form 
and other curvature terms are developed. The derivatives are taken
with respect to the immersion inducing these objects. 

The last section contains numerical experiments for geodesics. 
We do only  boundary value problems and no 
initial value problems (it is not clear if the initial value problem is well posed).
We use Mathematica to set up the triangulation of the surfaces, feed this into AMPL (a modeling 
software developed for optimization), and use the solver IPOPT. 
The numerical results are tested on the totally geodesic subspace of concentric spheres 
where we also have analytic solutions. 
Then we study translations and deformations of surfaces for various metrics and discuss the 
appearing phenomena.

For the sake of simplicity we have restricted ourselves to the shape space of hyper-surfaces in $\R^n$. 
The more general case of arbitrary co-dimension and $\R^n$ replaced by a non-flat 
Riemannian manifold $(N,\g)$ will be treated in another paper. 

\medskip
\section{Shape space and the Hamiltonian approach}\label{hamiltonian}{\hfill}\par
\smallskip
The aim of this chapter is to develop a rigorous notion of shape space, to derive 
the geodesic equation on shape space, and to calculate the conserved momenta. 

\subsection{Manifolds of immersions and embeddings and the diffeomorphism group}\label{sh:im}
Mathematically, parametrized surfaces will be modeled as immersions or embeddings 
of one manifold into another.
We call immersions and embeddings parametrized since a change in their parametrization 
(i.e., applying a diffeomorphism on the domain of the function) 
results in a different object. 
We will deal with the following sets of functions:
\begin{equation}\label{sh:im:eq1}
\Emb(M,\R^n) \subset \Imm(M,\R^n) \subset C^\infty(M,\R^n).
\end{equation}
$C^\infty(M,\R^n)$ is the set of smooth functions from $M$ to $\R^n$. 
$\Imm(M,\R^n)$ is the set of all \emph{immersions} of $M$ into $\R^n$, i.e.,
all functions $f \in C^\infty(M,\R^n)$ such that 
$T_x f$ is injective for all $x \in M$. 
$\Emb(M,\R^n)$ is the set of all \emph{embeddings} of $M$ into $\R^n$, i.e.,
all immersions $f$ that are a homeomorphism onto their image. 
In most cases, immersions will be used since this is the most general 
setting. Working with embeddings instead of immersions makes a difference in 
section~\ref{geodesic_distance}. 

Since $M$ is compact, by assumption it follows that $C^\infty(M,\R^n)$ is a 
\emph{Fr\'echet manifold} \cite[section~42.3]{MichorG}. 
All inclusions in \eqref{sh:im:eq1} are inclusions of open subsets.
Therefore, all function spaces in \eqref{sh:im:eq1} are Fr\'echet manifolds as well. 

The tangent bundle of the manifold of immersions is
$$T\on{Imm}(M,\mathbb R^n)=\on{Imm}(M,\R^n)\x C^\infty(M,\R^n) \subset C^\infty(M,\R^n \x \R^n),$$
and the cotangent bundle is
$T^*\on{Imm}(M,\mathbb R^n)=\on{Imm}(M,\R^n)\x \mathcal D'(M)^n$,
where the second factor consists of $n$-tuples of distributions in
$\mathcal D'(M)=C^\infty(M)'$, which is the space of distributional sections of
the density bundle.

By $\on{Diff}(M)$ we will denote the group of all 
smooth diffeomorphisms.  $\on{Diff}(M)$ is a  Fr\'echet manifold as well, 
since it is an open subset of  $C^\infty(M,M)$. It is an infinite dimensional 
Lie group  in the sense of \cite[section~43]{MichorG}.
The diffeomorphism group $\Diff(M)$ acts smoothly on $C^\infty(M,\R^n)$ and its subspaces
$\Imm(M,\R^n)$ and $\Emb(M,\R^n)$ by composition from the right. The action is given by the mapping
$$\Imm(M,\R^n) \x \Diff(M) \to \Imm(M,\R^n), \qquad (f,\ph) \mapsto r(f,\ph) = r^\ph(f)= f \o \ph.$$
The tangent prolongation of this group action is given by the mapping
\begin{align*}
T(r^\ph): T\Imm(M,\R^n) \x \Diff(M) &\to T\Imm(M,\R^n), \\
(f,h,\ph) &\mapsto (f\o\ph,h \o \ph).
\end{align*}

We will sometimes use the abbreviations $\Emb, \Imm$, and $\Diff$
when the domain and co-domain of the functions are clear from the context.

\subsection{Riemannian metrics on the manifold of immersions}\label{hamiltonian:setting}
In this work we consider smooth Riemannian metrics on $\on{Imm}(M,\mathbb R^n)$, i.e.,
smooth mappings
\begin{align*}
&G:\on{Imm}(M,\mathbb R^n)\x C^\infty(M,\mathbb R^n)
  \x C^\infty(M,\mathbb R^n)\to \mathbb R, \\&
(f,h,k)\mapsto G_f(h,k),\quad \text{ bilinear in }h,k,
\\& G_f(h,h)>0 \quad\text{  for }h\ne0.
\end{align*}
Each such metric is {\it weak} in the sense that $G_f$, viewed as bounded
linear mapping
\begin{align*}
G_f:T_f\on{Imm}(M,\mathbb R^n)=C^\infty(M,\mathbb R^n)&\to
T_f^*\on{Imm}(M,\mathbb R^n)=\mathcal D'(M)^n,
\\
G: T\on{Imm}(M,\mathbb R^n)&\to T^*\on{Imm}(M,\mathbb R^n),
\\
G(f,h)&=(f,G_f(h,\;.\;)),
\end{align*}
is injective but can never be surjective.
We shall need also its tangent
mapping
\begin{align*}
TG:T(T\on{Imm}(M,\mathbb R^n))&\to T(T^*\on{Imm}(M,\mathbb R^n)).
\end{align*}
We write a tangent vector to $T\on{Imm}(M,\mathbb R^n)$ as
$(f,h;k,v)$, where $(f,h)\in T\on{Imm}(M,\mathbb R^n)$ is its footpoint,
$k$ is its vector component in the $\on{Imm}(M,\mathbb R^n)$-direction,
and where $v$ is
its component in the $C^\infty(M,\mathbb R^n)$-direction.

Then $TG$ is given by
$$
TG(f,h;k,v) = (f, G_f(h,\;.\;);k, D_{(f,k)}G_f(h,\;.\;)+G_f(v,\;.\;)).
$$

Note that only these smooth functions on $\on{Imm}(M,\mathbb R^n)$ whose
derivatives lie in the image of $G$ in the cotangent bundle have
$G$-gradients.
This requirement has only to be satisfied for the first
derivative; for the higher ones it follows (see \cite{MichorG}).
We shall denote
by $C^\infty_G(\on{Imm}(M,\mathbb R^n))$ the space of such smooth
functions.

In what follows we shall further assume that that {\it the weak Riemannian
metric $G$ itself admits $G$-gradients with respect to the variable $f$ in
the following sense:}
\begin{align*}
&\boxed{
D_{(f,m)}G_f(h,k) = G_f(m,H_f(h,k)) = G_f(K_f(m,h),k),
}\quad\text{  where }
\\
&H,K:\on{Imm}(M,\mathbb R^n)\x C^\infty(M,\mathbb R^n)
  \x C^\infty(M,\mathbb R^n)\to C^\infty(M,\mathbb R^n)
\\& \hspace*{2.5in}(f,h,k)\mapsto H_f(h,k),K_f(h,k)
\\&
\text{are smooth and bilinear in }h,k.
\end{align*}
Note that $H$ and $K$ could be expressed in abstract index notation as
$g_{ij,k}g^{kl}$ and $g_{ij,k}g^{il}$.
We will check and compute these gradients for several concrete metrics below.

\subsection{The fundamental symplectic form on $T\on{Imm}(M,\mathbb R^n)$
induced by a weak Riemannian metric}\label{hamiltonian:symplectic_form}
The basis of Hamiltonian theory is the natural 1-form on the cotangent
bundle $T^*\on{Imm}(M,\mathbb R^n)$
given by
\begin{gather*}
\Th: T(T^*\on{Imm}(M,\mathbb R^n))=\on{Imm}(M,\mathbb R^n)\x
\mathcal D'(M)^n\x C^\infty(M,\mathbb R^n)\x \mathcal D'(M)^n \to \mathbb R,
\\
(f,\al;h,\be) \mapsto \langle  \al,h\rangle.
\end{gather*}
The pullback via the mapping $G: T\on{Imm}(M,\mathbb R^n)\to
T^*\on{Imm}(M,\mathbb R^n)$ of $\Th$ is
$$ (G^*\Th)_{(f,h)}(f,h;k,v)=G_f(h,k).$$
Thus the symplectic form $\om=-dG^*\Th$ on $ T\on{Imm}(M,\mathbb R^n)$
can be computed as follows, where we use the constant vector fields
$(f,h)\mapsto (f,h;k,v)$:
\begin{align*}
\om_{(f,h)}&((k_1,v_1),(k_2,v_2)) =
-d(G^*\Th)((k_1,v_1),(k_2,v_2))|_{(f,h)} \\&
=-D_{(f,k_1)}G_f(h,k_2)- G_f(v_1,k_2)
+D_{(f,k_2)}G_f(h,k_1)+ G_f(v_2,k_1)
\\&
=G_f\big(k_2,H_f(h,k_1)-K_f(k_1,h)\big) + G_f(v_2,k_1) - G_f(v_1,k_2).
\tag1
\end{align*}

\subsection{The Hamiltonian vector field mapping}\label{hamiltonian:hamiltonian_vector field}
Here we compute the Hamiltonian vector field $\on{grad}^\om(F)$ associated
to a smooth function $F$ on the tangent space $T\on{Imm}(M,\mathbb R^n)$;
that is $F \in C^\infty_G(\on{Imm}(M,\mathbb R^n)\x C^\infty(M,\mathbb R^n))$
assuming that it has smooth $G$-gradients in both factors. See \cite[section~48]{MichorG}.
Using the explicit formulas in section \ref{hamiltonian:symplectic_form}, we have
\begin{align*}
&\om_{(f,h)}\left(\on{grad}^\om(F)(f,h),(k,v) \right)\\
&\quad=\om_{(f,h)}\left((\on{grad}_1^\om(F)(f,h),\on{grad}_2^\om(F)(f,h)),(k,v) \right) 
\\&\quad
= G_f\big(k, H_f\big(h,\on{grad}_1^\om(F)(f,h)\big)\big)
-G_f(K_f(\on{grad}_1^\om(F)(f,h),h),k)
\\&\qquad
 + G_f(v,\on{grad}_1^\om(F)(f,h))
 - G_f(\on{grad}_2^\om(F)(f,h),k).
\end{align*}
On the other hand, by the definition of the $\om$-gradient we have
\begin{align*}
&\om_{(f,h)}\left(\on{grad}^\om(F)(f,h),(k,v) \right)
=dF(f,h)(k,v)= D_{(f,k)}F(f,h) + D_{(h,v)}F(f,h)
\\&
=G_f(\on{grad}_1^{G}(F)(f,h),k)
+G_f(\on{grad}_2^{G}(F)(f,h),v),
\end{align*}
and we get the expression of the Hamiltonian vector field:
$$
\boxed{
\begin{aligned}
\on{grad}_1^\om(F)(f,h) &= \on{grad}_2^{G}(F)(f,h),\\
\on{grad}_2^\om(F)(f,h) &= \!-\! \on{grad}_1^{G}(F)(f,h)
\\&\quad
+ H_f\big(h,\on{grad}_2^{G}(F)(f,h)\big)
- K_f(\on{grad}_2^G(F)(f,h),h).
\end{aligned}
}
$$
Note that for a smooth function $F$ on $T\on{Imm}(M,\mathbb R^n)$ the
$\om$-gradient exists if and only if both $G$-gradients exist.

\subsection{The geodesic equation on the manifold of immersions} \label{hamiltonian:geodesic_equation}

The geodesic flow is defined by a vector field on $T\on{Imm}(M,\mathbb
R^n)$.
One way to define this vector field is as the Hamiltonian vector
field of the energy function
$$ E(f,h)=\frac12 G_f(h,h),\qquad
E:\on{Imm}(M,\mathbb R^n)\x C^\infty(M,\mathbb R^n)\to \mathbb R .$$
The two partial $G$-gradients are
\begin{align*}
G_f(\on{grad}_2^{G}(E)(f,h),v) &= d_2E(f,h)(v) = G_f(h,v),
\\
\on{grad}_2^{G}(E)(f,h)&=h,
\\
G_f(\on{grad}_1^{G}(E)(f,h),k) &= d_1E(f,h)(k)
= \tfrac12 D_{(f,k)}G_f(h,h)
\\&
=\tfrac12 G_f(k,H_f(h,h)),
\\
\on{grad}_1^{G}(E)(f,h) &= \tfrac12 H_f(h,h).
\end{align*}
Thus the geodesic vector field is
\begin{align*}
\on{grad}_1^\om(E)(f,h) &= h
\\
\on{grad}_2^\om(E)(f,h) &=  \tfrac12 H_f(h,h) - K_f(h,h)
\end{align*}
and the geodesic equation becomes
\begin{align*}
&\begin{cases}
f_t &= h,
\\
h_t &=  \tfrac12 H_f(h,h) - K_f(h,h),
\end{cases}
\quad\text{  or }\quad
\boxed{
f_{tt}= \tfrac12 H_f(f_t,f_t) - K_f(f_t,f_t).
}
\end{align*}
This is nothing but the usual formula for the geodesic flow using the
Christoffel symbols expanded out using the first derivatives of the metric
tensor.

\subsection{The momentum mapping for a $G$-isometric group action}\label{hamiltonian:momentum_mappings}
We consider now a (possibly infinite dimensional regular) Lie group with
Lie algebra $\mathfrak g$ with a right action $g\mapsto r^g$ by isometries
on $\on{Imm}(M,\mathbb R^n)$.
Denote by $\X(\on{Imm}(M,\mathbb R^n))$ 
the set of vector fields on $\on{Imm}(M,\mathbb R^n)$. Then we can specify
this action by the fundamental vector field mapping $\ze:\mathfrak g\to
\X(\on{Imm}(M,\mathbb R^n))$, which will be a bounded Lie algebra
homomorphism.
The fundamental vector field $\ze_X, X \in \mathfrak g$, is
the infinitesimal action in the sense that
$$\ze_X(f)=\p_t|_0 r^{\exp(tX)}(f).$$
We also consider the tangent prolongation of this action on
$T\on{Imm}(M,\mathbb R^n)$, where the fundamental vector field is given by
\begin{equation*}
\ze_X^{T\on{Imm}}:(f,h)\mapsto (f,h;\ze_X(f),D_{(f,h)}(\ze_X)(f)=:\ze'_X(f,h)).
\end{equation*}
The basic assumption is that the action is by isometries,
$$ G_f(h,k)=((r^g)^*G)_f(h,k)=G_{r^g(f)}(T_f(r^g)h,T_f(r^g)k).
$$
Differentiating this equation at $g=e$ in the direction $X\in \mathfrak g$,
we get
\begin{equation*}
0=D_{(f,\ze_X(f))}G_f(h,k)+G_f(\ze'_X(f,h),k)+G_f(h,\ze'_X(f,k)).
\tag{1}\end{equation*}
The key to the Hamiltonian approach is to define the group action by
Hamiltonian flows.
We define the {\it momentum map}
$j:\mathfrak g\to C^\infty_G (T\on{Imm} (M,\mathbb R^n),\mathbb R)$ by
$$\boxed{j_X(f,h) = G_f(\ze_X(f),h).}$$
Equivalently, since this map is linear, it is often written as a map
$$\mathcal J:T\on{Imm}(M,\mathbb R^n) \rightarrow \mathfrak g',\qquad
  \langle \mathcal J(f,h),X \rangle = j_X(f,h).$$
The main property of the momentum map is that it fits into the following
commuta\-ti\-ve diagram and is a homomorphism of Lie algebras:
\begin{displaymath}
\cgaps{0.8;0.7;0.7;0.8}\xymatrix{
H^0(T\on{Imm}) \ar[r]^{i} &
     C^\infty_G(T\on{Imm},\mathbb R)  \ar[rr]^{\on{grad}^{\om}} & &
     \X(T\on{Imm},\om) \ar[r] & H^1(T\on{Imm})
\\
 & & \mathfrak g \ar[lu]^{j} \ar[ru]_{\ze^{T\on{Imm}}} & &
}\end{displaymath}
where $\X(T\on{Imm},\om)$ is the space of vector fields on $T\on{Imm}$
whose flow leaves $\om$ fixed.
Note also that $\mathcal J$ is equivariant for the group action.
See \cite{Michor107} for more details.

By Noether's theorem, along any geodesic $t\mapsto f(t,\;.\;)$ this
momentum mapping is constant; thus for any $X\in\mathfrak g$ we have
\begin{displaymath}
\boxed{\quad
\langle \mathcal J(f,f_t),X \rangle = j_X(f,f_t) = G_f(\ze_X(f),f_t) \quad\text{  is
constant in }t.
\quad}
\end{displaymath}

We can apply this construction to the following group actions:
\begin{itemize}
\item
The smooth right action of the group
$\on{Diff}(M)$ on $\on{Imm}(M,\mathbb R^n)$,
given by composition from the right:
$f\mapsto f\o \ph$ for $\ph\in\on{Diff}(M)$.

For $X\in\X(M)$ the fundamental vector field is then given by
$$ \ze^{\on{Diff}}_X(f) = \ze_X(f) = X(f) = \p_t|_0 (f\o \Fl^X_t) = df.X, $$
where $\Fl^X_t$ denotes the flow of $X$.
The {\it reparametrization momentum}, for any vector field $X$ on $M$ is thus
$$ j_X(f,h) = G_f(df.X, h).$$
Assuming the metric is reparametrization invariant, it follows that on any
geodesic $f(x,t)$, the expression $G_f(df.X,f_t)$ is constant for all $X$.
\item
The left action of the Euclidean motion group
$\mathbb R^n\rtimes SO(n)$ on $\on{Imm}(M,\mathbb R^n)$ given by
$f\mapsto Af+ B$ for $(B,A)\in \mathbb R^n\x SO(n)$.
The fundamental vector field mapping is
\begin{align*}
\ze_{(B,X)}(f)&= Xf+B.
\end{align*}
The {\it linear-momentum} is thus $G_f(B,h), B \in \mathbb R^n$, and if the
metric is trans\-la\-tion invariant, $G_f(B,f_t)$ will be constant along
geodesics for every $B\in \mathbb R^n$.
The {\it angular-momentum} is similarly
$G_f(X.f,h), X\in \mathfrak s\mathfrak o(n)$, and if the
metric is rotation-invariant, then $G_f(X.f,f_t)$ will be constant along
geodesics for each $X\in \mathfrak s\mathfrak o(n)$.

\item
The action of the scaling group of $\mathbb R$ given by $c\mapsto e^r f$,
with fundamental vector field $\ze_a(f)=a.f$.

If the metric is scale-invariant, then
the {\it scaling momentum} $G_f(f,f_t)$ will also be invariant along
geodesics.
\end{itemize}

\subsection{Shape space}
As discussed in the introduction, by a shape we mean a
smoothly immersed or embedded  hypersurface in $\R^n$ which is diffeomorphic to a 
fixed compact, connected, and oriented
manifold $M$ of dimension $n-1$. 
The space of these shapes will be  denoted $B_i(M,\mathbb R^n)$ or $B_e(M,\mathbb R^n)$ and viewed as the 
quotient 
$$ B_e(M,\mathbb R^n) = \on{Emb}(M,\mathbb R^n)/\on{Diff}(M)\text{ or }  
B_i(M,\mathbb R^n) = \on{Imm}(M,\mathbb R^n)/\on{Diff}(M).$$ 
In \cite[section~44.1]{MichorG} it is shown that $B_e(M,\R^n)$ is a manifold again. 
$B_i(M,\mathbb R^n)$  is, however, no longer a manifold
but an orbifold with finite isotropy groups; see \cite{Michor40}.
We will sometimes use the abbreviations $B_i$ and $B_e$ 
when it is clear what the domain and co-domain of the functions are.

More generally, a shape will be an element of the Cauchy completion (i.e., the metric completion 
for the geodesic distance) of $B_i(M,\mathbb R^n)$ with 
respect to a suitably chosen Riemannian metric. This will allow for corners.

\subsection{Riemannian submersions and the metric on shape space}\label{sh:sub}

We will always assume that a $\Diff(M)$-invariant Riemannian metric on $\Imm(M,\R^n)$ is given. 
Then there is a unique Riemannian metric on the quotient space $B_i(M,\R^n)$ such that the
quotient map $\pi:\Imm(M,\R^n)\rightarrow B_i(M,\R^n)$ is a 
\emph{Riemannian submersion}. This is the construction 
that we use to induce a metric on shape space.

Let $\on{ker}(T\pi) \subset T\Imm(M,\R^n)$
be the \emph{vertical bundle}. The \emph{horizontal bundle} is
its orthogonal complement with respect to the metric $G$.
Then $T_{\pi(f)}B_i(M,\R^n)$ is isometric to the horizontal bundle at $f \in \Imm(M,\R^n)$.
Note that the horizontal bundle 
depends on the definition of the metric. For almost local metrics, it consists of vector 
fields along $f$ that are everywhere normal to $f$; see section~\ref{geodesic_equation_Bi:horizontal} .

By the conservation of the reparametrization momentum, geodesics in the space of immersions
with horizontal initial velocity stay horizontal for all time. Such geodesics 
project down to geodesics in shape space because $\pi$ is a Riemannian submersion.
See \cite[section~26]{MichorH} for a proof of this fact. 
We will show in section~\ref{geodesic_equation_Bi:horizontal} 
that almost local metrics have the property that 
any curve in shape space can be lifted to a horizontal curve of immersions. 
This implies that instead of solving the geodesic equation on shape space one can equivalently solve
the equation for horizontal geodesics in the space of immersions.

\medskip
\section{Differential geometry of surfaces and notation}\label{no}{\hfill}\par
\smallskip
In this section we will present and develop the differential geometric tools 
that are needed to deal with immersed surfaces. 
The most important point is a rigorous treatment of the covariant derivative 
and related concepts. 

In \cite[section 2]{Michor119} one can find some parts of this section in a more general setting.
We use the notation of \cite{MichorH}. Some of the definitions can also be found in \cite{Kobayashi1996a}.

\subsection{Tensor bundles and tensor fields}\label{no:te}

We will deal with the \emph{tensor bundles}
\begin{equation*}\xymatrix{
T^r_s M \ar[d] & 
T^r_s M \otimes f^*T\R^n \ar[d] \\
M & M 
}\end{equation*}
Here $T^r_sM$ denotes the bundle of 
$\left(\begin{smallmatrix}r\\s\end{smallmatrix}\right)$-tensors on $M$, i.e.,
$$T^r_sM=\bigotimes^r TM \otimes \bigotimes^s T^*M,$$
and $f^*T\R^n$ is the pullback of the bundle $T\R^n$ via $f$; see \cite[section~17.5]{MichorH}. 
A \emph{tensor field} is a section of a tensor bundle. Generally, when $E$ is a bundle, 
the space of its sections will be denoted by $\Ga(E)$. 

To clarify the notation that will be used later, 
some examples of tensor bundles and tensor fields are given now.
\begin{itemize}
\item $\on{End}(TM)=L(TM,TM)=T^1_1M$ is the bundle of \emph{endomorphisms of $TM$}. 
\item $S^k T^*M = L^k_{\on{sym}}(TM; \R)$ is the bundle of \emph{symmetric 
$\left(\begin{smallmatrix}0\\k\end{smallmatrix}\right)$-tensors}.
\item $S^2_{>0} T^*M$ is the bundle of \emph{symmetric positive definite
$\left(\begin{smallmatrix}0\\2\end{smallmatrix}\right)$-tensors}.
\item $\La^k T^*M = L^k_{\on{alt}}(TM; \R)$ is the bundle of \emph{alternating
$\left(\begin{smallmatrix}0\\k\end{smallmatrix}\right)$-tensors}.
\item $\Om^r(M)=\Ga(\La^r T^*M)$ is the space of \emph{differential forms}.
\item $\X(M)=\Ga(TM)$ is the space of \emph{vector fields}.
\item $\Ga(f^*T\R^n)\cong \big\{ h \in C^\infty(M,T\R^n): \pi_N \o h = f \big\}$ is the space of 
\emph{vector fields along $f$}.
\end{itemize}

For $X \in \X(M)$ the insertion $\io_X$ will always insert $X$ into the leftmost covariant 
entry of a tensor.

\subsection{Metric on tensor spaces}\label{no:me}

Let $\g \in \Gamma(S^2_{>0} T^*\R^n)$ denote the Euclidean metric on $\R^n$. 
The \emph{metric induced on $M$ by $f \in \Imm(M,\R^n)$} is the pullback metric 
\begin{align*}
g=f^*\g \in \Gamma(S^2_{>0} T^*M), \qquad g(X,Y)=(f^*\g)(X,Y) = \g(Tf.X,Tf.Y),
\end{align*}
where $X,Y$ are vector fields on $M$.
The dependence of $g$ on the immersion $f$ should be kept in mind.
Let $(u,U)$ be a fixed chart 
on $M$ with $\p_i=\frac{\p}{\p u^i}$. In these coordinates the pullback metric is given by
$$g|_U=\sum_{i,j}^{n-1}g_{ij}du^i \otimes
du^j=\sum_{i,j}^{n-1}\g\big(\p_if,\p_j f \big) du^i \otimes du^j.$$
The metric can be seen as a mapping
$$g: TM \to T^*M, \qquad X \mapsto g(X)=:X^\flat$$ 
with inverse
$$g\i: T^*M \to TM, \qquad \al \mapsto g\i(\al)=:\al^\sharp.$$ 
This defines a metric on the cotangent bundle $T^0_1M=T^*M$ via
$$g^0_1(\alpha,\beta)=g\i(\alpha,\beta)=\alpha(\beta^\sharp)=g(\alpha^\sharp,\beta^\sharp)$$
for $\alpha,\beta \in T^*M$. 
The product metric 
$$g^r_s = \bigotimes^r g \otimes \bigotimes^s g\i$$
extends $g$ to all tensor spaces $T^r_s M$, and 
$g^r_s \otimes \g$ yields a metric on $T^r_s M \otimes f^*T\R^n$.

\subsection{Traces}\label{no:tr}

The \emph{trace} contracts pairs of vectors and covectors in a tensor product: 
\begin{align*}
\Tr:\; T^*M \otimes TM = L(TM,TM) \to M \x \R
\end{align*}
A special case of this is the operator
$\io_X$ inserting a vector $X$ into a covector or into a covariant factor of a tensor product.
The inverse of the metric $g$ can be used to define a trace 
$$\Tr^g: T^*M \otimes T^*M \to M \x \R$$
contracting pairs of covectors.
Note that $\Tr^g$ depends on the metric whereas $\Tr$ does not. 
The following lemma will be useful in many calculations (see \cite[section 2]{Michor119}).

\subsection*{Lemma}{\em
$g^0_2(B,C)= \on{Tr}(g\i B g\i C)$ for $B,C \in T^0_2M$ if $B$ or $C$ is symmetric.
}

\noindent
In the expression under the trace, $B$ and $C$ are seen as maps $TM \to T^*M$.

\subsection{Volume density}\label{no:vo}

Let $\Vol(M)$ be the \emph{density bundle} over $M$; see \cite[section~10.2]{MichorH}.
The \emph{volume density} on $M$ induced by $f \in \Imm(M,\R^n)$ is 
$$\vol(g)=\vol(f^*\g) \in \Ga\big(\Vol(M)\big).$$
In a chart $(u,U)$ the volume density reads as
$$\on{vol}(g)=\sqrt{\on{det}(\g\left(\p_if,\p_j f
\right))}\ |du^1\wedge\dots\wedge du^{n-1}|.$$ 
The \emph{volume} of the immersion is given by
$$\Vol(f)=\int_M \vol(f^*\g)=\int_M \vol(g).$$
The integral is welldefined since $M$ is compact. Since $M$ is oriented, we may identify the volume 
density with a differential form.

\subsection{Covariant derivative}\label{no:co}
We will use covariant derivatives on vector bundles as explained in \cite[sections 19.12, 22.9]{MichorH}.
Let $\nabla^g, \nabla^{\g}$ be the \emph{Levi--Civita covariant derivatives} on $(M,g)$
and $(\R^n,\g)$, respectively. 
For any manifold $Q$ and vector field $X$ on $Q$ one has
\begin{align*}
\nabla^g_X:C^\infty(Q,TM) &\to C^\infty(Q,TM), & h &\mapsto \nabla^g_X h, \\
\nabla^{\g}_X: C^\infty(Q,T\R^n) &\to C^\infty(Q,T\R^n), & h &\mapsto \nabla^{\g}_X h.
\end{align*}
Usually we will simply write $\nabla$ for all covariant derivatives.
It should be kept in mind that $\nabla^g$ depends on the metric $g=f^*\g$ and therefore also on 
the immersion $f$. The $\R^n$ covariant derivative $\nabla^{\g}_Xh$ equals the ordinary differential 
$dh(X)$ but remembers the footpoint $f$ of $h$,
i.e. ,  $\nabla^{\g}_X(f,h)=(f,dh(X))$ if we write $(f,h)$ instead of  $h$.
The following properties hold \cite[section 22.9]{MichorH}:
\begin{enumerate}
\item \label{no:co:ba}
$\nabla_X$ respects base points, i.e., 
$\pi \o \nabla_X h = \pi \o h$, where $\pi$ is the projection 
of the tangent space onto the base manifold. 
\item
$\nabla_X h$ is $C^\infty$-linear in $X$. So for a tangent vector $X_x \in T_xQ$, 
$\nabla_{X_x}h$ makes sense and equals $(\nabla_X h)(x)$.
\item
$\nabla_X h$ is $\R$-linear in $h$.
\item
$\nabla_X (a.h) = da(X).h + a.\nabla_X h$ for $a \in C^\infty(Q)$, the derivation property of $\nabla_X$.
\item \label{no:co:eq}
For any manifold $\widetilde Q$ and smooth mapping 
$q:\widetilde Q \to Q$ and $Y_y \in T_y \widetilde Q$ one has
$\nabla_{Tq.Y_y}h=\nabla_{Y_y}(h \o q)$. If $Y \in \X(Q_1)$ and $X \in \X(Q)$ are $q$-related, then 
$\nabla_Y(h \o q) = (\nabla_X h) \o q$.
\end{enumerate}
The two covariant derivatives $\nabla^g_X$ and $\nabla^{\g}_X$ 
can be combined to yield a covariant derivative $\nabla_X$ acting on
$C^\infty(Q,T^r_sM \otimes T\R^n)$ by additionally requiring the following properties 
\cite[section 22.12]{MichorH}:
\begin{enumerate}
\setcounter{enumi}{5} 
\item $\nabla_X$ does not change the grade of tensors, i.e., it induces mappings
\newline
$\nabla_X:C^\infty(Q,T^r_sM \otimes T\R^n) \to C^\infty(Q,T^r_sM \otimes T\R^n)$. 
\item $\nabla_X(h \otimes k) = (\nabla_X h) \otimes k + h \otimes (\nabla_X k)$, 
a derivation with respect to the tensor product.
\item $\nabla_X$ commutes with any kind of contraction (see \cite[section 8.18]{MichorH}). 
A special case of this is
$$\nabla_X\big(\alpha(Y)\big)=(\nabla_X \alpha)(Y)+\alpha(\nabla_X Y) \quad 
\text{for } \alpha\otimes Y :Q \to T^1_1M.$$
\end{enumerate}
Property \ref{no:co:ba} is important because it implies that $\nabla_X$ 
respects spaces of sections of bundles. 
For example, for $Q=M$ and $f \in C^\infty(M,\R^n)$, one gets
$$\nabla_X : \Ga(T^r_s M \otimes f^* T\R^n) \to \Ga(T^r_s M \otimes f^* T\R^n). $$
\subsection{Swapping covariant derivatives}\label{no:sw}

We will make repeated use of some formulas, allowing us to swap covariant derivatives. 
Let $f$ be an immersion, $h$ a vector field along $f$, and $X,Y$ vector fields on $M$. 
Since $\nabla$ is torsion-free, one has \cite[section~22.10]{MichorH}
\begin{equation}\label{no:sw:tor}
\nabla_X Tf.Y-\nabla_Y Tf.X -Tf.[X,Y] = 0.
\end{equation}
Furthermore, one has \cite[section~24.5]{MichorH}
\begin{equation}\label{no:sw:r}
\nabla_X \nabla_Y h - \nabla_Y \nabla_X h - \nabla_{[X,Y]} h 
=0,
\end{equation}
since $\R^n$ is flat.
These formulas also hold when $f:\R \x M \to \R^n$ is a path of immersions, 
$h:\R \x M \to T\R^n$ is a vector field along $f$, and
the vector fields are vector fields on $\R \x M$. 
A case of special importance is when one of the vector fields is $(\p_t,0_M)$ and the 
other $(0_\R,Y)$, where $Y$ is a vector field on $M$. 
Since the Lie bracket of these vector fields vanishes, 
\eqref{no:sw:tor} and \eqref{no:sw:r} yield
\begin{equation}\label{no:sw:tor_with_t}
\nabla_{(\p_t,0_M)} Tf.(0_{\R},Y)-\nabla_{(0_{\R},Y)} Tf.{(\p_t,0_M)} = 0
\end{equation}
and
\begin{equation}\label{no:sw:r_with_t}
\nabla_{(\p_t,0_M)} \nabla_{(0_\R,Y)} h - \nabla_{(0_\R,Y)} \nabla_{(\p_t,0_M)} h
\\= 0 .
\end{equation}
If the context is clear, we shall write $\p_t$ instead of the more detailed notation $(\p_t,0_M)$
and $Y$ instead of $(0_\R,Y)$.

\subsection{Higher covariant derivatives and the Laplace operator}\label{no:co2}

When the covariant derivative is seen as a mapping
$$\nabla: \Gamma(T^r_s M) \to \Gamma(T^r_{s+1}M)\quad \text{or} \quad
\nabla : \Gamma(T^r_sM \otimes f^*T\R^n) \to \Gamma(T^r_{s+1}M \otimes f^*T\R^n),$$
then the \emph{second covariant derivative} is simply $\nabla\nabla=\nabla^2$.
Since the covariant derivative commutes with contractions,
$\nabla^2$ can be expressed as
$$\nabla^2_{X,Y} :=\iota_Y \iota_X \nabla^2 =
\iota_Y \nabla_X \nabla =
\nabla_X\nabla_Y -\nabla_{\nabla_XY} \qquad \text{for $X,Y\in \X(M)$.}$$
Higher covariant derivates are defined as $\nabla^k$, $k \geq 0$. 
We can use the second covariant derivative to define the \emph{Laplace--Bochner operator}. 
It can act on all tensor fields $B$ and is defined as
$$\Delta B =  - \on{Tr}^g(\nabla^2 B).$$ 
For $h=(h^1,\ldots,h^n):M \to \R^n$ one has $\De h = (\De h^1,\ldots, \De h^n)$.
\subsection{Normal bundle}\label{no:no}

The \emph{normal bundle} $\Nor(f)$ of an immersion $f$ is a subbundle of $f^*T\R^n$ 
whose fibers consist of all vectors that are orthogonal to the image of $f$, i.e.,
$$\Nor(f)_x = \big\{ Y \in T_{f(x)}\R^n : \forall X \in T_xM : \g(Y,Tf.X)=0  \big\}.$$
Any vector field $h$ along $f$ can be decomposed uniquely 
into parts {\it tangential} and {\it normal} to $f$ as
$$h=Tf.h^\top + h^\bot,$$ 
where $h^\top$ is a vector field on $M$ and $h^\bot$ is a section of the normal bundle $\Nor(f)$. 
When $f$ is orientable, then
the unit normal field $\nu$ of $f$ can be defined. It is a section of the normal bundle 
with constant $\g$-length one which is chosen such that
$$\big(\nu(x), T_xf.X_1,T_xf.X_2,\dots, T_xf.X_{n-1}\big)$$ 
is a positive oriented basis in $T_{f(x)}\R^n$ if $X_1, \ldots, X_{n-1}$ is a positive oriented basis in $T_xM$.
In this notation the decomposition of a vector field $h$ along $f$ reads as
$$h=Tf.h^\top+a.\nu.$$
The two parts are defined by the relations
\begin{align*}
&a =\g( h,\nu )\in C^\infty(M), \\
&h^\top \in \X(M), \text{ such that } g(h^\top,X)=\g(h,Tf.X) \text{ for  all } X\in \X(M).
\end{align*}

\subsection{Second fundamental form and Weingarten mapping}\label{no:we}

Let $X$ and $Y$ be vector fields on $M$. 
Then the covariant derivative $\nabla_X Tf.Y$ splits into tangential and normal parts as
$$\nabla_X Tf.Y=Tf.(\nabla_X Tf.Y)^\top + (\nabla_X Tf.Y)^\bot = Tf.\nabla_X Y + S(X,Y).$$
$S=S^f$ is the \emph{second fundamental form of $f$}. 
It is a symmetric bilinear form with values in the normal bundle of $f$. 
When $Tf$ is seen as a section of $T^*M \otimes f^*T\R^n$, one has $S=\nabla Tf$ since
$$S(X,Y) = \nabla_X Tf.Y - Tf.\nabla_X Y = (\nabla Tf)(X,Y).$$
Taking the trace of $S$ yields the \emph{vector valued mean curvature} $$\Tr^g(S) \in \Ga\big(\Nor(f)\big).$$
One can define the \emph{scalar second fundamental form $s=s^f$} as
$$s(X,Y) = \g\big(S(X,Y),\nu\big).$$
Moreover, there is the \emph{Weingarten mapping} or \emph{shape operator} $L=L^f=g\i s$.
It is a $g$-symmetric bundle mapping defined by $$s(X,Y) =  g( LX,Y ).$$ 
The eigenvalues of $L$ are called \emph{principal curvatures} and
the eigenvectors \emph{principal curvature directions}. 
$\on{Tr}(L)=\Tr^g(s)$ is the \emph{scalar mean curvature} and for surfaces in $\R^3$  
the \emph{Gau{\ss} curvature} is given by $\det(L)$. The covariant derivative $\nabla_X \nu$ of the normal vector is related to $L$ 
by the \emph{Weingarten equation}
$$ \nabla_X \nu = - Tf. L.X. $$
In a chart $(u,U)$ the second fundamental form is given by
$$s_{ij}=s(\p_i,\p_j)=\g(
\nabla_{\p_i}Tf.\p_j,\nu)=\g\Big(\frac{\p^2
f}{\p_i\p_j},\nu\Big),$$ 
and the mean curvature by
$ \on{Tr}(L)=\sum_{i,j} g^{ij}s_{ij}$.

\subsection{Directional derivatives of functions}\label{no:di}

We will use the following ways to denote directional derivatives of functions, in particular in 
infinite dimensions.
Given a function $F(x,y)$, for instance,
we will write
$$ D_{(x,h)}F \text{ as shorthand for } \partial_t|_0 F(x+th,y).$$
Here $(x,h)$ in the subscript denotes the tangent vector with footpoint $x$ and direction $h$. 
If $F$ takes values in some linear space, we will identify this linear space and its tangent space. 

\medskip
\section{Variational formulas}\label{variation}{\hfill}\par
\smallskip

Recall that many operators such as
$$g=f^*\g, \quad S=S^f, \quad L=L^f,\quad \on{vol}(g), \quad \nabla=\nabla^g, \quad \Delta=\Delta^g$$ 
depend on the immersion $f$. We want to calculate their derivative 
with respect to $f$, which we call \emph{the first  variation.} 
We will use these formulas to calculate the metric gradients that are needed for the geodesic equation.

Some of the formulas 
can also be found in \cite{Michor119,Besse2008, Michor102, Verpoort2008,Harms2010}.

\subsection{Paths of immersions}\label{variation:paths}

All of the concepts introduced in section \ref{no}
can be recast for a path of immersions instead of a fixed immersion. 
This allows us to study variations of immersions. 
So let $f:\R \to \on{Imm}(M,\R^n)$ be a path of immersions. By convenient calculus
\cite{MichorG}, $f$ can equivalently be seen as $f:\R \x M \to \R^n$ 
such that $f(t,\cdot)$ is an immersion for each $t$. 
We can replace bundles over $M$ by bundles over $\R \x M$:
\begin{equation*}\xymatrix{
\on{pr}_2^* T^r_s M \ar[d] & 
\on{pr}_2^* T^r_s M \otimes f^*T\R^n \ar[d] &
\Nor(f) \ar[d]\\
\R \x M & \R \x M & \R \x M
}\end{equation*}
Here $\on{pr}_2$ denotes the projection $\on{pr}_2:\R \x M \to M$.
The covariant derivative $\nabla_Z h$ is now defined for vector fields $Z$ on $\R \x M$ 
and sections $h$ of the above bundles. 
The vector fields $(\p_t, 0_M)$ and $(0_{\R}, X)$, where $X$ is a vector field on $M$, are of
special importance. Let
$$\on{ins}_t : M \to \R \x M, \qquad x \mapsto (t,x) .$$
Then by \cite[22.9.6]{MichorH} one has for vector fields $X,Y$ on $M$
\begin{align*}
\nabla_X Tf(t,\cdot).Y &= \nabla_X T(f \o \on{ins}_t) \o Y = \nabla_X Tf \o T\on{ins}_t \o Y
\\&= \nabla_X Tf \o (0_\R,Y) \o \on{ins}_t 
= \nabla_{T\on{ins}_t \o X} Tf \o (0_\R,Y)\\&
 = \big(\nabla_{(0_\R,X)} Tf \o (0_\R,Y)\big) \o \on{ins}_t .
\end{align*}
This shows that one can recover the static situation at $t$ by using vector fields on $\R \x M$ 
with vanishing $\R$-component and evaluating at $t$. 

\subsection{Setting for first variations}\label{variation:setting_for_1st}
In the remainder of this section, let $f$ be an immersion and $f_t \in T_f\Imm$ a tangent vector to $f$. 
The reason for calling the tangent vector $f_t$ is that in calculations  
it will often be the derivative of a curve of immersions through $f$. 
Using the same symbol $f$ for the fixed immersion
and for the path of immersions through it, one has in fact that
$$D_{(f,f_t)} F = \p_t F(f(t)).$$ 
For the sake of brevity we will write
$\p_t$ instead of $(\p_t,0_M)$ and $X$ instead of $(0_\R,X)$, where $X$ is a vector field on $M$.

Let the smooth mapping $F:\on{Imm}(M,N) \to \Gamma(T^r_s M)$ take values in some space of tensor fields over $M$, 
or more generally in any natural bundle over $M$; see \cite{MichorF}.
\subsection{Lemma (tangential variation of equivariant tensor fields)}{\em 
\label{variation:tangential}
If $F$ is equivariant 
with respect to pullbacks by 
diffeomorphisms of $M$, i.e., 
$$F(f)=(\ph^* F)(f)=\ph^* \Big(F\big((\ph\i)^*f\big)\Big) $$ 
for all $\ph \in \on{Diff}(M)$ and $f \in \on{Imm}(M,N)$,
then the tangential variation of $F$ is its Lie derivative:
\begin{align*}
D_{(f,Tf.f_t^\top)} F&=
\p_t|_0 F\Big(f \o \Fl^{f_t^\top}_t\Big)=
\p_t|_0 F\Big((\Fl^{f_t^\top}_t)^* f\Big)\\&=
\p_t|_0 \Big(\Fl_t^{f_t^\top}\Big)^* \big(F(f)\big) = \L_{f_t^\top}\big(F(f)\big).
\end{align*}
}

Here $\Fl^{f_t^\top}_t$ denotes the flow of $f_t^\top$ and $\L_{f_t^\top}$ denotes 
the Lie derivative along the vector field $f_t^\top$ on $M$.
This allows us to calculate the tangential variation of the pullback metric and 
the volume density, because these tensor fields are natural with respect to pullbacks by diffeomorphisms.

\subsection{Lemma (variation of the metric)}{\em 
\label{variation:metric}
The differential of the pullback metric
\begin{equation*}\left\{ \begin{array}{ccl}
\on{Imm} &\to &\Gamma(S^2_{>0} T^*M),\\
f &\mapsto &g=f^*\g
\end{array}\right.\end{equation*}
is given by
\begin{equation*}
D_{(f,f_t)} g=  -2 \g(f_t,\nu).s + \mathcal L_{f_t^\top}(g).
\end{equation*}
}

{\em Proof}.
Let $f:\R \x M \to \R^n$ be a path of immersions. Swapping covariant derivatives as in 
section \ref{no:sw}, formula \eqref{no:sw:tor_with_t}, one gets
\begin{align*}
\p_t\big(g(X,Y)\big) &= \p_t\big( \g( Tf.X,Tf.Y ) \big)
= \g( \nabla_{\p_t}Tf.X,Tf.Y ) + \g( Tf.X, \nabla_{\p_t}Tf.Y )\\
&=\g( \nabla_X f_t,Tf.Y ) + \g( Tf.X, \nabla_Y f_t ) = \big(2 \on{Sym}\g(\nabla f_t,Tf)\big)(X,Y).
\end{align*}
Splitting $f_t$ into its normal and tangential parts yields
\begin{align*}
2 \on{Sym}\g(\nabla f_t,Tf) &=
2 \on{Sym}\g(\nabla f_t^\bot + \nabla Tf.f_t^\top,Tf) \\&=
-2 \on{Sym}\g(f_t^\bot,\nabla Tf)+2 \on{Sym} g(\nabla f_t^\top,\cdot) \\&=
-2 \g(f_t^\bot,S)+2 \on{Sym} \nabla (f_t^\top)^\flat .
\end{align*}
Finally, the relation
$$D_{(f,Tf.f_t^\top)} g = 2 \on{Sym} \nabla (f_t^\top)^\flat = \L_{f_t^\top} g $$
follows from the equivariance of $g$ (see section \ref{variation:tangential}). 
\qquad\qedhere\par

\subsection{Lemma (variation of the inverse of the metric)}{\em 
\label{variation:inverse_metric}
The differential of the inverse of the pullback metric
\begin{equation*}\left\{ \begin{array}{ccl}
\on{Imm} &\to &\Gamma(L(T^*M,TM)),\\
f &\mapsto &g\i=(f^*\g)\i
\end{array}\right.\end{equation*}
is given by
\begin{equation*}
D_{(f,f_t)} g\i=  -2 \g(f_t,\nu).L.g\i + \mathcal L_{f_t^\top}(g\i).
\end{equation*}
}

{\em Proof}.
\begin{align*}
\p_t g\i &= - g\i (\p_t g ) g\i
 = -g\i \big(-2 \g(f_t^\bot,S)+ \L_{f_t^\top} g\big) g\i \\
& = 2 g\i \g(f_t^\bot,S) g\i -g\i (\L_{f_t^\top} g) g\i
= 2 \g(f_t^\bot,g\i S g\i)+ \L_{f_t^\top} (g\i).\qquad \qedhere
\end{align*}

\subsection{Lemma (variation of the volume density)}{\em 
\label{variation:volume_form}
The differential of the volume density
\begin{equation*}
\left\{ \begin{array}{ccl}
\on{Imm} &\to &\Ga(\on{Vol}(M)),\\
f &\mapsto &\on{vol}(g)=\on{vol}(f^*\g)
\end{array}\right.\end{equation*}
is given by
\begin{equation*}
D_{(f,f_t)} \on{vol}(g)= \Big(\on{div}^{g}(f_t^{\top})-\g(f_t^{\bot},\nu).\Tr(L)\Big) \on{vol}(g).
\end{equation*}
}
\par {\em Proof}.
Let $g(t) \in \Ga(S^2_{>0}T^*M)$ be any curve of Riemannian metrics. Then
$$\p_t \on{vol}(g)=\frac{1}{2}\on{Tr}(g\i.\p_t g)\on{vol}(g).$$
This follows from the formula for $\on{vol}(g)$ in a local oriented chart
$(u^1,\ldots, u^n)$ on $M$:
\begin{align*}
\p_t\on{vol}(g)&=\p_t \sqrt{\det( (g_{ij})_{ij})}\ du^1\wedge\cdots\wedge du^{n-1}\\
&=\frac{1}{2\sqrt{\det ((g_{ij})_{ij})}}\on{Tr}(\on{adj}(g) \p_t g)\
du^1\wedge\cdots\wedge du^{n-1}\\
&=\frac{1}{2\sqrt{\det ((g_{ij})_{ij})}}\on{Tr}(\det((g_{ij})_{ij})g^{-1}\p_t g)\
du^1\wedge\cdots\wedge du^{n-1}\\
&=\frac{1}{2}\on{Tr}(g\i.\p_t g)\on{vol}(g).
\end{align*}
Now we can set $g = f^*\g$ and plug in the formula
for $\p_t g=\p_t(f^*\g)$. This yields
\begin{align*}
&\p_t \on{vol}(g) =\tfrac12\on{Tr}(g\i(-2 \g(f_t,\nu) .s + \mathcal
L_{h^\top}g)).\on{vol}(g)\\
&\qquad= -\g(f_t,\nu)\on{Tr}(g\i . s).\on{vol}(g) +
\tfrac12\on{Tr}(g\i\mathcal L_{h^\top}g).\on{vol}(g).
\end{align*}
The same calculation as above with $\p_t$ replaced by $\L_{h^\top}$ shows that
$$\L_{h^\top} \on{vol}(g) =\tfrac12\on{Tr}(g\i\mathcal L_{h^\top}g).\on{vol}(g).$$
Therefore,
\begin{align*}
\p_t \on{vol}(g) &= -\g(f_t,\nu) \on{Tr}(L).\on{vol}(g) + \mathcal L_{h^\top}(\on{vol}(g)) \\ &= 
-\g(f_t,\nu) \on{Tr}(L).\on{vol}(g)+\on{div}^{g}(h^\top) \on{vol}(g).\qquad\qedhere
\end{align*}

\subsection{Lemma (variation of the volume)}{\em 
\label{variation:volume}
The differential of the total volume
\begin{equation*}
\left\{ \begin{array}{ccl}
\on{Imm} &\to &\R,\\
f &\mapsto &\on{Vol}(f)=\int_M\on{vol}(f^*\g)
\end{array}\right.\end{equation*}
is given by
\begin{equation*}
D_{(f,f_t)} \on{vol}(g)= -\int_M \g(f_t^{\bot},\nu).\Tr(L) \on{vol}(g).
\end{equation*}
}

\subsection{Lemma (variation of the  second fundamental form)}{\em 
\label{variation:second_fundamental}
The differential of the second fundamental form
\begin{equation*}
\left\{ \begin{array}{ccl}
\on{Imm} &\to &\Gamma(S^2 T^*M),\\
f &\mapsto &s^f
\end{array}\right.\end{equation*}
is given by
\begin{align*}
D_{(f,f_t)} s  &=\g(\nabla^2 f_t,\nu)
= \nabla^2 \g(f_t,\nu) - \g(f_t,\nu).g(L,L) + \L_{f_t^\top}.s.
\end{align*}
}
\par {\em Proof}.
By definition $s(X,Y)=\g\big(S(X,Y),\nu\big)=\g\big(\nabla_X (Tf.Y)-Tf.\nabla_X Y,\nu\big).$ Interchanging covariant derivatives as in section \ref{no:sw}, formulas \eqref{no:sw:tor_with_t} 
and \eqref{no:sw:r_with_t}, yields
\begin{align*}
\p_t s(X,Y)&=\g\big(\p_t S(X,Y),\nu\big)+\g\big(S(X,Y),\p_t\nu\big)\\
&=\g\big(\nabla_X\nabla_Y Tf.\p_t-\nabla_{\nabla_X Y}Tf.\p_t,\nu\big)+0=
\g\big(\nabla^2_{X,Y} f_t,\nu\big),
\end{align*}
where the term $\g\big(S(X,Y),\p_t\nu\big)$ vanishes since $\p_t \nu$ is tangential (see section \ref{variation:normal_vector}). For the normal part this yields the following: 
To get the second formula we calculate
\begin{align*}
(D_{(f,\g(f_t,\nu).\nu)} s)(X,Y)&=\g\Big(\nabla_{X,Y}^2 \big(\g(f_t,\nu).\nu\big),\nu\Big)\\&=
\nabla_{X,Y}^2 \g(f_t,\nu)+0+\g(f_t,\nu).\g\big(\nabla_{X,Y}^2 \nu,\nu\big)\\&=
\nabla_{X,Y}^2 \g(f_t,\nu)-\g(f_t,\nu).\g\big(\nabla_{X} \nu,\nabla_Y \nu\big)+0\\&=
\nabla_{X,Y}^2 \g(f_t,\nu)-\g(f_t,\nu).g\big(LX ,LY\big)
\end{align*}
By section \ref{variation:tangential}, the formula for the tangential variation 
follows from the equivariance of the second fundamental form:
\begin{align*}
s^{f \o \phi}(X,Y) &= \g\big(\nabla_X T(f \o \phi) \o Y, \nu^{f \o \phi}\big) =
\g\big( \nabla_X (Tf \o (\phi_* Y) \o \phi), \nu^f \o \phi) \\ & =
\g\big(\nabla_{T\phi \o X} Tf \o (\phi_* Y), \nu^f \o \phi ) \o \phi \\ &=
\g\big( \nabla_{\phi_* X} Tf \o (\phi_* Y), \nu^f )\o \phi =
s^f(\phi_*X,\phi_*Y) \o \phi = (\phi^* s^f)(X,Y).
\qquad\qedhere\end{align*}

\subsection{Lemma (variation of the Weingarten map)}{\em 
\label{variation:weingarten}
The differential of the Weingarten map
\begin{equation*}
\left\{ \begin{array}{ccl}
\on{Imm} &\to &\Gamma(\on{End}(TM)),\\
f &\mapsto &L^f
\end{array}\right.\end{equation*}
is given by
\begin{equation*}\begin{aligned}
D_{(f,f_t)} L =g\i .\nabla^2\big(\g(f_t,\nu)\big) + \g(f_t,\nu) L^2  +\L_{f_t^\top}(L).
\end{aligned}\end{equation*}
}
\par {\em Proof}.
From $L=g\i . s$ follows
\begin{align*}
\p_t L&=g\i \p_t s +\p_t(g\i) s\\
&=g\i  \Big(\nabla^2(\g(f_t,\nu)) - \g(f_t,\nu) g  L^2+\L_{f_t^\top}(s)\Big)
+\Big( 2 \g(f_t,\nu) L  g\i + \mathcal L_{f_t^\top}(g\i) \Big).s\\
&=g\i \nabla^2\big(\g(f_t,\nu)\big)+\g(f_t,\nu) L^2 +\L_{f_t^\top}(L).
\qquad\qedhere\end{align*}

\subsection{Lemma (variation of the mean curvature)}{\em 
\label{variation:scalar_mean_curvature}
The differential of the mean curvature
\begin{equation*}
\left\{ \begin{array}{ccl}
\on{Imm} &\to &C^\infty(M),\\
f &\mapsto &\on{Tr}(L^f)
\end{array}\right.\end{equation*}
is given by
\begin{align*}
D_{(f,f_t)} \on{Tr}(L) &= -\Delta\big(\g(f_t,\nu)\big) + \g(f_t,\nu).\on{Tr}\big(L^2\big)       +d\big(\on{Tr}(L)\big)(f_t^\top).
\end{align*}
}
\par {\em Proof}.
This statement follows from the linearity of the trace operator and 
from the previous equation for $D_{(f,f_t)} L$.
\qquad\qedhere\par

\subsection{Lemma (variation of the normal vector field)}{\em 
\label{variation:normal_vector}
When $f$ is a curve of immersions, the normal vector field to $f$ is a smooth map $\nu:\R \x M \to \R^n$.
Therefore, as explained in section~\ref{no:co},
we can take its covariant derivative in the direction of vector fields on $\R \x M$. 
Identifying $\p_t$ with the vector field $(\p_t,0_M)$ on $\R \x M$, we get
\begin{equation*}
\nabla_{\p_t} \nu = -Tf.\Big(L f_t^\top +
\on{grad}^{g}\big( \g(f_t,\nu)\big)\Big).
\end{equation*}
}

\par {\em Proof}.
$\nabla_{\p_t}\nu$ is tangential because
$\g(\nabla_{\p_t}\nu,\nu)=\tfrac12 \p_t \g(\nu,\nu)= 0$.
Therefore, one can write
$\nabla_{\partial_t}\nu= Tf. (\nabla_{\partial_t}\nu)^{\top}$.
Then for all $X\in \X (M)$ we have
\begin{align*}
g((&\nabla_{\p_t}\nu)^{\top} ,X)
=\g\big(\nabla_{\p_t}\nu,Tf.X\big)
=0-\g\big( \nu, \nabla_{\p_t}Tf.X\big)
=-\g(\nu, \nabla_X Tf.\p_t),
\end{align*}
where in the last step we swapped $X$ and $\p_t$ as
in section~\ref{no:sw}, formula \eqref{no:sw:tor_with_t}.
Splitting into normal and tangential parts yields
\begin{align*}
g((\nabla_{\p_t}\nu)^{\top} ,X)
&=-\g(\nu, \nabla_X f_t)=-\g\Big(\nu, \nabla_X(Tf.f_t^\top+\g(f_t,\nu).\nu)\Big)\\&
=-\g\Big(\nu, \nabla_X(Tf.f_t^\top+\g(f_t,\nu).\nu)\Big)\\&
=-s(X,f_t^\top)-\nabla_X\big(\g(f_t,\nu)\big)-0\\&
=-g\Big(L f_t^\top+\on{grad}^g\g(f_t,\nu),X\Big).
\qquad\qedhere\end{align*}

\subsection{Lemma (variation of the covariant derivative)}{\em 
\label{variation:co}
Let $\nabla=\nabla^g=\nabla^{f^*\g}$ be the 
Levi Civita covariant derivative acting on vector fields on $M$. 
Since any two covariant derivatives on $M$ differ by a tensor field, 
the first variation of $\nabla^{f^*\g}$ is tensorial. It is given by the 
tensor field $D_{(f,f_t)} \nabla^{f^*\g} \in \Ga(T^1_2 M)$,
which is determined by the following relation
holding for vector fields $X,Y,Z$ on $M$:
\begin{multline*}
g\big((D_{(f,f_t)} \nabla)(X, Y),Z\big) = 
\frac12 (\nabla D_{(f,f_t)} g)\big( X \otimes Y \otimes Z
+ Y \otimes X \otimes Z - Z \otimes X \otimes Y \big).
\end{multline*}
}

\par {\em Proof}.
The defining formula for the covariant derivative is
\begin{align*}
g(\nabla_X Y,Z)&= \frac12 \Big[ Xg(Y,Z)+Yg(Z,X)-Zg(X,Y)\\&\qquad
-g(X,[Y,Z])+g(Y,[Z,X])+g(Z,[X,Y]) \Big].
\end{align*}
Taking the derivative $D_{(f,f_t)}$ yields
\begin{multline*}
(D_{(f,f_t)}g)(\nabla_X Y,Z)+g\big((D_{(f,f_t)}\nabla)(X, Y),Z\big)\\ 
\begin{aligned}
=\frac12 \Big[ & X\big((D_{(f,f_t)}g)(Y,Z)\big)+Y\big((D_{(f,f_t)}g)(Z,X)\big)-Z\big((D_{(f,f_t)}g)(X,Y)\big)\\&
-(D_{(f,f_t)}g)(X,[Y,Z])+(D_{(f,f_t)}g)(Y,[Z,X])+(D_{(f,f_t)}g)(Z,[X,Y]) \Big].
\end{aligned}
\end{multline*}
Then the result follows by replacing all Lie brackets in the above formula by covariant derivatives using 
$[X,Y]=\nabla_X Y - \nabla_Y X$
and by expanding all terms of the form $X\big((D_{(f,f_t)}g)(Y,Z)\big)$ using
\begin{align*}
&X\big((D_{(f,f_t)}g)(Y,Z)\big)\\&\qquad\qquad
=(\nabla_X D_{(f,f_t)}g)(Y,Z)
+(D_{(f,f_t)}g)(\nabla_X Y,Z)
+(D_{(f,f_t)}g)(Y,\nabla_X Z).
\qquad\qedhere\end{align*}

\subsection{Setting for second variations}\label{variation:setting_for_2nd}

All formulas for second derivatives will be used in
section~\ref{sectional_curvature:2nd_derivative}.
There we consider a curve of immersions
$$f(t,x) =f_0(x)+t.a(x).\nu^{f_0}(x)$$
for a fixed immersion $f_0$. This curve of immersions has the
property that at $t=0$ its first derivative and the covariant derivative of the first derivative 
are both horizontal, i.e.,
\begin{equation}\label{variation:setting_for_2nd:two_properties}
f|_{t=0}=f_0, \quad
\p_t|_0 f=a.\nu^{f_0} , \quad\text{and}\quad
\nabla_{\p_t} Tf.\p_t |_{t=0} = 0.
\end{equation}
In all calculations 
of second variations we will assume that  the above properties hold.

\subsection{Lemma (second variation of the metric)}{\em 
\label{variation:2nd_metric}
The second derivative of the pullback metric
\begin{equation*}\left\{ \begin{array}{ccl}
\on{Imm} &\to &\Gamma(S^2_{>0} T^*M),\\
f &\mapsto &g=f^*\g
\end{array}\right.\end{equation*}
along a curve of immersions $f$ satisfying 
property~\eqref{variation:setting_for_2nd:two_properties} from section~\ref{variation:setting_for_2nd} is given by
\begin{equation*}
\p_t^2|_0 f^*\g=2 (da \otimes da)+2a^2 g_0(L^{f_0}, L^{f_0}).
\end{equation*}
}
\par {\em Proof}.
Since $\nabla_{\p_t}Tf.\p_t|_0=0$, we have
\begin{align*}
\p_t^2|_0 g(X,Y)&=\p_t^2|_0 \g( Tf.X,Tf.Y) \\
&= \p_t|_0 \g(\nabla_{\p_t}Tf.X,Tf.Y ) + \p_t|_0 \g(
Tf.X,\nabla_{\p_t}Tf.Y ) \\
&= 2\g(\nabla_{\p_t}Tf.X|_0,\nabla_{\p_t}Tf.Y|_0 )+0+0
 = 2\g( \nabla_X Tf.\p_t,\nabla_Y Tf.\p_t ).
\end{align*}
Using $Tf.\p_t=a.\nu^{f_0}$, we get
\begin{align*}
\p_t^2|_0 (g(X,Y))
&=2 da(X).da(Y)+2a^2 \g( \nabla_X \nu^{f_0},\nabla_Y
\nu^{f_0})\\
&=2 (da \otimes da)(X,Y)+2a^2 .g_0(L^{f_0}X,L^{f_0}Y).
\qquad\qedhere\end{align*}

\subsection{Lemma (second variation of the inverse metric)}{\em 
\label{variation:2nd_inverse_metric}
The second derivative of the inverse of the pullback metric
\begin{equation*}\left\{ \begin{array}{ccl}
\on{Imm} &\to &\Gamma(L(T^*M,TM)),\\
f &\mapsto &g\i=(f^*\g)\i
\end{array}\right.\end{equation*}
along a curve of immersions $f$ satisfying 
property~\eqref{variation:setting_for_2nd:two_properties} from section~\ref{variation:setting_for_2nd} is given by
\begin{align*}
\p_t^2|_0 (f^*\g)\i =
6a^2 (L^{f_0})^2.g_0\i
-2 g_0\i(da \otimes da)g_0\i.
\end{align*}
}
\par {\em Proof}.
We look at $g=f^*\g$ as a bundle map from $TM$ to $T^*M$. Then
\begin{align*}
\p_t^2|_0 (g\i) &= \p_t|_0 \big( -g\i.\p_t g.g\i \big) =
2 g_0\i.\p_t|_0 g.g_0\i.\p_t|_0 g.g_0\i -g_0\i.\p_t^2|_0 g.g_0\i \\ &=
2 (-2a L^{f_0})^2.g_0\i -g_0\i.\big(2 (da \otimes da)+2a^2 g_0 \o (L^{f_0} \otimes L^{f_0})\big).g_0\i \\ &=
8 a^2 (L^{f_0})^2.g_0\i -2g_0\i(da \otimes da)g_0\i-2a^2 (L^{f_0})^2 .g_0\i. \qquad\qedhere
\end{align*}

\subsection{Lemma (second variation of the volume form)}{\em 
\label{variation:2nd_volume_form}
The second derivative of the volume form
\begin{equation*}\left\{ \begin{array}{ccl}
\on{Imm} &\to &\Ga(\Vol(M)),\\
f &\mapsto &\on{vol}(g)=\on{vol}(f^*\g)
\end{array}\right.\end{equation*}
along a curve of immersions $f$ satisfying 
property~\eqref{variation:setting_for_2nd:two_properties} from section~\ref{variation:setting_for_2nd} is given by
\begin{equation*}\begin{aligned}
\p_t^2|_0 \on{vol}(g) 
=\Big[ a^2 \on{Tr}(L^{f_0})^2- a^2 \on{Tr}\big( (L^{f_0})^2 \big)+\norm{da}^2_{g_0\i} \Big] \on{vol}(g_0).
\end{aligned}\end{equation*}
}

{\em Proof}.
In section \ref{variation:volume_form} we showed that for any curve  $g(t) \in \Ga(S^2_{>0}T^*M)$ of
Riemannian metrics, we have
$$\p_t \on{vol}(g)=\frac{1}{2}\on{Tr}\big(g\i.\p_t g\big)\on{vol}(g).$$
Therefore,
\begin{align*}
\p_t^2 \on{vol}(g)&=\p_t \frac{1}{2}\on{Tr}\big(g\i.\p_t g\big)\on{vol}(g)
= \frac{1}{2}\on{Tr}\big(\p_t(g\i).\p_t g\big) \on{vol}(g) \\
&+ \frac{1}{2}\on{Tr}\big(g\i.\p_t^2 g\big) \on{vol}(g)
+ \frac{1}{2}\on{Tr}\big(g\i.\p_t g\big) \p_t \on{vol}(g).
\end{align*}
Evaluating at $t=0$ and setting $g(t)=f^*\g$, we get
\begin{align*}
\p_t^2|_0 \on{vol}(g)
=& \frac{1}{2}\on{Tr}\big((2 a L^{f_0}  g_0\i).(-2 a .s^{f_0})\big)
\on{vol}(g_0) \\
\qquad&+  \frac{1}{2}\on{Tr}\big(g_0\i.2 (da\otimes da)\big) \on{vol}(g_0)\\
\qquad&+  \frac{1}{2}\on{Tr}\big(g_0\i.2a^2
g_0.(L^{f_0})^2\big) \on{vol}(g_0)\\
\qquad&+  \frac{1}{2}\on{Tr}\big(g_0\i.(-2 a .s^{f_0})\big) (-
\on{Tr}(L^{f_0}) . a) \on{vol}(g_0) \\
=&  \Big[ a^2 \on{Tr}(L^{f_0})^2- a^2 \on{Tr}\big( (L^{f_0})^2 \big)+\norm{da}^2_{g\i} \Big] \on{vol}(g_0).
\qquad \qedhere \end{align*}

\subsection{Lemma (Second variation of the second fundamental form)}{\em 
\label{variation:2nd_second_fundamental_form}
The second derivative of  the second fundamental form
\begin{equation*}
\left\{ \begin{array}{ccl}
\on{Imm} &\to &\Gamma(S^2 T^*M),\\
f &\mapsto &s^f
\end{array}\right.\end{equation*}
along a curve of immersions $f$ satisfying property \ref{variation:setting_for_2nd:two_properties} from section \ref{variation:setting_for_2nd} is given by
\begin{align*}
\p_t^2|_0 s &=2(da\otimes da)\bigl(\on{Id}\otimes L^{f_0} + L^{f_0}\otimes \on{Id}\bigr)
-\norm{da}^2_{g_0\i}.s^{f_0}\\
&\qquad+2.a(\nabla_{\grad^{g_0}(a)}s^{f_0}).
\end{align*}
}
\par {\em Proof}.
From section~\ref{variation:second_fundamental} we have
$$\p_t s(X,Y)= \g( \nabla^2_{X,Y} Tf.\p_t, \nu )
=\g( \nabla^2_{X,Y} f_t, \nu ).$$
Using $\nabla_{\p_t}f_t=0$, we get
\begin{align*}
&\p_t^2 s(X,Y)= \g( \nabla^2_{X,Y} f_t, \nabla_{\p_t}\nu ) +
\g( \nabla_{\p_t} \nabla_X \nabla_Y f_t 
- \nabla_{\p_t} \nabla_{\nabla_X Y} f_t, \nu ) 
\\&\qquad=
\g( \nabla^2_{X,Y} f_t, \nabla_{\p_t}\nu )+ 
0 - \g(  \nabla_{\nabla_X Y} \nabla_{\p_t} f_t
+\nabla_{[\p_t,\nabla_X Y]}f_t,\nu)
\\&\qquad=
\g( \nabla^2_{X,Y} f_t, \nabla_{\p_t}\nu )+
0- \g(\nabla_{[\p_t,\nabla_X Y]}f_t,\nu)
\\&\qquad=
\g( \nabla^2_{X,Y} f_t, \nabla_{\p_t}\nu )
-\g(\nabla_{(D_{(f,f_t)}\nabla)(X,Y)}f_t,\nu).
\end{align*}
In the last step we used
$$
[\p_t,\nabla_X Y] = [(\p_t,0_M),(0_\R,\nabla_X Y)]
=\big(0_\R,(D_{(f,f_t)}\nabla)(X,Y)\big)=(D_{(f,f_t)}\nabla)(X,Y).
$$
Evaluating at $t=0$ yields
\begin{align*}
&\p_t^2 |_0s(X,Y) \\&\quad=
\g( \nabla^2_{X,Y} (a.\nu^{f_0}), -Tf_0.\on{grad}^{g_0} a) 
-\g(\nabla_{(D_{(f_0,a.\nu^{f_0})}\nabla)(X,Y)}(a.\nu^{f_0}),\nu^{f_0})
\\&\quad=0+\g( da(X).\nabla_Y\nu^{f_0}+da(Y).\nabla_{X}\nu^{f_0},
 -Tf_0.\on{grad}^{g_0} a) \\&\qquad +
 \g( a.\nabla^2_{X,Y} (\nu^{f_0}), -Tf_0.\on{grad}^{g_0} a)
- da\big((D_{(f_0,a.\nu^{f_0})}\nabla)(X,Y)\big)+0.
\end{align*}
We will treat the three terms separately. Using
$\nabla_Z\nu=-Tf.L.Z$, one gets for the first term
\begin{align*}
&\g( da(X).\nabla_Y\nu^{f_0}+da(Y).\nabla_{X}\nu^{f_0},
 -Tf_0.\on{grad}^{g_0} a)
\\&\qquad=g_0(da(X)L^{f_0}Y+da(Y)L^{f_0}X,\on{grad}^{g_0} a)\\&\qquad
=da(X). da(L^{f_0}Y) + da(Y).da(L^{f_0}X).
\end{align*}
For the second term one gets
\begin{align*}
&\g( a.\nabla^2_{X,Y} (\nu^{f_0}), -Tf_0.\on{grad}^{g_0} a)
=-a\g( \nabla_X\nabla_Y \nu^{f_0}-
\nabla_{\nabla_X Y}\nu^{f_0},Tf_0.\on{grad}^{g_0} a)
\\&\qquad
=-a\g( -\nabla_X(Tf_0L^{f_0}Y)+
Tf_0L^{f_0}\nabla_X Y,Tf_0.\on{grad}^{g_0} a)
\\&\qquad
=-a\g( -(\nabla Tf_0)(X,L^{f_0}Y)-Tf_0\nabla_X(L^{f_0}Y)+
Tf_0L^{f_0}\nabla_X Y,Tf_0.\on{grad}^{g_0} a)
\\&\qquad
=0+a\g( Tf_0(\nabla_X L^{f_0})(Y),Tf_0.\on{grad}^{g_0} a)
=ag_0\big((\nabla_X L^{f_0})(Y),\on{grad}^{g_0} a\big)
\\&\qquad
=a\nabla_X\big(g_0(L^{f_0}Y,\on{grad}^{g_0} a)\big)
-ag_0(L^{f_0}\nabla_X Y,\on{grad}^{g_0} a)\\&\qquad\qquad
-ag_0(L^{f_0}Y,\nabla_X\on{grad}^{g_0} a)
\\&\qquad
=a\nabla_X\big(s^{f_0}(Y,\on{grad}^{g_0} a)\big)
-as^{f_0}(\nabla_X Y,\on{grad}^{g_0} a)
-as^{f_0}(Y,\nabla_X\on{grad}^{g_0} a)
\\&\qquad
=a(\nabla_X s)(Y,\on{grad}^{g_0} a).
\end{align*}
$\nabla^2_{X,Y}\nu$ is symmetric in $X,Y$ because the ambient space $\R^n$ is flat. 
Therefore, the last formula and the symmetry of $s$ imply that
\begin{align*}
a(\nabla_X s)(Y,\on{grad}^{g_0} a)=a(\nabla_Y s)(X,\on{grad}^{g_0} a)=
a(\nabla_{\on{grad}^{g_0} a} s)(X,Y).
\end{align*}
The third term yields, using the formula in section~\ref{variation:co}
\begin{align*}
&-g_0\big((D_{(f,a.\nu^{f_0})}\nabla)(X,Y),\grad^{g_0}(a)\big)
\\&\quad=
-\frac12(\nabla(-2a.s^{f_0}))(X, Y, \grad^{g_0}(a))
-\frac12(\nabla(-2a.s^{f_0}))(Y, X,\grad^{g_0}(a))
\\&\quad\qquad
+\frac12(\nabla(-2a.s^{f_0}))(\grad^{g_0}(a), X, Y)
\\&\quad=
da(X).s^{f_0}(Y , \grad^{g_0}(a))+
a.(\nabla_Xs^{f_0})(Y, \grad^{g_0}(a))
\\&\quad\qquad+
da(Y).s^{f_0}(X , \grad^{g_0}(a))+
a.(\nabla_Y s^{f_0})(X, \grad^{g_0}(a))
\\&\quad\qquad
-da(\grad^{g_0}(a)).s^{f_0}(X ,Y)-
a.(\nabla_{\grad^{g_0}(a)}s^{f_0})(X,Y)
\\&\quad=
da(X).da(L^{f_0}Y)+da(Y).da(L^{f_0}X)
\\&\quad\qquad
-\norm{da}_{g\i}^2.s^{f_0}(X ,Y)+
a.(\nabla_{\grad^{g_0}(a)}s^{f_0})(X,Y).\qquad\qedhere
\end{align*}

\subsection{Lemma (second variation of the mean curvature)}{\em 
\label{variation:2nd_mean_curvature}
The second derivative of the mean curvature
\begin{equation*}
\left\{ \begin{array}{ccl}
\on{Imm} &\to &C^\infty(M),\\
f &\mapsto &\on{Tr}(L^f)
\end{array}\right.\end{equation*}
along a curve of immersions $f$  satisfying property \ref{variation:setting_for_2nd:two_properties} from section \ref{variation:setting_for_2nd} is given by
\begin{align*}
\p_t^2|_0 \on{Tr}(L) &=
2a^2\on{Tr}\big( (L^{f_0})^3\big)+4a\on{Tr}\big( L^{f_0}  g_0\i.\nabla^2 a\big)
+2\Tr(g\i(da\otimes da)L^{f_0})\\&\quad-\norm{da}^2_{g_0\i}\Tr(L^{f_0})
+2a\Tr^{g_0}\big(\nabla_{\grad^{g_0}a}s^{f_0}\big).
\end{align*}
}
\par {\em Proof}.
From $\on{Tr}(L)=\on{Tr}(g\i.s)$ we get
\begin{align*}
\p_t^2 \on{Tr}(L) &= \on{Tr}\big(\p_t^2(g\i).s\big)
+2\on{Tr}\big(\p_t(g\i).\p_t s\big)
+\on{Tr}\big(g\i.\p_t^2 s\big).
\end{align*}
Evaluating at $t=0$, we get
\begin{align*}
\p_t^2|_0 \on{Tr}(L) &=
\on{Tr}\big(6a^2 (L^{f_0})^2.g_0\i.s^{f_0}\big)
+\on{Tr}\big(-2g_0\i. (da\otimes da).g_0\i.s^{f_0}\big)\\
&\qquad+2\on{Tr}\big(2 a L^{f_0}  g_0\i.\nabla^2 a \big)
+2\on{Tr}\big(2 a L^{f_0}  g_0\i.(-a g_0(L^{f_0})^2)\big)\\
&\qquad+2.\on{Tr}\big(g_0\i.((da\otimes da\circ L^{f_0}) + (da\circ
L^{f_0}\otimes da))\big)\\&\qquad
-\norm{da}_{g_0\i}^2\Tr(L^{f_0})
+2a\Tr^{g_0}\big(\nabla_{\grad^{g_0}a}s^{f_0}\big)
\\
&=2a^2\on{Tr}\big( (L^{f_0})^3\big)
-2\on{Tr}\big(g_0\i. (da\otimes da).L^{f_0}\big) \\
&\qquad+4a\on{Tr}\big( L^{f_0}  g_0\i.\nabla^2 a \big)
+4\on{Tr}\big(g_0\i.(da\otimes da).(L^{f_0}) \big)
\\&\qquad
-\norm{da}_{g_0\i}^2\Tr(L^{f_0})
+2a\Tr^{g_0}\big(\nabla_{\grad^{g_0}a}s^{f_0}\big).\qquad\qedhere
\end{align*}

\medskip
\section{The geodesic equation on $\on{Imm}(M,\mathbb R^n)$}\label{geodesic_equation_imm}
{\hfill}\par
\smallskip

We recall the definition of the $G^{\Ph}$-metric from section~\ref{introduction:wherethispaper},
$$ G^\Ph_f(h,k) =  \int_{M} \Ph\big(\on{Vol},\on{Tr}(L)(x)\big) \g\big( h(x),k(x)\big)\vol(g)(x).$$ 
We will write $\Ph(\v,\m)$ for the function $\Ph$ and its arguments $(\v,\m) \in \R^2$ corresponding to 
the volume and mean curvature.
%
\subsection{The geodesic equation on $\on{Imm}(M,\mathbb R^n)$}\label{geodesic_equation_imm:geodesic_equation}

We use the method of section~\ref{hamiltonian:geodesic_equation} to calculate the 
geodesic equation. So we need to compute the metric gradients.
The calculation at the same time shows the existence of the gradients. 
Let $m \in T_f\on{Imm}(M,\mathbb R^n)$ with $$\boxed{m=a.\nu^f+Tf.m^\top.}$$ 
To shorten the notation, we will not always note the dependence on $f$ in expressions as $\nu^f, L^f,\ldots.$ 
\begin{align*}
D_{(f,m)} G^\Ph_f(h,k) &=
               \int_M (\p_\v \Ph) (D_{(f,m)} \on{Vol}) \g( h,k ) \on{vol}(g) \\ 
& \qquad           + \int_M (\p_\m \Ph) (D_{(f,m)} \on{Tr}(L))         \g( h,k ) \on{vol}(g) \\ 
& \qquad           + \int_M  \Ph   .                         \g( h,k ) (D_{(f,m)} \on{vol}(g)).
\end{align*}
To read off the $K$-gradient of the metric, we write this expression as
\begin{align*}
\int_M \Ph .\g\left( \Big[ \frac{\p_\v\Ph}{\Ph} (D_{(f,m)} \on{Vol}) 
+ \frac{\p_\m\Ph}{\Ph} (D_{(f,m)} \on{Tr}(L)) + \frac{D_{(f,m)} \on{vol}(g)}{\on{vol}(g)} \Big] h,k \right) \on{vol}(g).
\end{align*}
Therefore, using the formulas from section~\ref{variation}, we can calculate the $K$--gradient:
\begin{align*}
K_f(m,h)&=\Big[ \frac{\p_\v\Ph}{\Ph} (D_{(f,m)} \on{Vol}) 
+ \frac{\p_\m\Ph}{\Ph} (D_{(f,m)} \on{Tr}(L)) +\frac{D_{(f,m)} \on{vol}(g)}{\on{vol}(g)} \Big] h\\
&= \Big[ \frac{\p_\v\Ph}{\Ph} \Big(\int_M -\on{Tr}(L).a \on{vol}(g)\Big)\\
&\qquad + \frac{\p_\m\Ph}{\Ph} \big(-\Delta a+a \on{Tr}(L^2)+d\on{Tr}(L)({m^\top})\big)\\
&\qquad+ \on{div}^g(m^\top)- \on{Tr}(L).a \Big] h.
\end{align*}
To calculate the $H$-gradient, we treat the three summands of $D_{(f,m)} G^\Ph_f(h,k)$ separately. 
The first summand is
\begin{align*}
&\int_M (\p_\v \Ph) (D_{(f,m)} \on{Vol}(x)) \g( h(x),k(x) ) \on{vol}(g)(x)\\ 
&=\int_{x\in M} (\p_\v \Ph)\Big( \int_{y\in M} -\on{Tr}(L(y)).a(y) \on{vol}(g)(y)\Big) \g( h(x),k(x) ) \on{vol}(g)(x) \\
&=\int_{y\in M}\!  \g\left( a(y).\nu(y),-\on{Tr}(L(y)) \int_{x\in M}\!  (\p_\v \Ph)  \g( h(x),k(x) ) \on{vol}(g)(x) . \nu(y) \right)  \on{vol}(g)(y) \\
&=G^\Ph_f\Big(m,-\frac{1}{\Ph} \on{Tr}(L) \int_M (\p_\v \Ph)  \g( h,k ) \on{vol}(g) . \nu \Big).
\end{align*}
By the symmetry of the Laplacian
$$\int_M \Delta(a).b\on{vol}(g)=\int_M a.\Delta(b)\on{vol}(g) \quad \text{for $a,b\in C^\infty(M)$}$$
one gets for the second summand
\begin{align*}
&\int_M (\p_\m \Ph) (D_{(f,m)} \on{Tr}(L))\ \g( h,k ) \on{vol}(g) \\
&\qquad =\int_M (\p_\m \Ph) \big(-\Delta a + a \on{Tr}(L^2)+d\on{Tr}(L)({m^\top})\big) \g( h,k ) \on{vol}(g) \\
&\qquad=\int_M -a. \Delta\big((\p_\m\Ph)\g( h,k ) \big) \on{vol}(g) 
+ \int_M  a. (\p_\m \Ph) \on{Tr}(L^2)\g( h,k ) \on{vol}(g) \\
&\qquad\qquad +\int_M (\p_\m \Ph) g\big(\on{grad}^g(\on{Tr}(L)),m^\top\big)\g( h,k )\on{vol}(g)\\
&\qquad=\int_M -a. \Delta\big((\p_\m\Ph)\g( h,k ) \big) \on{vol}(g) 
+ \int_M  a. (\p_\m \Ph) \on{Tr}(L^2)\g( h,k ) \on{vol}(g) \\
&\qquad\qquad +\int_M (\p_\m \Ph) \g( Tf.\on{grad}^g(\on{Tr}(L)),Tf.m^\top)\g( h,k )\on{vol}(g)\\
&\qquad=      G^\Ph_f\Big(m,-\frac{1}{\Ph} \Delta\big((\p_\m\Ph)\g( h,k )\big) .\nu\Big) 
+G^\Ph_f\Big(m, \frac{1}{\Ph} (\p_\m\Ph)\on{Tr}(L^2)\g( h,k ) . \nu \Big)\\
&\qquad\qquad +G^\Ph_f\Big(m,\frac{1}{\Ph}(\p_\m \Ph)\g( h,k ) Tf.\on{grad}^g(\on{Tr}(L))\Big).
\end{align*}
In the calculation of the third term of the $H_f(m,h)$--gradient, we will make use of the following formula,
which is valid for $\phi\in C^\infty(M)$ and $X \in \X(M)$: 
\begin{align*}
&0=\int_M \on{div}(\phi.X).\on{vol}(g) = \int_M \L_{\phi.X} \on{vol}(g) \\
&\qquad= \int_M (d \o i_{\phi.X} + i_{\phi.X} \o d) \on{vol}(g) 
= \int_M  d(\phi . i_{X} \on{vol}(g) )\\
&\qquad= \int_M d\phi \wedge i_X \on{vol}(g)+\int_M \phi \wedge d (i_X \on{vol}(g)) \\
&\qquad= \int_M \big(-i_X(d\phi \wedge \on{vol}(g))+i_X \o d\phi \wedge  \on{vol}(g)\big) 
+ \int_M \phi . \L_X  \on{vol}(g)  \\
&\qquad= \int_M d\phi(X)  \on{vol}(g)  +\int_M \phi . \on{div}(X)  \on{vol}(g).  
\end{align*}
Therefore, we can calculate the third summand, which is given by
\begin{align*}
&\int_M \Ph \g( h,k ) (D_{(f,m)} \on{vol}(g)) 
=\int_M \Ph \g( h,k ) \big(\on{div}^g(m^\top)- \on{Tr}(L).a\big)  \on{vol}(g)\\ 
&\qquad = -\int_M d(\Ph\g( h,k))(m^\top)\on{vol}(g) +G^\Ph_f(m, -\g( h,k ) \on{Tr}(L) .\nu)\\
&\qquad=  -\int_M \g\Big( Tf.\on{grad}^g(\Ph\g( h,k)),m\Big)\on{vol}(g) +G^\Ph_f(m, -\g( h,k ) \on{Tr}(L) .\nu)\\
&\qquad=  G^\Ph_f\Big(m, -\frac{1}{\Ph} Tf.\on{grad}^g(\Ph\g( h,k))  -\g( h,k ) \on{Tr}(L) .\nu\Big).
\end{align*}
Summing up all the terms the $H$-gradient is given by
\begin{align*}
H_f(h,k)&=\Big[-\frac{1}{\Ph} \on{Tr}(L) \int_M (\p_\v \Ph)  \g( h,k ) \on{vol}(g)-\frac{1}{\Ph} \Delta\big((\p_\m\Ph)\g( h,k )\big)\\
&\qquad+\frac{1}{\Ph} (\p_\m\Ph)\on{Tr}(L^2)\g( h,k )-\g( h,k ) \on{Tr}(L)\Big] \nu^f \\
&\qquad +\frac{1}{\Ph}Tf.\Big[(\p_\m \Ph)\g( h,k )\on{grad}^g(\on{Tr}(L))-\on{grad}^g(\Ph\g( h,k))\Big].
\end{align*}
\subsection{Theorem}{\em 
The geodesic equation for the almost local metrics $G^\Ph$ on $\on{Imm}(M,\mathbb R^n)$ is then given by
$$\boxed{\begin{aligned}
f_t &= h=a.\nu^f+Tf.h^\top,   \\
h_t&=\frac{1}{2}\Big[-\frac{1}{\Ph} \on{Tr}(L) \int_M \p_\v \Ph  \norm{h}^2 \on{vol}(g)-\frac{1}{\Ph} \Delta\big((\p_\m\Ph) \norm{h}^2 \big)\\
&\qquad\qquad+\frac{1}{\Ph} (\p_\m\Ph)\on{Tr}(L^2)\norm{h}^2-\norm{h}^2 \on{Tr}(L)\Big] \nu^f\\
&\qquad\qquad +\frac{1}{2\Ph}Tf.\Big[(\p_\m \Ph)\norm{h}^2\on{grad}^g(\on{Tr}(L))- \on{grad}^g(\Ph\norm{h}^2)\Big]\\
&\qquad-\big[ \frac{\p_\v\Ph}{\Ph} \int_M -\on{Tr}(L).a \on{vol}(g)\\
&\qquad\qquad + \frac{\p_\m\Ph}{\Ph} \big(-\Delta a+a \on{Tr}(L^2)+d\on{Tr}(L)({h^\top})\big)\\
&\qquad\qquad+ \on{div}^g(h^\top)- \on{Tr}(L).a \big] h.
\end{aligned}}$$
}

\subsection{Momentum mappings}\label{geodesic_equation_imm:momentum_mappings}

The metric $G^\Ph$ is invariant under the action of the re\-para\-me\-triza\-tion group
$\on{Diff}(M)$ and under the Euclidean motion group $\mathbb R^n \rtimes
\on{SO}(n)$. According to section~\ref{hamiltonian:momentum_mappings}, the momentum mappings for these group
actions are constant along any geodesic in $\on{Imm}(M,\mathbb R^n)$:
$$\boxed
{\begin{aligned}
\forall X\in\X(M): \int_M \Ph(\on{Vol(f)},\on{Tr}(L^f)) \g( Tf.X,f_t ) \on{vol}(g) && \text{reparam. momenta}\\
\text{or } \Ph(\on{Vol(f)},\on{Tr}(L^f))  g(f_t^\top)  \on{vol}(g) \in\Ga(T^*M\otimes_M\on{Vol}(M))
        && \text{reparam. momentum}\\
\int_M \Ph(\on{Vol(f)},\on{Tr}(L^f)) f_t \on{vol}(g) && \text{linear momentum}\\
\forall X\in \mathfrak{so}(n): \int_M \Ph(\on{Vol(f)},\on{Tr}(L^f)) \g( X.f,f_t )
\on{vol}(g) && \text{angular momenta} \\
\text{or }\int_M \Ph(\on{Vol(f)},\on{Tr}(L^f)) (f\wedge f_t) 
\on{vol}(g) \in {\textstyle\bigwedge^2}\mathbb R^n&& \text{angular momentum.}
\end{aligned}}$$
Here $\Ga(T^*M\otimes_M\on{Vol}(M))\subset \Ga(TM)'$ is the space of cotangent bundle valued densities contained 
in the dual of the Lie algebra $\Ga(TM)$. The name angular momentum is justified by the natural identification 
${\textstyle\bigwedge^2}\mathbb R^n\cong \mathfrak s\mathfrak o(n)\cong \mathfrak s\mathfrak o(n)^* $.

\medskip
\section{The geodesic equation on $B_i(M,\mathbb R^n)$}\label{geodesic_equation_Bi}{\hfill}\par
\smallskip

\subsection{The horizontal bundle and the metric on the quotient space}\label{geodesic_equation_Bi:horizontal}

Since $\on{vol}(f^*\g)$ and $\on{Tr}(L^f)$ react
equivariantly to the action of the group $\on{Diff}(M)$, every $G^\Ph$-metric is
$\on{Diff}(M)$-invariant. As described in
Section~\ref{sh:sub}, it induces a
Riemannian metric on $B_i$ (off the singularities) such that the projection $\pi:\on{Imm}\to B_i$ is a Riemannian 
submersion. 

By definition, a tangent vector $h$ to $f \in \on{Imm}(M,\mathbb R^n)$ is horizontal if
and only if it is $G^\Ph$-perpendicular to
the $\on{Diff}(M)$-orbits. This is the case if and only if $\g( h(x), T_x f .X_x ) =
0$ at every point $x \in M$. Therefore, the horizontal bundle at the point $f$ equals the set of 
sections of the normal bundle (see Section~\ref{no:no}) along $f$. Thus the metric  
on the horizontal bundle is given by 
$$ G^\Ph_f(h^{\hor}, k^{\hor}) = G^\Ph_f(a\cdot \nu, b\cdot \nu) = \int_{M} \Ph(\on{Vol},\on{Tr}(L))\, a.b \on{vol}(g).$$ 

The following lemma shows that  every path in $B_i$ corresponds to exactly one horizontal path in $\on{Imm}$,
and therefore the calculation of the geodesic equation can be done on the horizontal bundle instead of on $B_i$. 
\subsection*{Lemma}{\em
For any smooth path $f$ in $\on{Imm}$ there exists a smooth path $\ph$ 
in $\on{Diff}(M)$ with $\varphi(t,\ )=\on{Id}_M$ depending smoothly on $f$ such that the path 
$f(t,\varphi(t,x))$ is horizontal, i.e., $\p_t f(t,\varphi(t,x))$ lies in the horizontal bundle.
}

The basic idea is to write the path $\varphi$ as the integral curve of a time dependent vector field. 
This method is called the Moser--Trick, (see \cite[Section 2.5]{Michor102}).

\subsection{The geodesic equation on $B_i(M,\mathbb R^n)$}\label{geodesic_equation_Bi:geodesic_equation}

As described in section~\ref{sh:sub}, 
geodesics in $B_i$ correspond to horizontal geodesics in
$\on{Imm}$. A horizontal geodesic $f$ in $\on{Imm}$
has $f_t=a . \nu^f$ with $a \in C^\infty(\mathbb R \x M)$. The geodesic equation is
given by 
$$f_{tt} = \underbrace{a_t . \nu}_{\text{normal}} + \underbrace{a .
\nu_t}_{\text{tang.}} =\frac{1}{2} H(a.\nu,a.\nu)-K(a.\nu,a.\nu);$$
see section~\ref{hamiltonian:geodesic_equation}. 
This equation splits into a normal part and a tangential part. 
From the conservation of the reparametrization momentum (see
section~\ref{hamiltonian:momentum_mappings} and the previous section) 
it follows that the tangential part of the geodesic equation is 
satisfied automatically. We will nevertheless check this by hand.
From section~\ref{geodesic_equation_imm:geodesic_equation}, where we calculated the metric gradients on $\on{Imm}$, we get
\begin{align*}
K_f(a.\nu,a\nu)&=\big[ -\frac{\p_\v\Ph}{\Ph} \int_M \on{Tr}(L).a \on{vol}(g) 
\\&\qquad + \frac{\p_\m\Ph}{\Ph} \big(-\Delta a+a \on{Tr}(L^2)\big)- \on{Tr}(L).a \big] a.\nu\\
H_f(a.\nu,a\nu)&=\frac{1}{\Ph}Tf.\Big[(\p_\m \Ph)a^2 \on{grad}^g(\on{Tr}(L))-\on{grad}^g(\Ph a^2)\Big]\\
&\qquad +\Big[-\frac{1}{\Ph} \on{Tr}(L) \int_M \p_\v \Ph  a^2 \on{vol}(g)-\frac{1}{\Ph} \Delta\big((\p_\m\Ph) a^2\big)\\
&\qquad+\frac{1}{\Ph} (\p_\m\Ph)\on{Tr}(L^2)a^2-a^2 \on{Tr}(L)\Big] \nu.
\end{align*}
From this we can easily read  the tangential part of the geodesic equation
$$a.\nu_t=\frac{1}{2\Ph}Tf.\Big[(\p_\m \Ph)a^2 \on{grad}^g(\on{Tr}(L))-\on{grad}^g(\Ph a^2)\Big].$$
We expand the right--hand side using a Leibnitz rule for the gradient, 
\begin{align*}
\on{grad}^g(a.b)=a\on{grad}^gb+b\on{grad}^ga \quad \text{for $a,b \in C^\infty(M)$}.
\end{align*}
This yields
\begin{align*}
a.\nu_t&=\frac{1}{2\Ph}Tf.\Big[(\p_\m \Ph)a^2 \on{grad}^g(\on{Tr}(L))-\on{grad}^g(\Ph a^2)\Big]\\
&=\frac{1}{2\Ph}Tf.\Big[(\p_\m \Ph)a^2 \on{grad}^g(\on{Tr}(L))-\Ph.\on{grad}^g(a^2)- a^2.\on{grad}^g(\Ph)\Big]\\
&=\frac{1}{2\Ph}Tf.\Big[(\p_\m \Ph)a^2 \on{grad}^g(\on{Tr}(L))-\Ph.\on{grad}^g(a^2)- (\p_\m \Ph)a^2 \on{grad}^g(\on{Tr}(L))\Big]\\
&=-\frac{1}{2\Ph}\Ph Tf.\on{grad}^g(a^2)=-Tf.a.\on{grad}^g(a).
\end{align*}
By the variational formula for $\nu$ in section~\ref{variation:normal_vector} this equation is satisfied automatically.
The normal part is given by
\begin{align*}
a_t &= \g\big( \frac12 H(a.\nu,a.\nu)-K(a.\nu,a.\nu) , \nu \big)\\&=
\frac{1}{\Ph} \big[\frac12 \Ph a^2 \on{Tr}(L^f) 
- \frac12 \on{Tr}(L^f) \int_M (\p_\v \Ph) a^2 \on{vol}(f^*\g)-\frac12 \Delta((\p_\m\Ph).a^2) \\  
&\qquad+ (\p_\v\Ph) a \int_M \on{Tr}(L^f).a \on{vol}(f^*\g)   
- \frac12 (  \p_\m\Ph)\on{Tr}((L^f)^2) a^2  + (\p_\m\Ph) a \Delta a  \big]. 
\end{align*}
We rewrite this equation by expanding Laplacians of products,  
$$\De(a.b) = (\De a)b-2 \Tr^g(da\otimes db) + a (\De b) \quad \text{for $a,b \in C^\infty(M)$}.$$
\subsection{Theorem}{\em 
The geodesic equation of the almost local metric $G^\Ph$ on $B_i$ reads as 
\begin{equation*}\boxed{\begin{aligned}
f_t&= a.\nu^f,\\
a_t&=\frac{1}{\Ph} \bigg[\frac12 \Ph a^2 \on{Tr}(L^f) 
- \frac12 \on{Tr}(L^f) \int_M (\p_\v \Ph) a^2 \on{vol}(f^*\g)
-\frac12 a^2 \Delta(\p_\m\Ph) \\ &\qquad
+ 2a \Tr^g(d(\p_\m\Ph)\otimes da) 
+ (\p_\m\Ph) \Tr^g(da\otimes da) \\  
&\qquad+ (\p_\v\Ph) a \int_M \on{Tr}(L^f).a \on{vol}(f^*\g)   
- \frac12 (  \p_\m\Ph)\on{Tr}((L^f)^2) a^2  \bigg]. 
\end{aligned}}\end{equation*}
}

For the case of curves immersed in $\R^2$, this formula specializes to the formula given in 
\cite[section~3.4]{Michor107}. (When verifying this, remember that $\Delta=-D_s^2$ in 
the notation of \cite{Michor107}.) 

\medskip
\section{Sectional curvature  on shape space}\label{sectional_curvature}{\hfill}\par
\smallskip
To compute the sectional curvature  we will use the following formula, which is valid in a chart 
at the center 0 of the chart:
\begin{equation*}\begin{aligned}
&R_0(a_1,a_2,a_1,a_2)=G^\Ph_0(R_0(a_1,a_2)a_1,a_2)\\
&\qquad=\frac12 d^2G^\Ph_0(a_1,a_1)(a_2,a_2)-d^2G^\Ph_0(a_1,a_2)(a_1,a_2)+\frac12 d^2G^\Ph_0(a_2,a_2)(a_1,a_1)\\
&\qquad+G^\Ph_0(\Gamma_0(a_1,a_1),\Gamma_0(a_2,a_2))-G^\Ph_0(\Gamma_0(a_1,a_2),\Gamma_0(a_1,a_2)).
\end{aligned}\end{equation*}
Sectional curvature is given by $$R_0(a_1,a_2,a_2,a_1)=-R_0(a_1,a_2,a_1,a_2).$$
Therefore, we have to calculate the  metric in a chart, calculate its second derivative, and the 
value $\Gamma_0$ at the center 0 of the Christoffel symbols $\Gamma$.
\subsection{The almost local metric $G^\Ph$ in a chart}\label{sectional_curvature:metric_in_chart}

In the following section we will follow the method of \cite{Michor102}.
First we will construct a local chart for $B_i$.
Let $f_0 :M\rightarrow \mathbb R^n$ be a fixed immersion, which will be the
center of our chart. Consider the mapping 
\begin{align*}
\psi&= \psi_{f_0} : C^{\infty}(M,(-\epsilon,\epsilon)) \rightarrow \on{Imm}(M,\mathbb R^n),\\
\psi(a)(x)&=\on{exp}^{\g}_{f_0(x)}(a(x).\nu^{f_0}(x))=f_0(x)+a(x).\nu^{f_0}(x),
\end{align*}
where $\on{exp}^{\g}$ is the exponential mapping on $(\R^n,\g)$ and
where $\epsilon$ is so small that $\psi(a)$ is an immersion for each $a$.

Denote by $\pi$ the projection from $\on{Imm}(M,\mathbb R^n)$ to $B_i(M,\mathbb
R^n)$. The inverse on its image of $\pi \circ \psi :
C^{\infty}(M,(-\epsilon,\epsilon))  \rightarrow B_i(M,\mathbb{R}^n)$ is then a
smooth chart on $B_i(M, \mathbb{R}^n)$. We want to calculate the induced metric
in this chart, i.e., $$((\pi \circ \psi)^*G^\Ph)_{a}(b_1,b_2)$$ for any
$a\in C^{\infty}(M,(-\epsilon,\epsilon))$ and $b_1,b_2 \in C^{\infty}(M)$. We shall
fix the function $a$ and work with the ray of points $t.a$ in this chart.
Everything will revolve around the map:
$$f(t,x)=\psi(t.a)(x)=f_0(x)+t.a(x).\nu^{f_0}(x).$$ We shall use a fixed chart $(u,U)$
on $M$ with $\p_i=\frac{\p}{\p u^i}$. 
To calculate the metric $G^\Ph$ in this
chart we have to understand how 
$$T_{t.a} \psi.b_1 = b_1(x).\nu^{f_0}(x)$$ 
splits into  tangential and
horizontal parts with respect to the immersion $f(t,~~)$. The
tangential part locally has the form
\begin{eqnarray*}
Tf.(T_{(t.a)}\psi .(b_1))^\top= \sum_{i=1}^{n-1} c^i \p_if(t,x),        
\end{eqnarray*}
where the coefficients $c^i$ are given by
\begin{eqnarray*}
c^i=\sum_{j=1}^{n-1}g^{ij}\g\left(  b_1(x)\nu^{f_0}(x), \p_j f(t,x) \right). 
\end{eqnarray*}
Thus the horizontal part is 
$$(T_{t.a}\psi .b_1).\nu=
(T_{t.a}\psi .b_1)-Tf.(T_{(t.a)}\psi .(b_1))^\top=
b_1(x)\nu^{f_0}(x) - \sum_{i=1}^{n-1}c^i \p_i f(t,x).$$
\subsection*{Lemma}{\em
Using the local  expression of section~\ref{no} the metric $G^\Ph$ in the chart $(\pi\circ \psi)\i$ reads as
(by an abuse of notation)
\begin{align*}
\big((\pi\circ\psi_{f_0})^*G^\Ph\big)&_{(t.a)}(b_1,b_2)=G^\Ph_{\pi(\psi(t.a))}\big(T_{t.a}(\pi\circ\psi).b_1,T_{t.a}(\pi\circ\psi)b_2 \big)\\
 &= G^\Ph_{\psi(t.a)}\big((T_{t.a}\psi.b_1)^\bot.\nu,(T_{t.a}\psi.b_2)^\bot.\nu \big)\\
 &= \int_M \Ph \g\left((T_{t.a}\psi.b_1)^\bot.\nu,(T_{t.a}\psi.b_2) \right) \on{vol}(g)\\
 &= \int_M \Ph \left(b_1.b_2- \sum_{i=1}^{n-1}c^i \g( \p_i f(t,x), b_2(x).\nu^{f_0}(x) )\right) \on{vol}(g).
\end{align*}
}
\subsection{Second derivative of the $G^\Ph$-metric in the chart}\label{sectional_curvature:2nd_derivative}

We will calculate
$$\p^2_t|_0 ((\pi\circ\psi_{f_0})^*G^\Ph)_{(t.a)}(b_1,b_2).$$
We will use the following arguments repeatedly:
\begin{align*}
&\p_t|_0\p_j f=\p_j\p_t|_0 f=\p_j(a.\nu^{f_0})=(\p_ja)\nu^{f_0} +a\underbrace{(\p_j \nu^{f_0})}_{\text{tang.}},\\
&\g\left(  b_1(x)\nu^{f_0}(x), \p_j f(t,x) \right) |_{t=0}=0,
\end{align*}
and consequently $c_i|_{t=0} = 0$. 
\begin{align*}
\p_t|_0 c_i 
&=\sum_{j=1}^{n-1}\p_t|_0(g^{ij}) .0
+ \sum_{j=1}^{n-1}g^{ij} \p_t|_0\g\left(  b_1 \nu^{f_0}, \p_j f \right) \\
&=\sum_{j=1}^{n-1}g^{ij}\g\left(  b_1\nu^{f_0}, \p_t|_0\p_j f \right) 
= \sum_{j=1}^{n-1}g^{ij}\g\left(  b_1\nu^{f_0}, \p_j (a.\nu^{f_0}) \right) 
= \sum_{j=1}^{n-1}g^{ij}b_1 \p_j a. 
\end{align*}
Therefore
\begin{align*} &
\Big(b_1.b_2- \sum_{i=1}^{n-1}c^i \g( \p_i f, b_2.\nu^{f_0} )\Big)\big|_{t=0} = b_1.b_2 \\&
\p_t|_0 \Big(b_1.b_2- \sum_{i=1}^{n-1}c^i \g( \p_i f, b_2.\nu^{f_0} )\Big)  \\&\quad=
- \sum_{i=1}^{n-1} (\p_t|_0 c^i).0 - \sum_{i=1}^{n-1}0.\p_t|_0 \g( \p_i f, b_2.\nu^{f_0} )=0,\\&
\p_t^2|_0  \Big(b_1.b_2- \sum_{i=1}^{n-1}c^i \g( \p_i f, b_2.\nu^{f_0} )\Big) = 
\\ &\quad= 
- \sum_{i=1}^{n-1}(\p_t^2|_0 c^i).0  
-2\sum_{i=1}^{n-1}(\p_t|_0 c^i) \p_t|_0 \g( \p_i f, b_2.\nu^{f_0} ) 
- \sum_{i=1}^{n-1}0.\p_t^2|_0 \g( \p_i f, b_2.\nu^{f_0} ) \\ &\quad= 
-2\sum_{i=1}^{n-1}(\p_t|_0 c^i) \g( \p_i (a.\nu^{f_0}), b_2.\nu^{f_0} ) =
-2\sum_{i=1}^{n-1}(\p_t|_0 c^i) (\p_i a) b_2 \\ & \quad =
-2 b_1b_2 \sum_{i=1}^{n-1}\sum_{j=1}^{n-1}g^{ij} \p_j a . \p_i a =
-2 b_1b_2 \norm{da}_{g\i}^2 .
\end{align*}
The derivatives of $\Ph$ are
\begin{align*} &\begin{aligned}
\p_t|_0 \big(\Ph \o (\on{Vol},\on{Tr}(L))\big) &= (\p_\v \Ph). (\p_t|_0 \on{Vol})+(\p_\m \Ph) .(\p_t|_0 \on{Tr}(L)) \\
\p_t^2|_0\big(\Ph \o (\on{Vol},\on{Tr}(L))\big) &=
(\p_\v^2\Ph). (\p_t|_0 \on{Vol})^2 + (\p_\m^2\Ph). (\p_t|_0 \on{Tr}(L))^2 \end{aligned} \\ & \qquad
+ 2 (\p_\v\p_\m\Ph). (\p_t|_0 \on{Vol}).(\p_t|_0 \on{Tr}(L)) 
+ (\p_\v\Ph) (\p_t^2|_0 \on{Vol}) 
+ (\p_\m\Ph) (\p_t^2|_0 \on{Tr}(L)).
\end{align*}
\subsection*{Lemma}{\em
The second derivative of the $G^\Ph$-metric in the chart $(\pi\circ\psi)\i$ is given by
\begin{equation}\label{sectional_curvature:2nd_derivative:2nd_derivative}\begin{aligned}
\p^2_t|_0 \big((\pi\circ\psi_{f_0})^*G^\Ph\big)_{(t.a)}(b_1,b_2)&=
\Big(d^2\big((\pi \o \psi_{f_0})^* G^\Ph\big)(0)(a,a)\Big)(b_1,b_2) \\ &=
\int_M \ldots\; b_1.b_2 \on{vol}(g)
\end{aligned}\end{equation}
over the following expression:
\begin{equation*}\begin{aligned}
\ldots &=
\Ph \Big(\frac{\p_t^2|_0 \on{vol}}{\on{vol}} -2 \norm{da}_{g\i}^2  \Big)
+ (\p_\v \Ph).\Big((\p_t^2|_0 \on{Vol})+2(\p_t|_0 \on{Vol})\frac{\p_t|_0 \on{vol}}{\on{vol}} \Big) \\ &\qquad
+ (\p_\m \Ph). \Big((\p_t^2|_0 \on{Tr}(L))+2(\p_t|_0 \on{Tr}(L))\frac{\p_t|_0 \on{vol}}{\on{vol}} \Big)
+ (\p_\v^2\Ph). (\p_t|_0 \on{Vol})^2 \\ &\qquad  
+ 2 (\p_\v\p_\m\Ph). (\p_t|_0 \on{Vol})(\p_t|_0 \on{Tr}(L)) + (\p_\m^2\Ph ).(\p_t|_0 \on{Tr}(L))^2. 
\end{aligned}\end{equation*}
}
\subsection{Sectional curvature on shape space}\label{sectional_curvature:sectional_curvature}

To understand the structure of the formulas for the sectional curvature tensor, we will use some facts
from linear algebra. 
\subsubsection*{Sublemma 1}{\em
Let $V=C^\infty(M)$, and let $P$ and $Q$ be bilinear and symmetric maps $V \x V \to V$. Then 
\begin{align*}
\bp(P, Q)(a_1 \wedge a_2, b_1 \wedge b_2) &:= \frac12 \big(
 P(a_1,b_1)Q(a_2,b_2)-P(a_1,b_2)Q(a_2,b_1) \\ &\qquad 
+P(a_2,b_2)Q(a_1,b_1)-P(a_2,b_1)Q(a_1,b_2) \big)
\end{align*} 
defines a symmetric, bilinear map $(V \wedge V)\otimes(V \wedge V) \to V$.
}

Also $\bp(P,Q)=\bp(Q,P)$. 
The symbol $\bp$ stands for the Young tableau encoding the symmetries; see \cite{Fulton1997}. We have
\begin{align*}
&\bp(P, Q)(a_1 \wedge a_2, a_1 \wedge a_2) \\ & \qquad = 
\frac12 P(a_1,a_1)Q(a_2,a_2)-P(a_1,a_2)Q(a_2,a_1)+\frac12 P(a_2,a_2)Q(a_1,a_1).
\end{align*}
$P$ is called positive semidefinite if for all $x \in M$ and $a \in C^\infty(M)$, $P(a,a)(x) \geq 0$.
$P$ is called negative semidefinite if $-P$ is positive semidefinite. 
We will write  $P\geq 0, P \leq 0,P\gs 0$ if $P$ is positive semidefinite, negative semidefinite, or indefinite. 
 
\subsubsection*{Sublemma 2}{\em
If $P$ and $Q$ are positive semidefinite  bilinear and symmetric maps $V \x V \to V$, 
then  $\bp(P,Q)$ also is a positive semidefinite symmetric, bilinear map.
} 
\par {\em Proof}.
To shorten notation, we will write, for instance, $P_{12}$ instead of $P(a_1,a_2)$. 
The Cauchy inequality applied to $P$ and $Q$ gives us
$$P_{12}Q_{12}\leq \sqrt{P_{11} P_{22} Q_{11} Q_{22}},$$
and therefore we have
\begin{align*}
\bp(P, Q)(a_1 \wedge a_2, a_1 \wedge a_2)&=\frac12 P_{11}Q_{11}-P_{12}Q_{12}+\frac12 P_{22}Q_{22}\\
& \geq \frac12 P_{11}Q_{22}-\sqrt{P_{11} P_{22} Q_{11} Q_{22}}+\frac12 P_{22}Q_{11}\\
&=\frac12 \Big( \sqrt{P_{11}Q_{22}}-\sqrt{P_{22} Q_{11}} \Big)^2 \geq 0. \qquad\qedhere
\end{align*}

Let $\lambda,\mu: V\to V$. Then the map $\lambda \otimes \mu: V\otimes V \to V$ is given by
$$(\lambda\otimes \mu)(a\otimes b)=\lambda(a).\mu(b),$$
where the multiplication is in $V=C^\infty(M)$. Denote by $\lambda \vee \mu$ the symmetrization of the tensor product given by
$$\lambda \vee \mu = \tfrac12(\lambda \otimes \mu + \mu \otimes \lambda).$$
We will make use of the following simplifications.
\subsubsection*{Sublemma 3}{\em
Let $\lambda,\beta,\mu,\nu,: V\to V$. Then the bilinear symmetric map $$\bp(\lambda \vee \beta, \mu \vee \nu)$$
satisfies the following properties:
\begin{align}\tag{S1} \label{sectional_curvature:sectional_curvature:simplification1}
&\bp(\lambda \vee \mu,\lambda \vee \nu)(a_1 \wedge a_2,a_1 \wedge a_2) =
-\frac14 (\lambda \otimes \mu)(a_1 \wedge a_2) . (\lambda \otimes \nu)(a_1 \wedge a_2),\\
\tag{S2}\label{sectional_curvature:sectional_curvature:simplification2}
&\bp(\lambda \vee \mu, \lambda \otimes \lambda) = 0,\\
\tag{S3}\label{sectional_curvature:sectional_curvature:simplification3}
&\bp(\lambda \otimes \lambda, \mu \vee \nu)(a_1 \wedge a_2,a_1 \wedge a_2) = 
\frac12 (\lambda \otimes \mu)(a_1 \wedge a_2).(\lambda \otimes \nu)(a_1 \wedge a_2).
\end{align}

}

{\em Proof}.
For the proof of simplification~\eqref{sectional_curvature:sectional_curvature:simplification1} we calculate
\begin{align*}
&\bp(\lambda \vee \mu, \lambda \vee \nu)(a_1 \wedge a_2,a_1 \wedge a_2) \\
&\qquad \begin{aligned} =\frac12(\lambda \otimes \mu \otimes \lambda \otimes \nu)\Big[
&a_1 \otimes a_1 \otimes a_2 \otimes a_2 + \phantom{\frac12} a_2 \otimes a_2 \otimes a_1 \otimes a_1 \\ - 
\frac12 &a_1 \otimes a_2 \otimes a_1 \otimes a_2 - \frac12  a_1 \otimes a_2 \otimes a_2 \otimes a_1 \\ - 
\frac12 &a_2 \otimes a_1 \otimes a_1 \otimes a_2 - \frac12  a_2 \otimes a_1 \otimes a_2 \otimes a_1
\Big]. \end{aligned}
\end{align*}
Using the symmetries of the quasilinear mapping $\lambda\otimes\mu\otimes\lambda\otimes\mu$, 
we can swap the first and third positions in the tensor product of the two summands in the first line. Then the expression inside the square brackets equals 
$-\tfrac12 (a_1 \wedge a_2) \otimes (a_1 \wedge a_2)$.

Since $\lambda \otimes \lambda$ vanishes when applied to elements of $V \wedge V$, 
simplification~\eqref{sectional_curvature:sectional_curvature:simplification2} is a
direct consequence of \eqref{sectional_curvature:sectional_curvature:simplification1}.

For the proof of simplification~\eqref{sectional_curvature:sectional_curvature:simplification3} we
calculate
\begin{align*}
&\bp(\lambda \otimes \lambda, \mu \vee \nu)(a_1 \wedge a_2,a_1 \wedge a_2) \\
&\qquad \begin{aligned} =\frac12 (\lambda \otimes \lambda \otimes \mu \otimes \nu)\Big[ 
&a_1 \otimes a_1 \otimes a_2 \otimes a_2 +
 a_2 \otimes a_2 \otimes a_1 \otimes a_1 \\ - 
&a_1 \otimes a_2 \otimes a_1 \otimes a_2 -  
 a_1 \otimes a_2 \otimes a_2 \otimes a_1 \Big]. \end{aligned}
\end{align*}
Using symmetries as above, we can replace the third  summand $a_1 \otimes a_2 \otimes a_1 \otimes a_2$
by $a_2 \otimes a_1 \otimes a_2 \otimes a_1$, because the first two tensor components of 
$\lambda \otimes \lambda \otimes \mu \otimes \nu$ are equal. Then, 
swapping the second and third positions in all tensor products, we get
\begin{align*}
&\bp(\lambda \otimes \lambda, \mu \otimes \nu)(a_1 \wedge a_2,a_1 \wedge a_2)  \\ 
&\qquad \begin{aligned}
=\frac12 (\lambda \otimes \mu \otimes \lambda \otimes \nu)\Big[  
&a_1 \otimes a_2 \otimes a_1 \otimes a_2 +
 a_2 \otimes a_1 \otimes a_2 \otimes a_1 \\ - 
&a_2 \otimes a_1 \otimes a_1 \otimes a_2  -
 a_1 \otimes a_2 \otimes a_2 \otimes a_1  
 \Big]. \end{aligned}
\end{align*}
The expression inside the square brackets equals $(a_1 \wedge a_2) \otimes (a_1 \wedge a_2)$.
\qquad\qedhere\par

For orthonormal $a_1,a_2$ the sectional curvature is the negative of the curvature tensor $R_0(a_1,a_2,a_1,a_2)$. 
We will use the following \emph{formula for the curvature tensor}:
\begin{equation}\label{sectional_curvature:sectional_curvature:sectional_curvature}\begin{aligned}
&R_0(a_1,a_2,a_1,a_2)=G^\Ph_0(R_0(a_1,a_2)a_1,a_2)\\
&\quad=\frac12 d^2G^\Ph_0(a_1,a_1)(a_2,a_2)-d^2G^\Ph_0(a_1,a_2)(a_1,a_2)+\frac12 d^2G^\Ph_0(a_2,a_2)(a_1,a_1)\\
&\quad+G^\Ph_0(\Gamma_0(a_1,a_1),\Gamma_0(a_2,a_2))-G^\Ph_0(\Gamma_0(a_1,a_2),\Gamma_0(a_1,a_2)).
\end{aligned}\end{equation}

Looking at formula \eqref{sectional_curvature:2nd_derivative:2nd_derivative} 
from section~\ref{sectional_curvature:2nd_derivative},
we can express the second derivative of the metric $G^\Ph$ in the chart as
\begin{align*}
&\Big(d^2 (\pi\circ\psi_{f_0})^*G^\Ph(0)(a_1,a_2)\Big)(b_1,b_2)\\
&=\int_M \Big(\Ph. P_1(a_1,a_2)+ (\p_\v \Ph) P_2(a_1,a_2)+(\p_\m\Ph)P_3(a_1,a_2)
+(\p_\v^2\Ph)P_4(a_1,a_2)\\ 
&\qquad\qquad+ (\p_\v\p_\m\Ph) P_5(a_1,a_2) + (\p_\m^2\Ph)P_6(a_1,a_2)\Big) P(b_1,b_2)\on{vol}(g),
\end{align*}
where $P(b_1,b_2)=b_1.b_2$, so $P=\on{id}\otimes\on{id}$, and 
where the $P_i$ are obtained by symmetrizing the terms in 
formula~\eqref{sectional_curvature:2nd_derivative:2nd_derivative} from section~\ref{sectional_curvature:2nd_derivative}.

For the rest of this section, we do not note the pullback via the chart anymore, 
writing $G^\Ph_0$ instead of $\big((\pi\circ\psi_{f_0})^*G^\Ph\big)(0)$, for example. 
To further shorten our notation, we write $L$ instead of $L^{f_0}$ and $g$ instead of $g_0$. 
The following terms are calculated using the variational formulas from section~\ref{variation}:

\begin{align*}
P_1(a,a)&= \frac{\p_t^2|_0 \on{vol}}{\on{vol}} -2 \norm{da}_{g\i}^2 
= a^2 \big(\on{Tr}(L)^2- \on{Tr}(L^2)\big)-\norm{da}^2_{g\i} \\
P_2(a,a)&= (\p_t^2|_0 \on{Vol})+2(\p_t|_0 \on{Vol})\frac{\p_t|_0 \on{vol}}{\on{vol}} \\
&= \int_M a^2 \big(\on{Tr}(L)^2- \on{Tr}(L^2)\big)+\int_M \norm{da}^2_{g\i} \on{vol}(g) \\
&\qquad+2 \on{Tr}(L) . a \int_M  \on{Tr}(L) . a \on{vol}(g),\\
P_3(a,a)&= (\p_t^2|_0 \on{Tr}(L))+2(\p_t|_0 \on{Tr}(L))\frac{\p_t|_0 \on{vol}}{\on{vol}} \\
&=
2a^2\on{Tr}(L^3)+4a\on{Tr}\big( L.g\i.\nabla^2 a\big)
+2\Tr(g\i(da\otimes da)L)\\&\qquad
-\norm{da}^2_{g\i}\Tr(L)
+2a\Tr^{g}\big(\nabla_{\grad a}s\big)
\\&\qquad
+2 \big(-\Delta a + a \on{Tr}(L^2)\big) (- \on{Tr}(L) . a) \\
&=2a^2\on{Tr}(L^3)+4a\on{Tr}\big( L.g\i.\nabla^2 a\big)
+2\Tr(g\i(da\otimes da)L)\\&\qquad
-\norm{da}^2_{g\i}\Tr(L)
+2a\Tr^{g}\big(\nabla_{\grad a}s\big)\\
&\qquad+2 \on{Tr}(L)  a \Delta a -2 \on{Tr}(L)\on{Tr}(L^2) . a^2,\\
P_4(a,a)&= (\p_t|_0 \on{Vol})^2 
=\Big(\int_M \on{Tr}(L) . a \on{vol}(g)\Big)^2,\\  
P_5(a,a)&= 2 (\p_t|_0 \on{Vol})(\p_t|_0 \on{Tr}(L)) 
=2 \int_M - \on{Tr}(L) . a \on{vol}(g) \big(-\Delta a + a \on{Tr}(L^2)\big)\\
&=2 \Delta a \int_M  \on{Tr}(L) . a \on{vol}(g)  
-2 \on{Tr}(L^2) a \int_M \on{Tr}(L) . a \on{vol}(g), \\
P_6(a,a)&= (\p_t|_0 \on{Tr}(L))^2 
=\big( -\Delta a + a \on{Tr}(L^2)\big)^2 \\
&=(\Delta a)^2 -2 a\Delta a \on{Tr}(L^2)+a^2 \on{Tr}(L^2)^2.
\end{align*}
Then the \emph{first part of the curvature tensor} is given by
\begin{align*}
&\frac12 d^2G^\Ph_0(a_1,a_1)(a_2,a_2)-d^2G^\Ph_0(a_1,a_2)(a_1,a_2) +\frac12 d^2G^\Ph_0(a_2,a_2)(a_1,a_1) \\
&\qquad=\int_M \big(\Ph. \bp(P_1,P)+(\p_\v \Ph) \bp(P_2,P)+(\p_\m\Ph)\bp(P_3,P) \\
&\qquad\qquad+(\p_\v^2\Ph)\bp(P_4,P)+(\p_\v\p_\m\Ph)\bp(P_5,P) + (\p_\m^2\Ph)\bp(P_6,P)\big) \on{vol}(g) \\
&\qquad\qquad \cdot(a_1 \wedge a_2,a_1 \wedge a_2).
\end{align*}

Note that $P$ is positive definite, so that $\bp(P_i,P)$ is positive semidefinite if $P_i$ is positive semidefinite.
We can always assume that $\Ph$ is positive because otherwise $G^\Ph$ would not be a Riemannian metric.

\begin{align*} 
\tag*{\boxed{\makebox[.8cm]{$P_1P$}}} &
P_1= P_1^1 + P_1^2, \\
\intertext{with} & 
P_1^1= (\on{Tr}(L)^2-\on{Tr}(L^2)) \on{id} \otimes \on{id}, \\ &
P_1^2=- \on{Tr}^g(d \otimes d). \\ 
\intertext{Applying simplification \eqref{sectional_curvature:sectional_curvature:simplification3} to $\bp(P_1^1,P)$ and $\bp(P_1^2,P)$, we get}&
\bp(P_1^1,P) = \frac12 (\on{Tr}(L)^2-\on{Tr}(L^2)) (\on{id} \otimes \on{id})^2 = 0 \\
\intertext{on $(V \wedge V) \otimes (V \wedge V)$ and}  &
\bp(P_1^2,P) = -\frac12 \on{Tr}^g \big((\on{id} \otimes \, d)^2\big), \\ &
\bp(P_1^2,P)(a_1\wedge a_2,a_1\wedge a_2) =-\frac12 \norm{a_1 da_2 - a_2 da_1}_{g\i}^2 \leq 0. \\
\intertext{Therefore, we have} &
\int_M \Ph.\bp(P_1,P)(a_1\wedge a_2,a_1\wedge a_2)\on{vol}(g)\leq 0.\\
\tag*{\boxed{\makebox[.8cm]{$P_2P$}}}
&P_2=P_2^1+P_2^2+P_2^3\\
\intertext{with}
&P_2^1=\int_M (\on{id} \otimes \on{id})(\on{Tr}(L)^2-\on{Tr}(L^2))\on{vol}(g)\\
&P_2^2=2\on{Tr}(L) (\on{id} \vee \int_M\on{Tr}(L)\on{id}\on{vol}(g))\\
&P_2^3=\int_M \on{Tr}^g( d \otimes d)\on{vol}(g)\\
\intertext{$P_2^1$ is indefinite. Applying simplification \eqref{sectional_curvature:sectional_curvature:simplification2} we get
$\bp(P_2^2,P)=0.$ $P_2^3$, and therefore also $\bp(P_2^3,P)$ is positive semidefinite. Therefore}
&\int_M (\p_\v \Ph)\bp(P_2^1,P)(a_1\wedge a_2,a_1\wedge a_2) \on{vol}(g)\gs 0, \\
&\int_M (\p_\v \Ph)\bp(P_2^3,P)(a_1\wedge a_2,a_1\wedge a_2) \on{vol}(g) \geq 0.\\
\tag*{\boxed{\makebox[.8cm]{$P_3P$}}}
&P_3=P_3^1+P_3^2+P_3^3,\\
\intertext{with}
&P_3^1=2 \id \vee \Big(\Tr(L^3)\id+2\Tr(L g\i \nabla^2(\id)) + \Tr^g(\nabla_{\grad id}s)
      \\ & \qquad\qquad\quad
      + \Tr(L) \De(\id) - \Tr(L) \Tr(L)^2 \id \Big),\\
&P_3^2=2\Tr\big(g\i (d\otimes d) L\big),\\
&P_3^3=-\Tr^g(d\otimes d)\Tr(L).
\intertext{Applying simplification \eqref{sectional_curvature:sectional_curvature:simplification2},  
we get that $\bp(P_3^1,P)$ vanishes. Furthermore, }
&\bp(P_3^2,P)(a_1\wedge a_2,a_1\wedge a_2)= a_1^2 \Tr(g\i (da_2 \otimes da_2) .L) 
      \\ &\qquad\qquad
      -2 a_1 a_2 \Tr(g\i(da_1 \otimes da_2).L)  \\ &\qquad\qquad
      +  a_2^2 \Tr(g\i (da_1\otimes da_1).L)\\
&\qquad =g_2^0\big((a_1da_2-a_2da_1)\otimes(a_1da_2-a_2da_1),s\big)\gs 0, \\
&\bp(P_3^3,P)(a_1\wedge a_2,a_1\wedge a_2)=-\frac12 \norm{a_1 da_2 - a_2 da_1}_{g\i}^2 \Tr(L) \gs 0.\\
\tag*{\boxed{\makebox[.8cm]{$P_4P$}}}
&P_4=\int_M \on{Tr}(L)\on{id}\on{vol}(g) \otimes \int_M \on{Tr}(L)\on{id}\on{vol}(g) \\
\intertext{Applying simplification \eqref{sectional_curvature:sectional_curvature:simplification3}, we get}
&\bp(P_4,P)=\frac12 \Big(\on{id} \otimes \int_M \on{Tr}(L)\on{id} \on{vol}(g)\Big)^2.\\
\intertext{ Therefore, if $\p_\v^2\Ph\geq 0$, then}
&\int_M (\p_\v^2\Ph) \bp(P_4,P)(a_1\wedge a_2,a_1\wedge a_2)\on{vol}(g)\geq 0,\\
\tag*{\boxed{\makebox[.8cm]{$P_5P$}}}
&P_5=P_5^1+P_5^2 \\
\intertext{with}
&P_5^1= 2 \Big(\Delta \vee \int_M \on{Tr}(L)\on{id}\on{vol}(g)\Big),\\
&P_5^2=-2\on{Tr}(L^2)\Big(\on{id} \vee \int_M \on{Tr}(L)a_2\on{vol}(g)\Big).\\
\intertext{Applying simplification \eqref{sectional_curvature:sectional_curvature:simplification3},  we get that $\bp(P_5^1,P)$ is the indefinite form given by}
&\bp(P_5^1,P) =\big(\on{id} \otimes \Delta\big)\otimes\big(\on{id}\otimes\int_M \on{Tr}(L)\on{id}\on{vol}(g)\big).\\
\intertext{Simplification \eqref{sectional_curvature:sectional_curvature:simplification2} gives $\bp(P_5^2,P) =0$. Therefore,}
&\int_M (\p_\v\p_\m\Ph)\bp(P_5^1,P)(a_1\wedge a_2,a_1\wedge a_2)\on{vol}(g) \gs 0. \\
\tag*{\boxed{\makebox[.8cm]{$P_6P$}}}
&P_6=P_6^1+P_6^2\\
\intertext{with}
&P_6^1=\Delta\otimes\Delta,\\
&P_6^2=\on{Tr}(L^2)^2 \on{id} \otimes \on{id},\\
&P_6^3=-2 \on{Tr}(L^2) \on{id} \vee \Delta. \\
\intertext{Applying simplification \eqref{sectional_curvature:sectional_curvature:simplification2}, 
we get that $\bp(P_6^2,P)$ and $\bp(P_6^3,P)$ vanish.
Simplification \eqref{sectional_curvature:sectional_curvature:simplification3} gives}
&\bp(P_6^1,P)=\frac12 (\on{id}\otimes\Delta)^2\\
\intertext{If $\p_\m^2\Ph\geq 0$, then we get}
&\int_M (\p_\m^2 \Ph)\bp(P_6,P)(a_1\wedge a_2,a_1\wedge a_2)\on{vol}(g)\geq 0.
\end{align*}

Now we come to the \emph{second part of the curvature tensor} $R_0(a_1,a_2,a_1,a_2)$, which is given by
\begin{align*}
&G_0(\Gamma_0(a_1,a_1),\Gamma_0(a_2,a_2))-G_0(\Gamma_0(a_1,a_2),\Gamma_0(a_1,a_2)).
\end{align*}
From the geodesic equation calculated in section~\ref{geodesic_equation_Bi}, which is given by
\begin{align*}
a_t=\Gamma_0(a,a)
&=\frac{1}{\Ph} \big[\frac12 \Ph a^2 \on{Tr}(L) 
- \frac12 \on{Tr}(L) \int_M (\p_\v \Ph) a^2 \on{vol}(g)
-\frac12 a^2 \Delta(\p_\m\Ph) \\ &\qquad
+ 2a \Tr^g(d(\p_\m\Ph)\otimes da) 
+ (\p_\m\Ph) \Tr^g(da\otimes da) \\  
&\qquad+ (\p_\v\Ph) a \int_M \on{Tr}(L).a \on{vol}(g)   
- \frac12 (  \p_\m\Ph)\on{Tr}(L^2) a^2  \big],
\end{align*}
we can extract the Christoffel symbol by symmetrization and get
\begin{align*}
\Gamma_0(a_1,a_2)&=
\frac{1}{\Ph} \sum_{i=1}^5 Q_i(a_1,a_2), 
\end{align*}
where $Q_1, \ldots, Q_5$ are the symmetrizations of the summands in the geodesic equation. 
$Q_i$ are given by
\begin{align*}
Q_1&=\frac12\Big( \Ph \on{Tr}(L)- \Delta(\p_\m\Ph) -  (\p_\m\Ph) \on{Tr}(L^2)\Big) \id \otimes \id, \\
Q_2&=-\frac12 \on{Tr}(L) \int_M (\p_\v \Ph) \id \otimes \id \on{vol}(g), \\
Q_3&=2 \id \vee \Tr^g(d(\p_\m\Ph)\otimes d),\\ 
Q_4&=(\p_\v\Ph) \id \vee \int_M \on{Tr}(L) \id \on{vol}(g), \\
Q_5&=(\p_\m\Ph) \Tr^g(d\otimes d).
\end{align*} 
Then 
\begin{align*}
&G_0(\Gamma_0(a_1,a_1),\Gamma_0(a_2,a_2))-G_0(\Gamma_0(a_1,a_2),\Gamma_0(a_1,a_2)) \\
&\qquad= \int_M \frac{1}{\Ph} \sum_{i} \bp(Q_i,Q_i)(a_1 \wedge a_2,a_1 \wedge a_2) \on{vol}(g) \\
&\qquad\qquad+\int_M \frac{2}{\Ph} \sum_{i<j} \bp(Q_i,Q_j)(a_1 \wedge a_2,a_1 \wedge a_2) \on{vol}(g).
\end{align*}

The contribution of the following terms to $R_0(a_1,a_2,a_1,a_2)$ is 
$\int_M \frac{1}{\Ph} \ldots \on{vol}(g)$ over the terms listed. 
\begin{align*} &
\tag*{\boxed{\makebox[.8cm]{$Q_1Q_1$}}} \bp(Q_1,Q_1) = 0\\
\intertext{according to simplification
\eqref{sectional_curvature:sectional_curvature:simplification2}.} &
\tag*{\boxed{\makebox[.8cm]{$Q_2Q_2$}}}
\bp(Q_2,Q_2)(a_1\wedge a_2,a_1\wedge a_2)=\frac{\on{Tr}(L)^2}{4}\Big[\\
& \quad 
\int_M (\p_\v \Ph) a_1^2 \on{vol}(g) . \int_M (\p_\v \Ph) a_2^2 \on{vol}(g) -
\big(\int_M (\p_\v \Ph) a_1a_2 \on{vol}(g)\big)^2\Big] \\
\intertext{which is positive by the Cauchy--Schwarz inequality, assuming that $\p_\v \Ph \geq 0$.} &
\tag*{\boxed{\makebox[.8cm]{$Q_3Q_3$}}}
\bp(Q_3,Q_3)(a_1\wedge a_2,a_1\wedge a_2)\\
&\qquad\qquad= - \Big(\big(\id \otimes \Tr^g(d(\p_\m\Ph)\otimes d)\big)(a_1\wedge a_2)\Big)^2 \\
&\qquad\qquad= -g\i\big(d(\p_\m\Ph),a_1da_2-a_2da_1\big)^2\leq 0\\
\intertext{according to simplification
\eqref{sectional_curvature:sectional_curvature:simplification1}.}&
\tag*{\boxed{\makebox[.8cm]{$Q_4Q_4$}}}
\bp(Q_4,Q_4)=-\frac14(\p_\v\Ph)^2 (\id\otimes\int_M\on{Tr}(L)\id
\on{vol}(g))^2\leq 0\\
\intertext{according to simplification
\eqref{sectional_curvature:sectional_curvature:simplification1}.} &
\tag*{\boxed{\makebox[.8cm]{$Q_5Q_5$}}}
\bp(Q_5,Q_5)= (\p_\m \Ph)^2 
\big( \norm{da_1}_{g\i}^2 \norm{da_2}_{g\i}^2 - g\i(da_1,da_2)^2 \big)\\
&\qquad\qquad\quad=(\p_\m \Ph)^2\norm{da_1\wedge da_2}^2_{g_2^0}\geq 0\\
\intertext{by the Cauchy--Schwarz inequality.}&
\end{align*}
The contribution of the following terms to $R_0(a_1,a_2,a_1,a_2)$ is 
$\int_M \frac{2}{\Ph} \ldots \on{vol}(g)$ over the terms listed: 
\begin{align*} &
\tag*{\boxed{\makebox[.8cm]{$Q_1Q_2$}}}
\bp(Q_1,Q_2)=-\frac14 \big(\Ph \on{Tr}(L)^2 - \Tr(L) \De(\p_\m \Ph) 
- \Tr(L)\Tr(L^2) (\p_\m \Ph)\big) .\\&\qquad\qquad\qquad\qquad.\bp\Big(\id\otimes\id,\int_M (\p_\v
\Ph)\id\otimes\id\on{vol}(g)\Big), \\
\intertext{where the second factor is $\geq 0$ assuming that $\p_\v \Ph \geq 0$.} &
\tag*{\boxed{\makebox[.8cm]{$Q_1Q_3$}}}
\bp(Q_1,Q_3)=0
\intertext{according to simplification
\eqref{sectional_curvature:sectional_curvature:simplification2}.} &
\tag*{\boxed{\makebox[.8cm]{$Q_1Q_4$}}}
\bp(Q_1,Q_4)=0\\
\intertext{according to simplification
\eqref{sectional_curvature:sectional_curvature:simplification2}.} &
\tag*{\boxed{\makebox[.8cm]{$Q_1Q_5$}}}
\bp(Q_1,Q_5)=\frac14 \big(\Ph \on{Tr}(L)(\p_\m \Ph) - (\p_\m \Ph) \De(\p_\m \Ph) 
- \Tr(L^2) (\p_\m \Ph)^2\big) .\\&\qquad\qquad\qquad\qquad.
\norm{a_1da_2-a_2da_1}_{g\i}^2 \\&
\tag*{\boxed{\makebox[.8cm]{$Q_2Q_3$}}}
\bp(Q_2,Q_3)\gs0\\&
\tag*{\boxed{\makebox[.8cm]{$Q_2Q_4$}}}
\bp(Q_2,Q_4)=-\frac12(\p_\v\Ph)\on{Tr}(L) \cdot \\ & \qquad \cdot
\bp\big(\int_M (\p_\v\Ph) \id \otimes \id \on{vol}(g), \id \vee \int_M \on{Tr}(L) \id \on{vol}(g)\big) \\ 
\intertext{This form is indefinite, but we have} &
\int_M \frac{2}{\Ph} \bp(Q_2,Q_4) \on{vol}(g) = - \bp(\tilde Q_2,\tilde Q_4),\\
\intertext{with the positive semidefinite form} &
\tilde Q_2= \int_M (\p_\v \Ph) \id \otimes \id \on{vol}(g)  \\
\intertext{and the form} &
\tilde Q_4= 
\int_M \on{Tr}(L) \frac{1}{\Ph}(\p_\v\Ph) \id \on{vol}(g) \vee \int_M \on{Tr}(L) \id \on{vol}(g),  \\
\intertext{which is positive semidefinite if $\tfrac{\p_\v\Ph}{\Ph}$ 
is a non negative constant.} &
\tag*{\boxed{\makebox[.8cm]{$Q_2Q_5$}}}
\bp(Q_2,Q_5) =
-\frac12 (\p_\m\Ph)\on{Tr}(L) \cdot \\ & \qquad \cdot
\bp\big(\int_M (\p_\v\Ph) \id \otimes \id \on{vol}(g),\Tr^g(d \otimes d) \big) \gs 0, \\
\intertext{because of the factor $(\p_\m\Ph)\on{Tr}(L)$. But the factor} &
\bp\big(\int_M (\p_\v\Ph) \id \otimes \id \on{vol}(g),\Tr^g(d \otimes d) \big) \\
\intertext{is positive definite.} &
\tag*{\boxed{\makebox[.8cm]{$Q_3Q_4$}}}
\bp(Q_3,Q_4)\gs0\\&
\tag*{\boxed{\makebox[.8cm]{$Q_3Q_5$}}}
\bp(Q_3,Q_5)(a_1\wedge a_2,a_1\wedge a_2)\\&
=(\p_\m\Ph)
\Big( a_1g\i(d(\p_\m\Ph),da_1)\norm{da_2}_{g\i}^2-
\big(a_1g\i(d(\p_\m\Ph),da_2)\\ &\quad+
a_2g\i(d(\p_\m\Ph),da_1)\big)g\i(da_1,da_2)
+a_2g\i(d(\p_\m\Ph),da_2)\norm{da_1}_{g\i}^2\Big)\\
&=(\p_\m\Ph)g_2^0\big(d(\p_\m\Ph)\otimes(a_1da_2-a_2da_1),da_1\wedge da_2\big)  \gs0,\\&
\tag*{\boxed{\makebox[.8cm]{$Q_4Q_5$}}}
\bp(Q_4,Q_5)\gs0.
\end{align*}

We are now able to compile a list of all negative, positive, and indefinite terms of the curvature $R_0(a_1,a_2,a_1,a_2)$. 
Remember that negative terms of $R_0(a_1,a_2,a_1,a_2)$ make a positive contribution to sectional curvature. 
Positive sectional curvature is connected to the vanishing of geodesic distance 
because the space wraps up on itself in tighter and tighter ways. 

\noindent
\boxed{\makebox[.8cm]{$P_4P$}} 
\boxed{\makebox[.8cm]{$P_6P$}} 
\boxed{\makebox[.8cm]{$Q_2Q_2$}} 
\boxed{\makebox[.8cm]{$Q_5Q_5$}} 
are positive, assuming  $\p_\v\Ph,\p_\v^2\Ph,\p_\m^2\Ph \geq 0$.

\noindent
\boxed{\makebox[.8cm]{$P_1P$}} 
\boxed{\makebox[.8cm]{$Q_3Q_3$}}
\boxed{\makebox[.8cm]{$Q_4Q_4$}}
\boxed{\makebox[.8cm]{$Q_1Q_2$}}
are negative, assuming  $\p_\v\Ph \geq 0$. 

\noindent
\boxed{\makebox[.8cm]{$Q_2Q_4$}}
is negative, assuming that $\tfrac{\p_\v\Ph}{\Ph}$ is a non negative constant, 
and indefinite otherwise.

\noindent
\boxed{\makebox[.8cm]{$Q_2Q_5$}}
is negative, assuming that $\Tr(L)(\p_\m\Ph)$ is positive, 
and indefinite otherwise.

\noindent
\boxed{\makebox[.8cm]{$P_2P$}} 
\boxed{\makebox[.8cm]{$P_3P$}} 
\boxed{\makebox[.8cm]{$P_5P$}} 
\boxed{\makebox[.8cm]{$Q_1Q_5$}}
\boxed{\makebox[.8cm]{$Q_2Q_3$}}  
\boxed{\makebox[.8cm]{$Q_3Q_4$}} 
\boxed{\makebox[.8cm]{$Q_3Q_5$}} 
\boxed{\makebox[.8cm]{$Q_4Q_5$}}   
are indefinite.

\medskip
\section{Geodesic distance on $B_i(M,\mathbb R^n)$}\label{geodesic_distance}{\hfill}\par
\smallskip

We will state some conditions on $\Ph$ ensuring that the almost local metric $G^\Ph$ induces
non vanishing geodesic distance on $B_i$. The proofs are based on a comparison 
between the $G^\Ph$-length of a path and its area swept out. 
In the last part we will use the vector space structure of $\R^n$  to define a Fr\'echet metric 
on shape space $B_i(M,\R^n)$. 
In section \ref{fr:fr} it is shown how this metric is related to an $L^\infty$ Finsler metric, 
and in section \ref{fr:fr2} the Fr\'echet metric is compared to almost local metrics. 

The main results are in sections~\ref{geodesic_distance:non_vanishing} and~\ref{fr:fr2}.

\subsection{Geodesic distance on $B_i$}\label{geodesic_distance:geodesic_distance}

Geodesic distance on $B_i$ is given by
$$\on{dist}^{G^\Ph}_{B_i}(f_0,f_1) := \inf_f \Len^{G^{\Ph}}_{B_i}(f),$$
where the infimum is taken over all $f :[0,1] \to B_i$ with $f(0)=f_0$ and $f(1)=f_1$.
$\Len^{G^{\Ph}}_{B_i}$ is the length of paths in $B_i$ given by
$$\Len^{G^{\Ph}}_{B_i}(f) = \int_0^1 \sqrt{G^\Ph_f(f_t,f_t)} dt \quad \text{for $f:[0,1] \to B_i$.}$$
Letting $\pi:\on{Imm} \to B_i$ denote the projection, we have 
$$\Len^{G^\Ph}_{B_i}(\pi \o f) =  \Len^{G^\Ph}_{\on{Imm}}(f) =
\int_0^1 \sqrt{G^\Ph_f(f_t,f_t)} dt \quad 
\text{for horizontal $f:[0,1]\to \on{Imm}$}.$$
By non vanishing geodesic distance on $B_i$ we mean that $\on{dist}^{G^\Ph}_{B_i}$ separates
points.

\subsection{Area swept out}\label{geodesic_distance:area_swept_out}

For $f:[0,1]\to \on{Imm}$ we have
$$(\text{area swept out by $f$})=\int_{[0,1]\x M} \on{vol}(f(\cdot,\cdot)^* \g).$$
If $f$ is horizontal, then this integral can be rewritten as
$$(\text{area swept out by $f$})=\int_0^1 \int_M \norm{f_t} \on{vol}(f(t,\cdot)^* \g) dt=:
\int_0^1 \int_M \norm{f_t} \on{vol}(g) dt.$$

\subsection{Lemma (first area swept out bound)}{\em 
\label{geodesic_distance:area_swept_out1}
For an almost local metric $G^\Ph$ satisfying 
\begin{align*}
\Ph(\v,\m) &\geq C_1 \qquad \text{for $C_1>0$} 
\end{align*}
and a horizontal path $f:[0,1] \to \on{Imm}$, we have the area swept out bound
\begin{align*}
\sqrt{C_1} \ (\text{area swept out by $f$})
\leq \max_t \sqrt{\on{Vol}\big(f(t)\big)} . \Len^{G^\Ph}_{\on{Imm}}(f).
\end{align*}
}
The proof is an adaptation of that given in \cite[section~3.4]{Michor102} for the $G^A$-metric.

\par {\em Proof}.
\begin{align*}
\Len^{G^\Ph}_{\on{Imm}}(f)&=
\int_0^1 \sqrt{G^\Ph_f(f_t,f_t)} dt \\&=
\int_0^1 \Big(\int_M \Ph \norm{f_t}^2 \on{vol}(g) \Big)^{\frac12} dt \geq
\sqrt{C_1} \int_0^1 \Big(\int_M \norm{f_t}^2 \on{vol}(g) \Big)^{\frac12} dt \\&\geq
\sqrt{C_1} \int_0^1 \Big(\int_M \on{vol}(g) \Big)^{-\frac12} \int_M 1.\norm{f_t} \on{vol}(g) dt\\&\geq
\sqrt{C_1} \min_t \Big(\int_M \on{vol}(g) \Big)^{-\frac12} \int_0^1 \int_M \norm{f_t} \on{vol}(g) dt \\&=
\sqrt{C_1} \Big(\max_t \int_M \on{vol}(g) \Big)^{-\frac12} \int_0^1 \int_M \norm{f_t} \on{vol}(g) dt.
\qquad\qedhere
\end{align*}

\subsection{Lemma (Lipschitz continuity of $\sqrt{\on{Vol}}$)}{\em 
\label{geodesic_distance:lipschitz}
For an almost local metric $G^\Ph$, the condition 
\begin{align*}
\Ph(\v,\m) \geq C_2 \m^2
\end{align*}
implies the Lipschitz continuity of the map
$$\sqrt{\on{Vol}}:(B_i,\on{dist}^{G^\Ph}_{B_i}) \to \mathbb R_{\geq 0}$$
by the inequality holding for $f_1$ and $f_2$ in $B_i$:
$$
\sqrt{\on{Vol}(f_1)}-\sqrt{\on{Vol}(f_2)} 
\leq \frac{1}{2 \sqrt{C_2}} \on{dist}^{G^\Ph}_{B_i}(f_1,f_2).
$$
}
The proof is an adaptation of that given in \cite[section~3.3]{Michor102} for the $G^A$-metric.

\par {\em Proof}.
Let $f:[0,1] \to \on{Imm}$ be a horizontal path, and let $f_t=a.\nu^f$ denote its derivative. 
Using the formula from section~\ref{variation:volume} for the variation of the volume, we get
\begin{align*}
\p_t \on{Vol}(f) &=
-\int_M \on{Tr}(L) a \on{vol}(g) \leq 
\left| \int_M \on{Tr}(L) a \on{vol}(g)\right| \\&\leq
\Big(\int_M 1^2 \on{vol}(g)\Big)^{\frac12} \Big(\int_M \on{Tr}(L)^2 a^2 \on{vol}(g)\Big)^{\frac12} \\&\leq
\sqrt{\on{Vol}(f)} \Big(\int_M \frac{\Ph}{C_2} a^2 \on{vol}(g)\Big)^{\frac12} \leq
\frac{1}{\sqrt{C_2}}\sqrt{\on{Vol}(f)} \sqrt{G^\Ph_f(f_t,f_t)}.
\end{align*}
Thus
\begin{align*}
\p_t \sqrt{\on{Vol}(f)}=\frac{\p_t \on{Vol}(f)}{2 \sqrt{\on{Vol}(f)}}\leq
\frac{1}{2\sqrt{C_2}} \sqrt{G^\Ph_f(f_t,f_t)}.
\end{align*}
By integration we get
\begin{align*}
\sqrt{\on{Vol}(f_1)}-\sqrt{\on{Vol}(f_0)} &= 
\int_0^1 \p_t \sqrt{\on{Vol}(f)}dt \\&\leq 
\int_0^1 \frac{1}{2\sqrt{C_2}} \sqrt{G^\Ph_f(f_t,f_t)} = 
\frac{1}{2\sqrt{C_2}}\ \Len^{G^\Ph}_{\on{Imm}}(f).
\end{align*}
Now take the infimum over all horizontal paths $f$ connecting $f_1$ and $f_2$. 
\qquad\qedhere\par

\subsection{Lemma (second area swept out bound)}{\em 
\label{geodesic_distance:area_swept_out2}
For an almost local metric $G^\Ph$ satisfying  
$$\Ph(\v,\m) \geq C \v \qquad \text{with $C>0$} $$
and a horizontal path $f:[0,1] \to \on{Imm}$, we get the area swept out bound
$$
\sqrt{C}\ (\text{area swept out by $f$}) \leq \Len^{G^\Ph}_{\on{Imm}}(f).
$$
}
The proof is adapted from proofs for the case of planar curves that can be found in 
\cite[section~3.7]{Michor107}, \cite[Lemma~3.2]{Shah2008}, \cite[Proposition~1]{YezziMennucci2005}
and \cite[Theorem~7.5]{YezziMennucci2004a}.

\par {\em Proof}.
\begin{align*}
\Len^{G^\Ph}_{\on{Imm}}(f)&=
\int_0^1 \sqrt{G^\Ph_f(f_t,f_t)} dt =
\int_0^1 \Big(\int_M \Ph \norm{f_t}^2 \on{vol}(g) \Big)^{\frac12} dt \\&\geq
\sqrt{C} \int_0^1 \sqrt{\on{Vol}(f)} \Big(\int_M \norm{f_t}^2 \on{vol}(g) \Big)^{\frac12} dt \geq
\sqrt{C} \int_0^1 \int_M 1.\norm{f_t} \on{vol}(g) dt\\&=
\sqrt{C} \int_{[0,1]\x M} \on{vol}(f(\cdot,\cdot)^*\g) dt=
\sqrt{C}\ (\text{area swept out by $f$}). \qquad\qedhere
\end{align*}

\subsection{Non vanishing geodesic distance}\label{geodesic_distance:non_vanishing}

Using the estimates proven above, we get the following result.

\subsection{Theorem}{\em 
At least on $B_e(M,\R^n)=\Emb(M,\R^n)/\Diff(M)$, the almost local metric $G^\Ph$ induces non vanishing geodesic distance if 
at least one of the two following conditions holds:
\begin{align}
\label{geodesic_distance:non_vanishing:condition1}
\Ph(\v,\m) &\geq C_1 + C_2 \m^2 &&\text{for $C_1,C_2 >0$,} \\
\label{geodesic_distance:non_vanishing:condition2}
\Ph(\v,\m) &\geq C \v &&\text{for $C>0$.}
\end{align}
}

\subsection{Fr\'echet distance and Finsler metric}\label{fr:fr}

The Fr\'echet distance on the shape space $B_i(M,\R^n)$ is defined as
\begin{align*}
\dist^{L^\infty}_{B_i}(F_0,F_1) = \inf_{f_0,f_1} \norm{f_0 - f_1}_{L^\infty},
\end{align*}
where the infimum is taken over all $f_0, f_1$ with $\pi(f_0)=F_0, \pi(f_1)=F_1$.
As before, $\pi$ denotes the projection $\pi:\Imm \to B_i$. 
Fixing $f_0$ and $f_1$, one has
\begin{align*}
\dist^{L^\infty}_{B_i}\big(\pi(f_0),\pi(f_1)\big) = \inf_{\varphi} \norm{f_0 \o \varphi - f_1}_{L^\infty},
\end{align*}
where the infimum is taken over all $\varphi \in \on{Diff}(M)$.
The Fr\'echet distance is related to the Finsler metric
\begin{align*}
G^\infty : T \Imm(M,\R^n) \rightarrow \R, \qquad h \mapsto \norm{h^\bot}_{L^\infty}.
\end{align*}

\subsection*{Lemma}{\em
The path length distance induced by the Finsler metric $G^\infty$ provides an upper bound for the Fr\'echet distance:
\begin{align*}
\dist^{L^\infty}_{B_i}(F_0,F_1) \leq \dist^{G^\infty}_{B_i}(F_0,F_1) = \inf_f \int_0^1 \norm{f_t}_{G^\infty} dt,
\end{align*}
where the infimum is taken over all paths 
$$f:[0,1] \to \Imm(M,\R^n)\quad\text{with}\quad\pi(f(0))=F_0,\quad \pi(f(1))=F_1.$$
}
\par {\em Proof}.
Since any path $f$ can be reparametrized such that $f_t$ is normal to $f$, one has
$$\inf_f \int_0^1 \norm{f_t^\bot}_{L^\infty} dt = \inf_f \int_0^1 \norm{f_t}_{L^\infty} dt, $$
where the infimum is taken over the same class of paths $f$ as described above. Therefore,
\begin{align*}
\dist^{L^\infty}_{B_i}(F_0,F_1) &= \inf_f \norm{f(1)-f(0)}_{L^\infty} 
= \inf_f \norm{ \int_0^1 f_t dt}_{L^\infty} 
\leq \inf_f \int_0^1 \norm{f_t}_{L^\infty} dt \\ &
= \inf_f \int_0^1 \norm{f_t^\bot}_{L^\infty} dt
= \dist^{G^\infty}_{B_i}(F_0,F_1). \qquad\qedhere
\end{align*}

It is claimed in \cite[Theorem~13]{MennucciYezzi2008} that $\dist^{L^\infty}_{B_i}=\dist^{G^\infty}_{B_i}$. 
Unfortunately, the proof is not correct because convex combinations of immersions are used, 
even though the space of immersions is not convex.

\subsection{Theorem (almost local versus Fr\'echet distance on shape space)}{\em 
\label{fr:fr2}
On $B_i(M,\mathbb R^n)$ the  $G^\Ph$ distance cannot be bounded from below by 
the Fr\'echet distance if 
any one of the following conditions holds:
\begin{align}
\label{gd:Frechet:condition1}
\Ph(\v,\m) &\leq C_1 + C_2 \m^{2k} &&\text{for $C_1,C_2 >0$ and $k< (\dim(M)+2)/2$,} \\
\label{gd:Frechet:condition2}
\Ph(\v,\m) &\leq C \v^k &&\text{for $C>0$,}\\
\Ph(\v,\m) &\leq C\on e^\v &&\text{for } C> 0.\label{gd:Frechet:condition3}
\end{align}
Indeed, then the identity map
$$\Id: \big(B_i(M,\R^n),\dist_{B_i}^{G^\Ph}\big) \to \big(B_i(M,\R^n),\dist_{B_i}^{G^\infty}\big)$$ 
is not continuous.
}

\par {\em Proof}.
Let $f_0$ be a fixed immersion of $M$ into $\R^n$, and let $f_1$ be a translation of $f_0$ by a vector $h$ of length $\ell$. 
We will show that the $G^{\Ph}$-distance between $\pi(f_0)$ and $\pi(f_1)$ is bounded by a constant that does 
not depend on $\ell$. It follows that the $G^{\Ph}$-distance cannot be bounded from below by the Fr\'echet distance, and this proves the claim. 

For small $r_0$, we calculate the $G^\Ph$-length of the following path of immersions:
First scale $f_0$ down to a factor $r_0$, then translate it by a vector $h$ of length $\ell$, and then
scale it up again around the new origin $h$ until it has reached $f_1$. Let $m=\dim(M)$.

For the scaling down part, let $r$ be a decreasing function such that $r(0)=1$ and $r(1)=r_0$. 
Then the length of the path $f(t):=r(t).f_0$ is 
\begin{align*}
&\Len_{\Imm}^{G^\Ph}(f)
\\&=\int_0^1\sqrt{\int_M\Ph\big(\Vol(r.f_0),\Tr(L^{r.f_0})\big) 
\g\big(r_t.f_0,r_t.f_0\big)\vol\big((r.f_0)^*\g\big)}dt\\
&=\int_0^1\sqrt{\int_Mr_t^2.\Ph\big(r^m\Vol(f_0),\frac{1}{r}\Tr(L^{f_0})\big) 
\g\big(f_0,f_0\big)r^m\vol\big(f_0^*\g\big)}dt\\
&=\int_{r_0}^1\sqrt{\int_M \Ph\big(r^m\Vol(f_0),\frac{1}{r}\Tr(L^{f_0})\big) 
\g\big(f_0,f_0\big)r^m\vol\big(f_0^*\g\big)}dr.
\end{align*}
The last integral converges for $r_0 \to 0$ under any of the above assumptions. 
So we see that the length of the shrinking part is bounded by
a constant that does not depend on $\ell$. 

The path $f(t):=r_0.f_0+t.h$ is a pure translation of the scaled immersion $r_0.f_0$ by the vector $h$ of length
$\ell$. The length of this path is
\begin{align*}
\Len^{G^\Ph}_{\Imm}(f)&=\int_0^1\int_M \Ph.\g(f_t,f_t)\vol(g)dt
\\&=\int_0^1\int_M \Ph.\ell^2\vol(g)dt= \ell^2\int_M\Ph\vol(g)\\&=
\begin{cases}
  O(r_0^{(m-2k)}) 
 & \text{if $\Ph$ satisfies \eqref{gd:Frechet:condition1}},\vspace{0.1cm}\\
 O(r_0^{m(k+1)})
 & \text{if $\Ph$ satisfies \eqref{gd:Frechet:condition2}},\vspace{0.1cm}\\
 O(e^{r_0.m}.r_0^{m})
 & \text{if $\Ph$ satisfies \eqref{gd:Frechet:condition3}}.
\end{cases}
\end{align*}
Under the above assumptions, this tends  to zero as $r_0$ tends to zero.

To scale the immersion $r_0.f_0+h$ back up to its original size, we use the path
$f(t):=r(1-t).f_0+h$ with $r(t)$ as in the shrinking part. It follows as before that the length of 
this path is bounded by a constant that does not depend on $\ell$.

Finally, we use
\begin{equation*}
\dist_{B_i}^{G^\Ph}\big(\pi(f_0),\pi(f_1)\big) \leq \dist_{\Imm}^{G^\Ph}(f_0,f_1). \qquad\qedhere
\end{equation*}

\medskip
\section{The set of concentric spheres}\label{con_spheres}{\hfill}\par
\smallskip

For an almost local metric, the set of spheres with common center $x \in \R^n$
is a totally geodesic subspace of $B_i$. The reason is that it is the fixed point set 
of a group of isometries acting on $B_i$, namely, the group of rotations of $\R^n$ around $x$. 
(We also have to assume uniqueness of solutions to the geodesic equation.)
For the $G^A$ metric where $\Ph=1+A\Tr(L)^2$ and plane curves, 
the set of concentric spheres has been studied in \cite{Michor98},
and  for Sobolev type metrics it has been studied in \cite{Michor119,Harms2010}.
Some work for the $G^0$-metric has also been done by \cite{Salvai2009}. 

We denote the $(n-1)$-dimensional volume of the $n-1$-dimensional unit sphere in $\R^n$ by
$$\om_{n-1} :=\frac{n\pi^{\frac n2}}{\Ga(1+\frac n2)} .$$

\subsection{Theorem}{\em 
Within a set of concentric spheres, any sphere is uniquely described by its radius $r$. 
Thus the geodesic equation within a set of concentric spheres reduces to an ordinary differential equation 
for the radius. It is given by
\begin{equation*}\boxed{\begin{aligned}
r_{tt}&=-r_t^2 \frac{n-1}{\Ph} \Big( \frac{1}{2r} \Ph 
+\frac{\p_\v \Ph}2 \,r^{n-2} \om_{n-1}  + \frac{1}{2 r^2} (\p_\m\Ph)  \Big). 
\end{aligned}}\end{equation*}
The space of concentric spheres is geodesically complete with respect to a $G^\Ph$ metric if and only if
\begin{align*}
\int_{0}^{r_1}r^{\frac{n-1}{2}} 
 \sqrt{\Ph\big(\om_{n-1} r^{n-1},-(n-1)/r\big)}   dr=\infty, \qquad r_1>0,
\intertext{and}
\int_{r_0}^{\infty}r^{\frac{n-1}{2}} 
 \sqrt{\Ph\big(\om_{n-1} r^{n-1},-(n-1)/r\big)}   dr=\infty, \qquad r_0>0.
\end{align*}
For the metrics studied in this work, this yields
\begin{align*}
&\Ph(\v,\m)=\v^k: &&\text{incomplete, } 
\\
&\Ph(\v,\m)=e^\v: && \text{incomplete, } 
\\
&\Ph(\v,\m)=1+A\m^{2k}: && \text{complete  iff } k \geq \frac{n+1}{2},
\\
&\Ph(\v,\m)=\v^\frac{1+n}{1-n}+A \frac{\m^2}{\v}: &&\text{complete.}
\end{align*}
}

\par {\em Proof}.
The differential equation for the radius  can be read from the geodesic equation in 
section~\ref{geodesic_equation_Bi:geodesic_equation} when 
it is taken into account that all functions are constant on each sphere and that
$$\Vol=\om_{n-1} r^{n-1}, \quad L=-\tfrac{1}{r}\on{Id}_{TM}, \quad \Tr(L^k)=(-1)^k\tfrac{n-1}{r^k}.$$
To determine whether the space of concentric spheres is complete, we calculate the length of a path $f$ connecting a sphere 
with radius $r_0$ to a sphere with radius $r_1$:
\begin{align*}
\Len^{G^{\Ph}}_{B_i}(f) &= \int_0^1 \sqrt{G^\Ph_f(f_t^{\bot},f_t^{\bot})} dt\\&=
\int_0^1\sqrt{\int_M \Ph\big(\Vol,\Tr(L)\big) r^2_t \vol(g)}  dt\\&=
\int_0^1 |r_t| \sqrt{\Ph\big(\om_{n-1} r^{n-1},-(n-1)/r\big) \om_{n-1} r^{n-1}}   
dt\\&=
\sqrt{\om_{n-1}}\int_{r_0}^{r_1}r^{\frac{n-1}{2}} 
 \sqrt{\Ph\big(\om_{n-1} r^{n-1},-(n-1)/r\big)}   dr.\qquad\qedhere
\end{align*}

\medskip
\section{Special cases of almost local metrics}\label{special_metrics}{\hfill}\par
\smallskip

\subsection{The $G^0$-metric}\label{special_metrics:g0_metric}

The $G^0$-metric is the special case of a $G^\Ph$-metric with $\Ph \equiv 1$. 
Thus its geodesic equation can be read  from section~\ref{geodesic_equation_imm:geodesic_equation}.
It reads as
$$\boxed{\begin{aligned}
f_t &=a.\nu+Tf.f_t^\top,\\
f_{tt}&=-\frac{1}{2}\left(
\|f_t\|^2 \on{Tr}(L) . \nu + Tf.\on{grad}^{g}(\|f_t\|^2) \right) +
\left(\on{Tr}(L) . a -\on{div}^{ g}(hf_t^\top) \right) . f_t.
\end{aligned}}$$
We have three conserved quantities, namely,
$$\boxed{\begin{aligned}
  g(f_t^\top) \on{vol}({g})\in\Ga(T^*M\otimes_M\on{vol}(M)) && 
\text{reparametrization momentum,}\\
 \int_M f_t \on{vol}({g}) && \text{linear momentum,}\\
\int_M (f\wedge f_t) \on{vol}({g})\in {\textstyle\bigwedge^2}\mathbb R^n\cong \mathfrak{so}(n)^* && 
\text{angular momentum.} \end{aligned}}$$
The geodesic equation on $B_i(M,\mathbb R^n)$ is well studied. 
We can read it from section~\ref{geodesic_equation_Bi}.
$$\boxed{\begin{aligned}
f_t=a.\nu, \qquad a_t=\frac{\on{Tr}(L).a^2}{2}.
\end{aligned}}$$
Sectional curvature is given by
$$\boxed{\begin{aligned}
R_0(a_1,a_2,a_2,a_1) &= 
\frac12 \int_M \norm{a_1 da_2 - a_2 da_1}_{g\i}^2 \on{vol}(g) \geq 0.
\end{aligned}}$$
This formula is in accordance with \cite[section~4.5]{Michor102} since 
we have codimension one and a flat ambient space, so that only $\text{term}(6)$ remains, 
and for the case of plain curves, it is in accordance with \cite[section~3.5]{Michor107}.

The $G^0$-metric induces vanishing geodesic distance; see section \ref{geodesic_distance}.

\subsection{The $G^A$-metric}\label{special_metrics:g_a_metric}

For a constant $A>0$, the $G^A$-metric is defined as
$$G^A_f(h,k)=\int_M \big(1+A \on{Tr}(L)^2\big) \g( h,k ) \on{vol}({g}).$$
This metric has been introduced by \cite{Michor98,Michor102,Michor107}.
It corresponds to an almost local metric $G^\Ph$ with $\Ph(\v,\m)=(1+A\m^2)$; thus its geodesic equation on $\on{Imm}(M,\mathbb R^n)$ is given 
by (see section~\ref{geodesic_equation_imm:geodesic_equation})
$$\boxed{\begin{aligned}
f_t &=a.\nu+Tf.f_t^\top,   \\
f_{tt}&=\frac{1}{2}\Big[-\frac{\Delta\big((2A\on{Tr}(L)) \norm{f_t}^2 \big)}{1+A\on{Tr}(L)^2}
 +\norm{f_t}^2.\on{Tr}(L)\big(\frac{2A\on{Tr}(L^2)}{1+A\on{Tr}(L)^2}-1\big) \Big] \nu\\
&\qquad +\frac{Tf.\Big[(2A\on{Tr}(L))\norm{f_t}^2\on{grad}^{ g}(\on{Tr}(L))- \on{grad}^{ g}((1+A\on{Tr}(L)^2)\norm{f_t}^2)\Big]}{2(1+A\on{Tr}(L)^2)}\\
&\qquad-\Big[  \frac{(2A\on{Tr}(L))}{1+A\on{Tr}(L)^2} \big(-\Delta a+a \on{Tr}(L^2)+d\on{Tr}(L)({f_t^\top})\big)\\
&\qquad+ \on{div}^{g}(f_t^\top)- \on{Tr}(L).a \Big] f_t.
\end{aligned}}$$
The conserved quantities have the form
$$\boxed{\begin{aligned}
(1+A\on{Tr}(L)^2)\, g(f_t^\top) \on{vol}({g}) \in \Ga(T^*M\otimes_M\on{vol}(M)) && 
\text{reparam. momentum,}\\
 \int_M (1+A\on{Tr}(L)^2) f_t \on{vol}({g}) && \text{linear momentum,}\\
\int_M (1+A\on{Tr}(L)^2) 
(f\wedge f_t)\on{vol}({g})\in {\textstyle\bigwedge^2}\mathbb R^n\cong \mathfrak{so}(n)^* && 
\text{angular momentum.} 
\end{aligned}}$$
The horizontal geodesic equation for the $G^A$--metric reduces to
$$\boxed{
\begin{aligned}
f_t &= a.\nu,\\
a_t&=\frac12 a^2 \on{Tr}(L)+\frac{-a^2A\Delta(\Tr(L))+4Aag\i(d\Tr(L),da)}{(1+A\on{Tr}(L)^2)}\\
&\qquad+\frac{2A\Tr(L)\norm{da}_{g\i}^2-A\Tr(L)\Tr(L^2)a^2}{(1+A\on{Tr}(L)^2).}
\end{aligned}}$$
For the case of curves immersed in $\R^2$, this formula specializes to the formula given in 
\cite[section~4.2]{Michor98}. (When verifying this, remember that $\Delta=-D_s^2$ in 
the notation of \cite{Michor98}.) 

The curvature tensor $R_0(a_1,a_2,a_1,a_2)$ is the sum of

\noindent
\boxed{\makebox[.8cm]{$P_1P$}} \boxed{\makebox[.8cm]{$Q_3Q_3$}} negative terms,

\noindent
\boxed{\makebox[.8cm]{$P_6P$}} \boxed{\makebox[.8cm]{$Q_5Q_5$}} positive terms, and

\noindent
\boxed{\makebox[.8cm]{$P_3P$}} \boxed{\makebox[.8cm]{$Q_1Q_5$}}  \boxed{\makebox[.8cm]{$Q_3Q_5$}}
indefinite terms.
\noindent
\noindent
$$\boxed{\begin{aligned}
&R_0(a_1,a_2,a_1,a_2)=\int_M A(a_1\Delta a_2-a_2\Delta a_1)^2 \vol(g)\\&\quad
+ \int_M 2A \Tr(L) g^0_2\big((a_1da_2-a_2da_1)\otimes(a_1da_2-a_2da_1),s\big) \vol(g) \\&\quad
+\int_M \frac{1}{1+A\Tr(L)^2}\bigg[-4A^2 g\i\big(d\Tr(L), a_1da_2-a_2da_1\big)^2\\
\\&\quad-
\Big(\frac{1}{2}\big(1+A\Tr(L)^2\big)^2+ 2A^2\Tr(L) \De(\Tr(L)) 
+2A^2\Tr(L^2) \Tr(L)^2\Big) \cdot\\&\qquad\quad\cdot
\norm{a_1da_2-a_2da_1}_{g\i}^2 +(2A^2\Tr(L)^2)
\norm{da_1\wedge da_2}_{g_0^2}^2 \\
&\quad+(8A^2\Tr(L))g_2^0\big(d\Tr(L)\otimes(a_1da_2-a_2da_1),da_1\wedge da_2\big)\bigg]\on{vol}({g}).
\end{aligned}}$$
We want to express the curvature in terms of the basic skew symmetric forms.
Therefore, mimicking the notation of \cite{Michor98,Michor107}, we define 
\begin{align*}
W_2&=a_1da_2-a_2da_1,\quad W_{22}=a_1\De a_2-a_2\De a_1,\quad W_{12}=da_1\wedge da_2.
\end{align*}
Then the above equation reads as
$$\boxed{\begin{aligned}
&R_0(a_1,a_2,a_1,a_2)=\int_M A W_{22}^2 \vol(g)
+ \int_M 2A \Tr(L) g^0_2\big(W_2\otimes W_2,s\big) \vol(g) \\&\quad
+\int_M \frac{1}{1+A\Tr(L)^2}\bigg[-4A^2 g\i\big(d\Tr(L),W_2\big)^2\\
\\&\quad-
\Big(\frac{1}{2}\big(1+A\Tr(L)^2\big)^2+ 2A^2\Tr(L) \De(\Tr(L)) 
+2A^2\Tr(L^2) \Tr(L)^2\Big) \norm{W_2}_{g\i}^2\\&\quad
 +(2A^2\Tr(L)^2)
\norm{W_{12}}_{g_0^2}^2 +(8A^2\Tr(L))g_2^0\big(d\Tr(L)\otimes W_2,W_{12}\big)\bigg]\on{vol}({g}).
\end{aligned}}$$

For the case of plain curves, this formula specializes to the formula given in 
\cite[section~3.6]{Michor107}.

The $G^A$-metric satisfies condition \eqref{geodesic_distance:non_vanishing:condition1} from 
section~\ref{geodesic_distance:non_vanishing}; thus it induces non-vanish\-ing geodesic distance.

\subsection{Conformal metrics} \label{special_metrics:conformal_metric}

The conformal metrics correspond to  almost local metrics $G^\Ph$ where $\Ph$ depends only on the volume 
and not on the mean curvature. 
For the case of planar curves these metrics have been treated in \cite{YezziMennucci2004,YezziMennucci2004a,YezziMennucci2005,Shah2008}.
Then \cite{Shah2008} provides very interesting estimates on geodesic distance induced by metrics with 
$\Ph(\v)=\v$ and $\Ph(\v)=e^{\v}$. 
The geodesic equation on $\on{Imm}(M,\mathbb R^n)$ is given by
$$\boxed{\begin{aligned}
f_t &= h=a.\nu+Tf.h^\top,   \\
h_t&=-\frac{1}{2}\Big[\frac{\Ph'}{\Ph}\left(\int_M \|h\|^2 \on{vol}(g) \right) \on{Tr}(L).\nu \\ 
&\qquad+\|h\|^2 \on{Tr}(L) . \nu + Tf.\on{grad}^{g}(\|h\|^2) \Big] \\
&\qquad+\left[ \frac{\Ph'}{\Ph} \left(\int_M \on{Tr}(L).a \on{vol}({g})\right)
+\on{Tr}(L) . a -\on{div}^{g}(h^\top) \right] . h.
\end{aligned}}$$
The conserved quantities are given by
$$\boxed{
\begin{aligned}
\Ph(\on{Vol}) g(f_t^\top)\on{vol}({g})\in\Ga(T^*M\otimes_M\on{vol}(M)) && 
\text{reparam. momentum,}\\
 \Ph(\on{Vol})\int_M f_t \on{vol}(g) && \text{linear momentum,}\\
\Ph(\on{Vol})\int_M (f\wedge f_t)\on{vol}({g})\in\textstyle{\bigwedge^2}\mathbb R^n 
\cong \mathfrak{so}(n)^*&& \text{angular momentum.} \end{aligned}}$$
The horizontal part of the geodesic equation is given by
\begin{align*}
a_t&=\g\Big( \frac{1}{2} H(a.\nu,a.\nu)-K(a.\nu,a.\nu),\nu \Big) \\
&= -\frac{\Ph'}{2\Ph} \left( \int_M
a^2 \on{vol}(g)\right) \on{Tr}(L)+\frac{1}{2}a^2 \on{Tr}(L) 
+\frac{\Ph'}{\Ph}\left(\int_M a.\on{Tr}(L)
\on{vol}(g)\right) a.
\end{align*}
To simplify this equation let $b(t)=\Ph(\on{Vol}).a(t)$. We get
\begin{align*}
b_t&=\Ph'.(D_{(f,a.\nu)}\on{Vol}).a+\Ph.a_t\\
&=-\Ph'.a.\int_M \on{Tr}(L). a.\on{vol}(g)
+\Ph\frac{1}{2}a^2.\on{Tr}(L) \\
&\qquad -\frac{1}{2}\Ph'\left(\int_M a^2 \on{vol}(g)\right)
        .\on{Tr}(L)+\Ph'.a.\int_M \on{Tr}(L) .a \on{vol}(g) \\
&=-\frac{1}{2}\Ph'\int_M a^2\on{vol}(g).\on{Tr}(L)+\frac{1}{2} \Ph a^2.\on{Tr}(L).
\end{align*}
Thus the geodesic equation of the conformal metric $G^\Ph$ on $B_i$ is
$$\boxed{
\begin{aligned}
f_t&=\frac{b(t)}{\Ph(\on{Vol})} \nu, \\
b_t&=\frac{\on{Tr}(L)}{2
\Ph(\on{Vol})}\left(b^2-\frac{\Ph'(\on{Vol})}{\Ph(\on{Vol})}
\int_M b^2 \on{vol}({g})\right).
\end{aligned}}$$
For the case of curves immersed in $\R^2$, this formula specializes to the formula given in 
\cite[section~3.7]{Michor107}. 

Assuming that $\Ph'$ and $\Ph''$ are non negative, the curvature tensor  consists of the following summands. 

\noindent
\boxed{\makebox[.8cm]{$P_4P$}} 
\boxed{\makebox[.8cm]{$Q_2Q_2$}} 
are the positive summands.

\noindent
\boxed{\makebox[.8cm]{$P_1P$}} 
\boxed{\makebox[.8cm]{$Q_4Q_4$}} 
\boxed{\makebox[.8cm]{$Q_1Q_2$}} 
are the negative summands. 

\noindent
\boxed{\makebox[.8cm]{$Q_2Q_4$}} 
is indefinite, but assuming that $\tfrac{\Ph'}{\Ph}$ is a non negative constant, it is negative. 
Solving the ODE $\tfrac{\Ph'}{\Ph}=C>0$ leads to $\Ph(\on{Vol})= e^{C.\on{Vol}}$. 
In the case of curves, conformal metrics of this type have been studied in 
\cite{MennucciYezzi2008} and \cite{Shah2008}.

\noindent
\boxed{\makebox[.8cm]{$P_2P$}} 
is indefinite.

Since the formula for sectional curvature with general $\Ph=\Ph(\on{Vol})$ is still too long, 
we will  print only the formula for $\Ph(\on{Vol})=\on{Vol}$. 
To shorten notation we will write $\ol{a}$ 
for the integral over $a \in C^\infty(M)$, i.e., $$\ol{a}=\int_M a\on{vol}(g).$$
Then the sectional curvature reads as
$$\boxed{
\begin{aligned}
&R_0(a_1,a_2,a_1,a_2)= -\frac12 \on{Vol}\int_M\norm{a_1 da_2 - a_2 da_1}_{g\i}^2\on{vol}(g)\\
&\quad+\frac{1}{4\on{Vol}}\ol{\on{Tr}(L)^2}\Big(\overline{a_1^2}.\ol{a_2^2}-\ol{a_1.a_2}^2\Big)\\
&\quad+\frac{1}{4}\Big(\overline{a_1^2}.\ol{\on{Tr}(L)^2 a_2^2}-2\ol{a_1.a_2}.\ol{\on{Tr}(L)^2 a_1.a_2}+\ol{a_2^2}.\ol{\on{Tr}(L)^2 a_1^2}\Big)\\
&\quad- \frac{3}{4\on{Vol}}\Big(\ol{a_1^2}.\ol{\on{Tr}(L)a_2}^2-2\ol{a_1.a_2}.\ol{\on{Tr}(L)a_1}.\ol{\on{Tr}(L)a_2} +\ol{a_2^2}.\ol{\on{Tr}(L)a_1}^2\Big)\\
&\quad+\frac{1}{2}\Big(\ol{a_1^2}.\ol{\on{Tr}^g((da_2)^2)}
-2\ol{a_1.a_2}.\ol{\on{Tr}^g(da_1.da_2)}+\ol{a_2^2}\ol{\on{Tr}^g((da_1)^2)}\Big)\\
&\quad-\frac12\Big(\ol{a_1^2}.\ol{a_2^2.\on{Tr}(L^2)}-2.\ol{a_1.a_2}.\ol{a_1.a_2.\on{Tr}(L^2)}+\ol{a_2^2}.\ol{a_1^2.\on{Tr}(L^2)}\Big).
\end{aligned}}$$
For the case of curves immersed in $\R^2$, this formula is in accordance with the formula given in 
\cite[section~3.7]{Michor107}. 

From condition~\eqref{geodesic_distance:non_vanishing:condition2} in 
section~\ref{geodesic_distance:non_vanishing} we read  that the conformal metrics induce non vanishing geodesic distance if
$\Ph(\on{Vol})\geq C.\on{Vol} \text{ for some constant } C > 0.$
\subsection{A scale--invariant metric}\label{special_metrics:scale_invariant}

For a constant $A>0$ we define the  metric
$$G^{SI}_f(h,k)= \int_M \Big(\on{Vol}^{\frac{1+n}{1-n}}+A \frac{\on{Tr}(L)^2}{\on{Vol}}\Big)
\g( h,k ) \on{vol}(g).$$
This is an almost local metric with $\Ph(\v,\m)=\v^{\frac{1+n}{1-n}}+A \frac{\m^2}{\v}$.
Scale--invariance means that this metric does not change when $f,h,k$ are replaced by 
$\lambda f, \lambda h, \lambda k$ for $\lambda >0$.
To see that $G^{SI}$ is scale--invariant, we calculate as in \cite{Michor107}
how the scaling factor $\lambda$ changes the metric, volume form, and volume and mean curvature. 
We fix an oriented chart $(u^1,\ldots, u^{n-1})$ on $M$. Then
\begin{align*}
(\lambda f)^*\g(\p_i,\p_j) & =
\g( T(\lambda f).\p_i ,T(\lambda f).\p_j ) =
\lambda^2.f^*\g (\p_i,\p_j), \\
\on{vol}((\lambda.f)^*\g) & = 
\sqrt{\on{det}(\lambda^2(f^*\g)|_U)}\ du^1\wedge\ldots\wedge d^{u-1} = 
\lambda^{n-1}\on{vol}(f^*\g), \\
\on{Tr}(L((\lambda f)^*\g)) & = 
((\lambda f)^*\g)^{ij}\g(\frac{\p^2(\lambda f)}{\p_i\p_j},\nu^{\lambda.f}) \\
&= \frac{\lambda}{\lambda^2}(f^*\g)^{ij}\g(\frac{\p^2 f}{\p_i\p_j},\nu^{f}) = 
\frac{1}{\lambda}\on{Tr}(L(f)).
\end{align*}
The scale-invariance of the  metric $G^{SI}$ follows.
Thus along geodesics we have an additional conserved quantity 
(see section~\ref{geodesic_equation_imm:momentum_mappings}), namely,
$$ \boxed{
\int_M \Big(\on{Vol}^{\frac{1+n}{1-n}}+A\frac{\on{Tr}(L)^2}{\on{Vol}}\Big)
\g( f,f_t)\on{vol}(g)\qquad\qquad \text{scaling momentum.}
} $$
From \ref{geodesic_equation_Bi} we can read  the geodesic equation for $G^{SI}$ on $B_i$:
\begin{equation*}\boxed{\begin{aligned}
f_t&= a.\nu,\\
a_t&=\frac12 a^2 \on{Tr}(L)+\frac{1}{\on{Vol}^{\frac{1+n}{1-n}}+A\frac{\on{Tr}(L)^2}{\on{Vol}}} \\ & \qquad 
\Big[ - \frac12 \on{Tr}(L) \int_M \Big(\tfrac{1+n}{1-n}\on{Vol}^\frac{2n}{1-n}
-A\frac{\on{Tr}(L)^2}{\on{Vol}^2}\Big) a^2 \on{vol}(g) 
-A\frac{\Delta(\on{Tr}(L)).a^2}{\on{Vol}} \\
&\qquad+\frac{4A.a}{\Vol}g\i(d\Tr(L),da)+\frac{2A\Tr(L)}{\Vol}\norm{da}^2_{g\i}\\&\qquad
 +\Big(\tfrac{1+n}{1-n}\on{Vol}^\frac{2n}{1-n}
-A\frac{\on{Tr}(L)^2}{\on{Vol}^2}\Big) a \int_M \on{Tr}(L).a \on{vol}(g)   
 - A\frac{\on{Tr}(L^2)\on{Tr}(L)}{\on{Vol}} a^2 \Big] .
\end{aligned}}\end{equation*}
For the case of curves immersed in $\R^2$, this formula specializes to the formula given in 
\cite[section~3.8]{Michor107}. (When verifying this, remember that $\Delta=-D_s^2$ in 
the notation of \cite{Michor107}.)

The metric $G^{SI}$ induces non vanishing geodesic distance. 
This follows from the fact that $\on{log}(\on{Vol})$ is Lipschitz; see \cite[section 3.8]{Michor107}.

\medskip
\section{Numerical results}\label{numerics}{\hfill}\par
\smallskip
\subsection{Discretizing the horizontal path energy}

We want to solve the boundary value problem for geodesics in shape space of surfaces in $\R^3$
with respect to several almost local metrics -- 
more specifically, with respect to $G^\Ph$-metrics 
with $$\Ph=\Vol^k,\quad \Ph=e^{\Vol},\quad \Ph=1+A\Tr(L)^{2k}$$
and the scale-invariant 
metric $$\Ph=\Vol^{\frac{1+3}{1-3}}+A \frac{\Tr(L)^2}{\Vol}.$$
In order to solve this infinite-dimensional problem 
numerically, we will reduce it to a finite-dimensional problem by approximating an immersed surface 
by a triangular mesh. 

One approach to solving the boundary value problem is by the method of geodesic shooting. 
This method is based on iteratively solving the initial value problem for geodesics while suitably
adapting the initial conditions. 

Another approach, and the approach we will follow, is to minimize horizontal path energy
$$E^{\on{hor}}(f) = \int_0^1 \int_M \Ph(\Vol,\Tr(L)) \g(f_t,\nu)^2 \vol(g)$$
over the set of paths $f$ of immersions with fixed endpoints. 
Note that, by definition, the horizontal path
energy does not depend on reparametrizations of the surface. Nevertheless we want
the triangular mesh to stay regular. This can be achieved by adding a penalty functional to the horizontal path energy. 

\subsection{Discrete path energy}\label{numerics:discrete}

To discretize the horizontal path energy
$$E^{\on{hor}}(f) = \int_0^1 \int_M \Ph(\Vol,\Tr(L)) \g(f_t,\nu)^2 \vol(g),$$
one has to find a discrete version of all the involved terms, notably the  mean curvature.
We will follow \cite{Sullivan2008} to do this. 
Let $V,E,F$ denote the vertices, edges, and faces of the triangular mesh, and let 
$\on{star}(p)$ be the set of faces surrounding vertex $p$. 
Then the discrete mean curvature at vertex $p$ can be defined as
$$\Tr(L)(p) = \frac{\norm{\text{vector mean curvature}}}{\norm{\text{vector area}}} =
\frac{\norm{\nabla_p (\text{surface area})}}{\norm{\nabla_p (\text{enclosed volume})}}.$$
Here $\nabla_p$ stands for a discrete gradient, and
$$(\text{vector mean curvature})_p=\nabla_p (\text{surface area}) = 
\sum_{(p,p_i) \in E} (\cot \alpha_i + \cot \beta_i) (p-p_i)$$
is the vector mean curvature defined by the cotangent formula. In this formula,
$\alpha_i$ and $\beta_i$ are the angles opposite the edge $(p,p_i)$ in the two adjacent triangles. 
For the numerical simulation it is advantageous to express this formula in terms of scalar and cross
products instead of the cotangents.  
Furthermore, 
$$(\text{vector area})_p=\nabla_p (\text{enclosed volume}) = \sum_{f \in \on{Star}(p)} \nu(f).(\text{surface area of $f$})$$
is the vector area at vertex $p$. 

We discretize the time by  $$0=t_1<\ldots< t_{N+1}=1.$$ Then the $(N-1)(\#V)$ free variables representing the path of immersions
$f$ are
$$f(t_i,p) \qquad \text{with } 2\leq i \leq N,\ p\in V.$$  $f(0,p)$ and $f(1,p)$ are not free variables, since they define the fixed boundary shapes.
$f_t$ can be approximated by either forward increments $$f^{fw}_t(t_i,p)=\frac{f(t_{i+1},p)-f(t_i,p)}{t_{i+1}-t_i}$$ or backward increments
$$f^{bw}_t(t_i,p)=\frac{f(t_{i},p)-f(t_{i-1},p)}{t_{i}-t_{i-1}}.$$ We use a combination of both to make path energy symmetric. (Instead of this we could
 have  used the central difference quotient. However, minimizing an energy functional depending on central differences favors oscillations, 
since they are not felt by the central differences.) 
Using the discrete definitions of the normal vector and increments we can calculate $f_t^\bot$ at every vertex $p$
and are now able to write  the discrete horizontal  path energy:
\begin{equation*}\boxed{
\begin{aligned}
&G^{\Ph}_f(h^\bot,k^\bot)=\sum_{p\in V}\sum_{F\ni p} \Ph\Big(\Vol,\Tr(L)(p)\Big)
\\&\qquad\qquad\qquad
.\g\big( h(p),\nu(F)\big) .\g\big(k(p),\nu(F) \big)\frac{\on{area}(\on{star}^f(p))}{3}\\
&E^{\on{hor}}(f) =
 \sum_{i=1}^N\frac{t_{i+1}-t_i}{2}\\&\qquad\qquad\qquad.\Big(G^{\Ph}_{f(t_i)}\big(f^{fw}_t(t_i),f^{fw}_t(t_i)\big)
+G^{\Ph}_{f(t_{i+1})}\big(f^{bw}_t(t_{i+1}),f^{bw}_t(t_{i+1})\big)\Big) .
\end{aligned}}
\end{equation*}

This is not the only way to discretize the energy functional. There are several ways to distribute the discrete energy on faces, vertices, and edges. Depending on how this was done, the minimizer converged faster, slower, or even not at all. However, if the minimizer converged to a smooth solution, the results were qualitatively the same. This increased our belief in the discretization. 
However, we do not guarantee the accuracy of the simulations in this section. 

This energy functional does not depend on the parametrization of the surface at each instant of time. 
So we are free to choose a suitable 
parametrization. We do this by adding to the energy functional a term penalizing  
 irregular meshes.  So instead of minimizing horizontal path energy, we minimize
the sum of horizontal path energy and a penalty term. The penalty term  measures
the deviation of angles from the ``perfect angle'' $2\pi$ divided by the number of surrounding
triangles, i.e.,
$$\sum_{t=2}^N \ \sum_{p \in V} \sum_{(p,q,r) \in \Delta} 
\Big|\sphericalangle(pq,pr)-(\text{perfect angle}) \Big|^k, \quad k \in \mathbb N.$$

\subsection{Numerical implementation}\label{numerics:implementation}

\begin{figure}[ht]
\centering
\includegraphics[width=.29\textwidth]{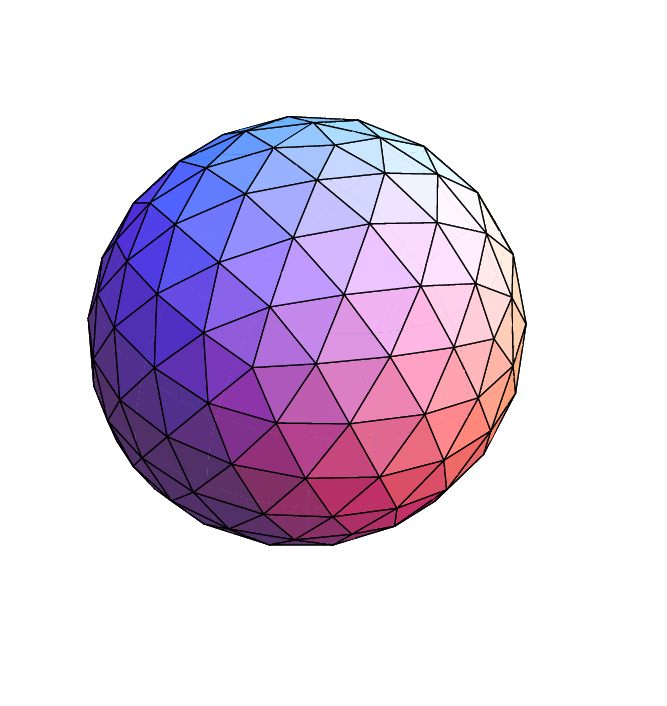}
\includegraphics[width=.29\textwidth]{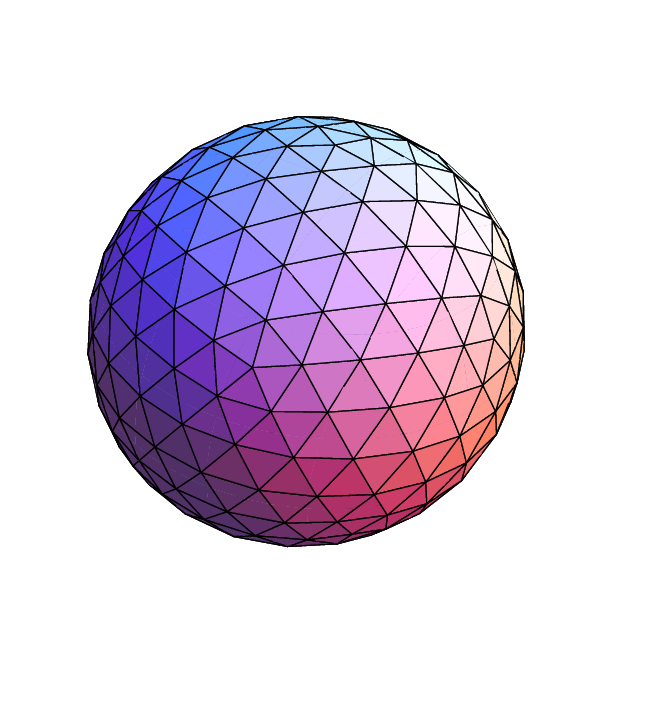}
\includegraphics[width=.29\textwidth]{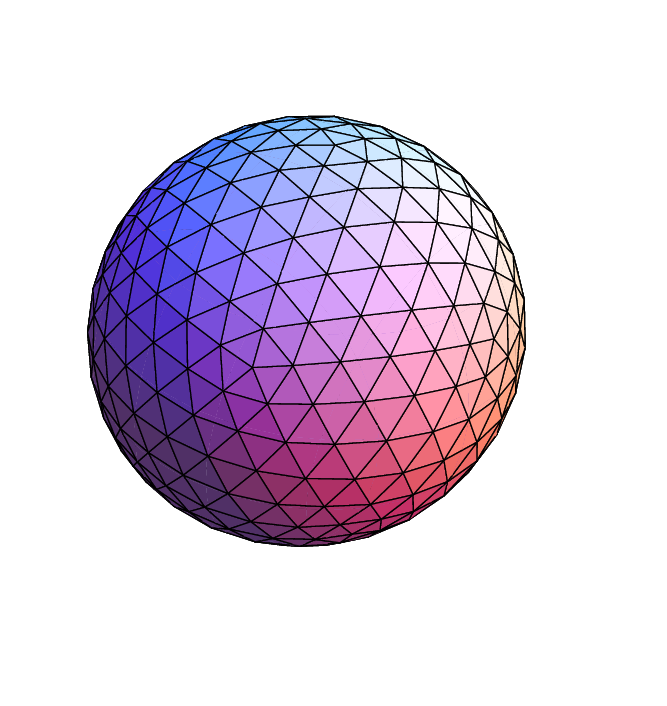}
\caption{Triangulations of a sphere with 320, 500 and 720 triangles, respectively.}
\label{fig:numerics:implementation:triangulation}
\end{figure}

Discrete path energy depends on a very high number of real variables, namely, three times the number of vertices times one 
less than the number of time steps. In the numerical experiments that we have done, there were between 5.000 and  50.000 variables. 
To solve this problem we used the nonlinear solver IPOPT (Interior Point OPTimizer \cite{Waechter2002}). 
IPOPT uses a filter based line search method  to compute the minimum. In this process it needs the gradient and the Hessian of the energy.  IPOPT was invoked by AMPL
(A Modeling Language for Mathematical Programming 
\cite{Fourer2002}). The advantage of using AMPL is that it is able to automatically and symbolically
calculate  the gradient and Hessian needed for the optimizer. 
All the user has to do is to write a model and data file for AMPL 
in a quite readable notation. The data file containing the definition of the combinatorics of the triangle mesh 
was automatically generated by the computer algebra system Mathematica. 
As an example, some discretizations of the sphere that we used can be seen in figure~\ref{fig:numerics:implementation:triangulation}.

\subsection{Scaling a sphere}\label{numerics:concentric_conformal}

In section \ref{con_spheres} we studied the set of concentric spheres in $n$ dimensions. In dimension three the 
geodesic equation for the radius simplifies to
\begin{equation*}\boxed{\begin{aligned}
r_{tt}&=-r_t^2 \frac{1}{\Ph} \big[ \frac{1}{r} \Ph 
+\p_\v \Ph 4r^2\pi   + \frac{1}{ r^2} (\p_\m\Ph)  +\frac{1}{r^2}(\p_3\Ph) \big]. 
\end{aligned}}\end{equation*}
{\it This equation is in accordance with the numerical results obtained by minimizing the discrete path energy defined in 
section~\ref{numerics:discrete}. As will be seen, the numerics show that the shortest path connecting two concentric spheres
in fact consists of spheres with the same center, and that the above differential equation is (at least qualitatively) satisfied. 
Furthermore, in our experiments the optimal paths obtained were independent of 
the initial path used as a starting value for the optimization. 
In all numerical experiments of this section we used 50 timesteps and a triangulation with 320 
triangles.}

\begin{figure}[htp]
\begin{center}
\begin{psfrags}
\def\PFGstripminus-#1{#1}%
\def\PFGshift(#1,#2)#3{\raisebox{#2}[\height][\depth]{\hbox{%
  \ifdim#1<0pt\kern#1 #3\kern\PFGstripminus#1\else\kern#1 #3\kern-#1\fi}}}%
\providecommand{\PFGstyle}{}%
%
\psfrag{Phi1Vol}[cl][cl]{\PFGstyle $\Ph=\frac{1}{\Vol}$}%
\psfrag{PhieVol}[cl][cl]{\PFGstyle $\Ph=e^{\Vol}$}%
\psfrag{PhiVol2}[cl][cl]{\PFGstyle $\Ph=\Vol^2$}%
\psfrag{PhiVol3}[cl][cl]{\PFGstyle $\Ph=\Vol^3$}%
\psfrag{PhiVol}[cl][cl]{\PFGstyle $\Ph=\Vol$}%
\psfrag{r}[bc][bc]{\PFGstyle $r$}%
\psfrag{t}[cl][cl]{\PFGstyle $t$}%
\psfrag{x0}[tc][tc]{\PFGstyle $0$}%
\psfrag{x11}[tc][tc]{\PFGstyle $1$}%
\psfrag{x1254}[tc][tc]{\PFGstyle $1250$}%
\psfrag{x14}[tc][tc]{\PFGstyle $1000$}%
\psfrag{x154}[tc][tc]{\PFGstyle $1500$}%
\psfrag{x1754}[tc][tc]{\PFGstyle $1750$}%
\psfrag{x21}[tc][tc]{\PFGstyle $2$}%
\psfrag{x253}[tc][tc]{\PFGstyle $250$}%
\psfrag{x2}[tc][tc]{\PFGstyle $0.2$}%
\psfrag{x41}[tc][tc]{\PFGstyle $4$}%
\psfrag{x4}[tc][tc]{\PFGstyle $0.4$}%
\psfrag{x53}[tc][tc]{\PFGstyle $500$}%
\psfrag{x61}[tc][tc]{\PFGstyle $6$}%
\psfrag{x6}[tc][tc]{\PFGstyle $0.6$}%
\psfrag{x753}[tc][tc]{\PFGstyle $750$}%
\psfrag{x81}[tc][tc]{\PFGstyle $8$}%
\psfrag{x8}[tc][tc]{\PFGstyle $0.8$}%
\psfrag{y0}[cr][cr]{\PFGstyle $0$}%
\psfrag{y11}[cr][cr]{\PFGstyle $1$}%
\psfrag{y151}[cr][cr]{\PFGstyle $1.5$}%
\psfrag{y1}[cr][cr]{\PFGstyle $0.1$}%
\psfrag{y21}[cr][cr]{\PFGstyle $2$}%
\psfrag{y2}[cr][cr]{\PFGstyle $0.2$}%
\psfrag{y3}[cr][cr]{\PFGstyle $0.3$}%
\psfrag{y4}[cr][cr]{\PFGstyle $0.4$}%
\psfrag{y5}[cr][cr]{\PFGstyle $0.5$}%
\psfrag{y6}[cr][cr]{\PFGstyle $0.6$}%
\psfrag{y7}[cr][cr]{\PFGstyle $0.7$}%
\psfrag{y8}[cr][cr]{\PFGstyle $0.8$}%
\psfrag{ym11}[cr][cr]{\PFGstyle $-1$}%
\psfrag{ym1253}[cr][cr]{\PFGstyle $-125$}%
\psfrag{ym13}[cr][cr]{\PFGstyle $-100$}%
\psfrag{ym151}[cr][cr]{\PFGstyle $-1.5$}%
\psfrag{ym153}[cr][cr]{\PFGstyle $-150$}%
\psfrag{ym21}[cr][cr]{\PFGstyle $-2$}%
\psfrag{ym252}[cr][cr]{\PFGstyle $-25$}%
\psfrag{ym52}[cr][cr]{\PFGstyle $-50$}%
\psfrag{ym5}[cr][cr]{\PFGstyle $-0.5$}%
\psfrag{ym752}[cr][cr]{\PFGstyle $-75$}%
\includegraphics[width=0.72\textwidth]{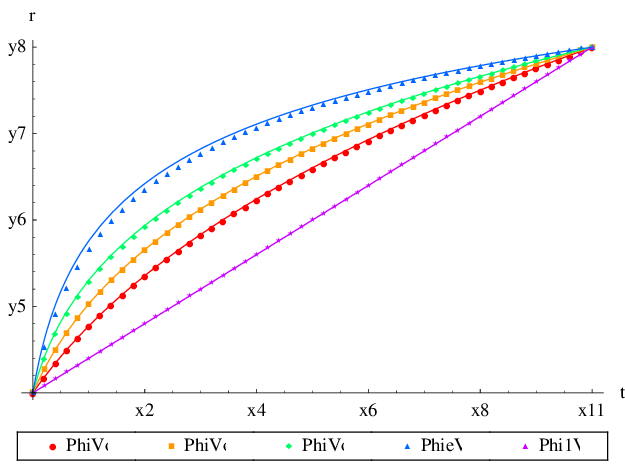}
\end{psfrags}
\caption{Geodesics between concentric spheres of radius $0.4$ to $0.8$ for several conformal metrics.
Solid lines are the exact solutions. }
\label{fig:numerics:concentric:bigradius}
\smallskip
\begin{psfrags}
\def\PFGstripminus-#1{#1}%
\def\PFGshift(#1,#2)#3{\raisebox{#2}[\height][\depth]{\hbox{%
  \ifdim#1<0pt\kern#1 #3\kern\PFGstripminus#1\else\kern#1 #3\kern-#1\fi}}}%
\providecommand{\PFGstyle}{}%
%
\psfrag{Phi1Vol}[cl][cl]{\PFGstyle $\Ph=\frac{1}{\Vol}$}%
\psfrag{PhieVol}[cl][cl]{\PFGstyle $\Ph=e^{\Vol}$}%
\psfrag{PhiVol2}[cl][cl]{\PFGstyle $\Ph=\Vol^2$}%
\psfrag{PhiVol3}[cl][cl]{\PFGstyle $\Ph=\Vol^3$}%
\psfrag{PhiVol}[cl][cl]{\PFGstyle $\Ph=\Vol$}%
\psfrag{r}[bc][bc]{\PFGstyle $r$}%
\psfrag{t}[cl][cl]{\PFGstyle $t$}%
\psfrag{x0}[tc][tc]{\PFGstyle $0$}%
\psfrag{x11}[tc][tc]{\PFGstyle $1$}%
\psfrag{x1254}[tc][tc]{\PFGstyle $1250$}%
\psfrag{x14}[tc][tc]{\PFGstyle $1000$}%
\psfrag{x154}[tc][tc]{\PFGstyle $1500$}%
\psfrag{x1754}[tc][tc]{\PFGstyle $1750$}%
\psfrag{x21}[tc][tc]{\PFGstyle $2$}%
\psfrag{x253}[tc][tc]{\PFGstyle $250$}%
\psfrag{x2}[tc][tc]{\PFGstyle $0.2$}%
\psfrag{x41}[tc][tc]{\PFGstyle $4$}%
\psfrag{x4}[tc][tc]{\PFGstyle $0.4$}%
\psfrag{x53}[tc][tc]{\PFGstyle $500$}%
\psfrag{x61}[tc][tc]{\PFGstyle $6$}%
\psfrag{x6}[tc][tc]{\PFGstyle $0.6$}%
\psfrag{x753}[tc][tc]{\PFGstyle $750$}%
\psfrag{x81}[tc][tc]{\PFGstyle $8$}%
\psfrag{x8}[tc][tc]{\PFGstyle $0.8$}%
\psfrag{y0}[cr][cr]{\PFGstyle $0$}%
\psfrag{y11}[cr][cr]{\PFGstyle $1$}%
\psfrag{y12}[cr][cr]{\PFGstyle $0.12$}%
\psfrag{y14}[cr][cr]{\PFGstyle $0.14$}%
\psfrag{y151}[cr][cr]{\PFGstyle $1.5$}%
\psfrag{y16}[cr][cr]{\PFGstyle $0.16$}%
\psfrag{y18}[cr][cr]{\PFGstyle $0.18$}%
\psfrag{y1}[cr][cr]{\PFGstyle $0.1$}%
\psfrag{y21}[cr][cr]{\PFGstyle $2$}%
\psfrag{y2}[cr][cr]{\PFGstyle $0.2$}%
\psfrag{y3}[cr][cr]{\PFGstyle $0.3$}%
\psfrag{y4}[cr][cr]{\PFGstyle $0.4$}%
\psfrag{y5}[cr][cr]{\PFGstyle $0.5$}%
\psfrag{y6}[cr][cr]{\PFGstyle $0.6$}%
\psfrag{ym11}[cr][cr]{\PFGstyle $-1$}%
\psfrag{ym1253}[cr][cr]{\PFGstyle $-125$}%
\psfrag{ym13}[cr][cr]{\PFGstyle $-100$}%
\psfrag{ym151}[cr][cr]{\PFGstyle $-1.5$}%
\psfrag{ym153}[cr][cr]{\PFGstyle $-150$}%
\psfrag{ym21}[cr][cr]{\PFGstyle $-2$}%
\psfrag{ym252}[cr][cr]{\PFGstyle $-25$}%
\psfrag{ym52}[cr][cr]{\PFGstyle $-50$}%
\psfrag{ym5}[cr][cr]{\PFGstyle $-0.5$}%
\psfrag{ym752}[cr][cr]{\PFGstyle $-75$}%
\includegraphics[width=0.72\textwidth]{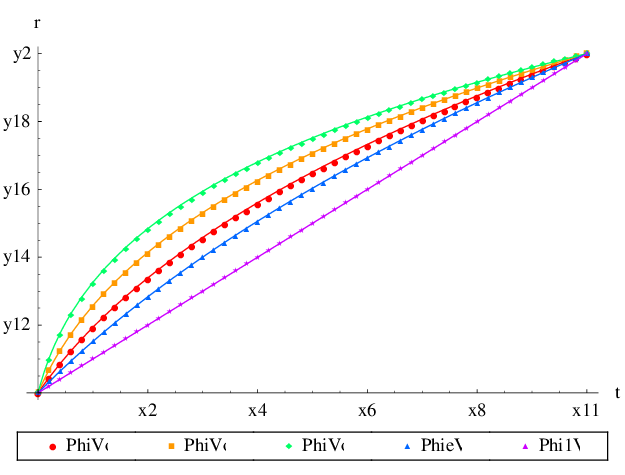}
\end{psfrags}
\caption{Geodesics between concentric spheres of radius $0.1$ to $0.2$ for several conformal metrics.
Solid lines are the exact solutions. }
\label{fig:numerics:concentric:smallradius}
\end{center}
\end{figure}

For conformal metrics of type $\Ph=\Vol^k$ and $\Ph=e^{\Vol}$, the ODE for the radius is
\begin{align*}
& \Ph=\Vol^k: && r_{tt}=-r_t^2  \frac{k+1}{r}, \\
& \Ph=e^{\Vol}: && r_{tt}=-r_t^2  \Big(\frac1r + 4 r \pi\Big).
\end{align*}
Note that the equation for $\Ph=\Vol^{-1}$ is $r_{tt}=0$. 
These equations have explicit solutions: 
\begin{align*}
& \Ph=\Vol^k: && r= C_1 \big((k+2)t-C_2 \big)^{\frac{1}{k+2}}, \\
& \Ph=e^{\Vol}: && r= \frac{1}{2\pi} \sqrt{\log\big(C_1 t + C_2\big)}.
\end{align*}
A comparison of the numerical results with the exact analytic solutions can be seen in 
Figures~\ref{fig:numerics:concentric:bigradius} and \ref{fig:numerics:concentric:smallradius}. 
The solid lines are the exact solutions. 
Note that for big radii as in Figure~\ref{fig:numerics:concentric:bigradius}, 
the solution for $\Ph=e^{\Vol}$ has a very steep ascent, is more curved, and lies above the solutions for $\Ph=\Vol, \Vol^2, \Vol^3$.
For small radii, it lies below these solutions, as can be seen in figure~\ref{fig:numerics:concentric:smallradius}.
Note also that when the ascent gets too steep, the discrete solution is somewhat inexact as in 
Figure~\ref{fig:numerics:concentric:bigradius}.

\begin{figure}[htp]
\begin{center}
\begin{psfrags}
\def\PFGstripminus-#1{#1}%
\def\PFGshift(#1,#2)#3{\raisebox{#2}[\height][\depth]{\hbox{%
  \ifdim#1<0pt\kern#1 #3\kern\PFGstripminus#1\else\kern#1 #3\kern-#1\fi}}}%
\providecommand{\PFGstyle}{}%
%
\psfrag{k2}[cl][cl]{\PFGstyle $k=2$}%
\psfrag{k4}[cl][cl]{\PFGstyle $k=4$}%
\psfrag{k6}[cl][cl]{\PFGstyle $k=6$}%
\psfrag{k8}[cl][cl]{\PFGstyle $k=8$}%
\psfrag{r}[bc][bc]{\PFGstyle $r$}%
\psfrag{t}[cl][cl]{\PFGstyle $t$}%
\psfrag{x0}[tc][tc]{\PFGstyle $0$}%
\psfrag{x11}[tc][tc]{\PFGstyle $1$}%
\psfrag{x124}[tc][tc]{\PFGstyle $1200$}%
\psfrag{x144}[tc][tc]{\PFGstyle $1400$}%
\psfrag{x14}[tc][tc]{\PFGstyle $1000$}%
\psfrag{x21}[tc][tc]{\PFGstyle $2$}%
\psfrag{x23}[tc][tc]{\PFGstyle $200$}%
\psfrag{x2}[tc][tc]{\PFGstyle $0.2$}%
\psfrag{x41}[tc][tc]{\PFGstyle $4$}%
\psfrag{x43}[tc][tc]{\PFGstyle $400$}%
\psfrag{x4}[tc][tc]{\PFGstyle $0.4$}%
\psfrag{x61}[tc][tc]{\PFGstyle $6$}%
\psfrag{x63}[tc][tc]{\PFGstyle $600$}%
\psfrag{x6}[tc][tc]{\PFGstyle $0.6$}%
\psfrag{x81}[tc][tc]{\PFGstyle $8$}%
\psfrag{x83}[tc][tc]{\PFGstyle $800$}%
\psfrag{x8}[tc][tc]{\PFGstyle $0.8$}%
\psfrag{y0}[cr][cr]{\PFGstyle $0$}%
\psfrag{y11}[cr][cr]{\PFGstyle $1$}%
\psfrag{y121}[cr][cr]{\PFGstyle $1.2$}%
\psfrag{y141}[cr][cr]{\PFGstyle $1.4$}%
\psfrag{y151}[cr][cr]{\PFGstyle $1.5$}%
\psfrag{y161}[cr][cr]{\PFGstyle $1.6$}%
\psfrag{y181}[cr][cr]{\PFGstyle $1.8$}%
\psfrag{y1}[cr][cr]{\PFGstyle $0.1$}%
\psfrag{y21}[cr][cr]{\PFGstyle $2$}%
\psfrag{y2}[cr][cr]{\PFGstyle $0.2$}%
\psfrag{y3}[cr][cr]{\PFGstyle $0.3$}%
\psfrag{y4}[cr][cr]{\PFGstyle $0.4$}%
\psfrag{y5}[cr][cr]{\PFGstyle $0.5$}%
\psfrag{y6}[cr][cr]{\PFGstyle $0.6$}%
\psfrag{ym11}[cr][cr]{\PFGstyle $-1$}%
\psfrag{ym1253}[cr][cr]{\PFGstyle $-125$}%
\psfrag{ym13}[cr][cr]{\PFGstyle $-100$}%
\psfrag{ym151}[cr][cr]{\PFGstyle $-1.5$}%
\psfrag{ym153}[cr][cr]{\PFGstyle $-150$}%
\psfrag{ym21}[cr][cr]{\PFGstyle $-2$}%
\psfrag{ym252}[cr][cr]{\PFGstyle $-25$}%
\psfrag{ym52}[cr][cr]{\PFGstyle $-50$}%
\psfrag{ym5}[cr][cr]{\PFGstyle $-0.5$}%
\psfrag{ym752}[cr][cr]{\PFGstyle $-75$}%
\includegraphics[width=0.72\textwidth]{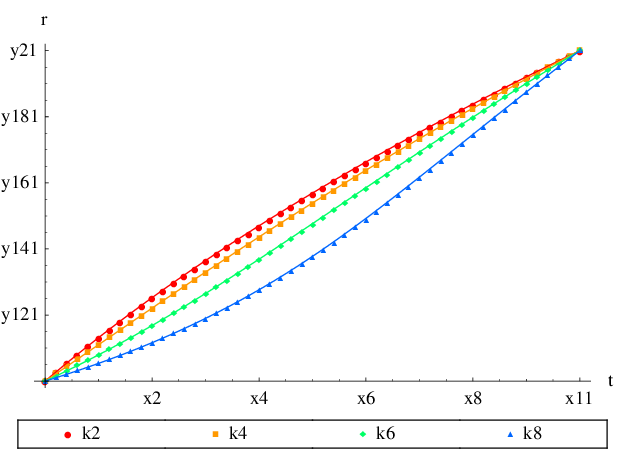}
\end{psfrags}
\caption{Geodesics between concentric spheres for $\Ph=1+0.1 \on{Tr}(L)^k$  and varying $k$.
Solid lines are the exact solutions. }
\label{fig:numerics:concentric:ga_varyingk}
\smallskip
\begin{psfrags}
\def\PFGstripminus-#1{#1}%
\def\PFGshift(#1,#2)#3{\raisebox{#2}[\height][\depth]{\hbox{%
  \ifdim#1<0pt\kern#1 #3\kern\PFGstripminus#1\else\kern#1 #3\kern-#1\fi}}}%
\providecommand{\PFGstyle}{}%
%
\psfrag{A001}[cl][cl]{\PFGstyle $A=0.01$}%
\psfrag{A01}[cl][cl]{\PFGstyle $A=0.1$}%
\psfrag{A0}[cl][cl]{\PFGstyle $A=0$}%
\psfrag{A1}[cl][cl]{\PFGstyle $A=1$}%
\psfrag{r}[bc][bc]{\PFGstyle $r$}%
\psfrag{t}[cl][cl]{\PFGstyle $t$}%
\psfrag{x0}[tc][tc]{\PFGstyle $0$}%
\psfrag{x11}[tc][tc]{\PFGstyle $1$}%
\psfrag{x124}[tc][tc]{\PFGstyle $1200$}%
\psfrag{x144}[tc][tc]{\PFGstyle $1400$}%
\psfrag{x14}[tc][tc]{\PFGstyle $1000$}%
\psfrag{x21}[tc][tc]{\PFGstyle $2$}%
\psfrag{x23}[tc][tc]{\PFGstyle $200$}%
\psfrag{x2}[tc][tc]{\PFGstyle $0.2$}%
\psfrag{x41}[tc][tc]{\PFGstyle $4$}%
\psfrag{x43}[tc][tc]{\PFGstyle $400$}%
\psfrag{x4}[tc][tc]{\PFGstyle $0.4$}%
\psfrag{x61}[tc][tc]{\PFGstyle $6$}%
\psfrag{x63}[tc][tc]{\PFGstyle $600$}%
\psfrag{x6}[tc][tc]{\PFGstyle $0.6$}%
\psfrag{x81}[tc][tc]{\PFGstyle $8$}%
\psfrag{x83}[tc][tc]{\PFGstyle $800$}%
\psfrag{x8}[tc][tc]{\PFGstyle $0.8$}%
\psfrag{y0}[cr][cr]{\PFGstyle $0$}%
\psfrag{y11}[cr][cr]{\PFGstyle $1$}%
\psfrag{y151}[cr][cr]{\PFGstyle $1.5$}%
\psfrag{y1}[cr][cr]{\PFGstyle $0.1$}%
\psfrag{y21}[cr][cr]{\PFGstyle $2$}%
\psfrag{y2}[cr][cr]{\PFGstyle $0.2$}%
\psfrag{y3}[cr][cr]{\PFGstyle $0.3$}%
\psfrag{y4}[cr][cr]{\PFGstyle $0.4$}%
\psfrag{y5}[cr][cr]{\PFGstyle $0.5$}%
\psfrag{y6}[cr][cr]{\PFGstyle $0.6$}%
\psfrag{y8}[cr][cr]{\PFGstyle $0.8$}%
\psfrag{ym11}[cr][cr]{\PFGstyle $-1$}%
\psfrag{ym1253}[cr][cr]{\PFGstyle $-125$}%
\psfrag{ym13}[cr][cr]{\PFGstyle $-100$}%
\psfrag{ym151}[cr][cr]{\PFGstyle $-1.5$}%
\psfrag{ym153}[cr][cr]{\PFGstyle $-150$}%
\psfrag{ym21}[cr][cr]{\PFGstyle $-2$}%
\psfrag{ym252}[cr][cr]{\PFGstyle $-25$}%
\psfrag{ym52}[cr][cr]{\PFGstyle $-50$}%
\psfrag{ym5}[cr][cr]{\PFGstyle $-0.5$}%
\psfrag{ym752}[cr][cr]{\PFGstyle $-75$}%
\includegraphics[width=0.72\textwidth]{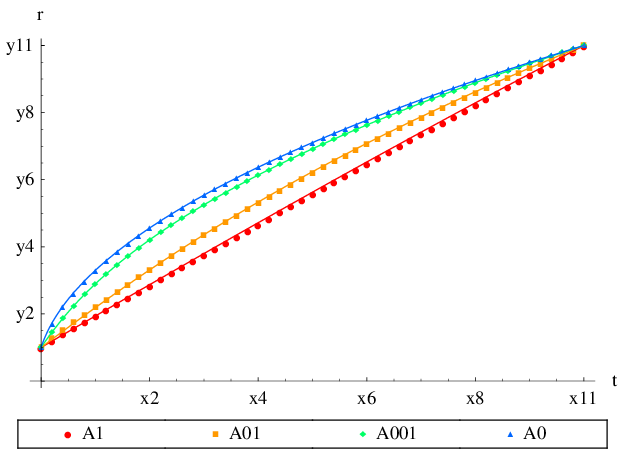}
\end{psfrags}
\caption{Geodesics between concentric spheres for $\Ph=1+A \on{Tr}(L)^2$ and varying $A$.
Solid lines are the exact solutions.}
\label{fig:numerics:concentric:ga_varyingA}
\end{center}
\end{figure}

\begin{figure}[htp]
\begin{center}
\begin{psfrags}
\def\PFGstripminus-#1{#1}%
\def\PFGshift(#1,#2)#3{\raisebox{#2}[\height][\depth]{\hbox{%
  \ifdim#1<0pt\kern#1 #3\kern\PFGstripminus#1\else\kern#1 #3\kern-#1\fi}}}%
\providecommand{\PFGstyle}{}%
%
\psfrag{Phisi}[cl][cl]{\PFGstyle $\Ph={\Vol^{-2}}+A\Tr(L)^2\Vol^{-1}$}%
\psfrag{r}[bc][bc]{\PFGstyle $r$}%
\psfrag{t}[cl][cl]{\PFGstyle $t$}%
\psfrag{x0}[tc][tc]{\PFGstyle $0$}%
\psfrag{x11}[tc][tc]{\PFGstyle $1$}%
\psfrag{x12}[tc][tc]{\PFGstyle $10$}%
\psfrag{x13}[tc][tc]{\PFGstyle $100$}%
\psfrag{x152}[tc][tc]{\PFGstyle $15$}%
\psfrag{x153}[tc][tc]{\PFGstyle $150$}%
\psfrag{x22}[tc][tc]{\PFGstyle $20$}%
\psfrag{x23}[tc][tc]{\PFGstyle $200$}%
\psfrag{x252}[tc][tc]{\PFGstyle $25$}%
\psfrag{x253}[tc][tc]{\PFGstyle $250$}%
\psfrag{x2}[tc][tc]{\PFGstyle $0.2$}%
\psfrag{x33}[tc][tc]{\PFGstyle $300$}%
\psfrag{x353}[tc][tc]{\PFGstyle $350$}%
\psfrag{x4}[tc][tc]{\PFGstyle $0.4$}%
\psfrag{x51}[tc][tc]{\PFGstyle $5$}%
\psfrag{x52}[tc][tc]{\PFGstyle $50$}%
\psfrag{x6}[tc][tc]{\PFGstyle $0.6$}%
\psfrag{x8}[tc][tc]{\PFGstyle $0.8$}%
\psfrag{y0}[cr][cr]{\PFGstyle $0$}%
\psfrag{y11}[cr][cr]{\PFGstyle $1$}%
\psfrag{y12}[cr][cr]{\PFGstyle $0.12$}%
\psfrag{y14}[cr][cr]{\PFGstyle $0.14$}%
\psfrag{y151}[cr][cr]{\PFGstyle $1.5$}%
\psfrag{y16}[cr][cr]{\PFGstyle $0.16$}%
\psfrag{y18}[cr][cr]{\PFGstyle $0.18$}%
\psfrag{y1}[cr][cr]{\PFGstyle $0.1$}%
\psfrag{y21}[cr][cr]{\PFGstyle $2$}%
\psfrag{y2}[cr][cr]{\PFGstyle $0.2$}%
\psfrag{y3}[cr][cr]{\PFGstyle $0.3$}%
\psfrag{y4}[cr][cr]{\PFGstyle $0.4$}%
\psfrag{y5}[cr][cr]{\PFGstyle $0.5$}%
\psfrag{y6}[cr][cr]{\PFGstyle $0.6$}%
\psfrag{ym11}[cr][cr]{\PFGstyle $-1$}%
\psfrag{ym12}[cr][cr]{\PFGstyle $-10$}%
\psfrag{ym151}[cr][cr]{\PFGstyle $-1.5$}%
\psfrag{ym21}[cr][cr]{\PFGstyle $-2$}%
\psfrag{ym22}[cr][cr]{\PFGstyle $-20$}%
\psfrag{ym32}[cr][cr]{\PFGstyle $-30$}%
\psfrag{ym42}[cr][cr]{\PFGstyle $-40$}%
\psfrag{ym52}[cr][cr]{\PFGstyle $-50$}%
\psfrag{ym5}[cr][cr]{\PFGstyle $-0.5$}%
\includegraphics[width=0.72\textwidth]{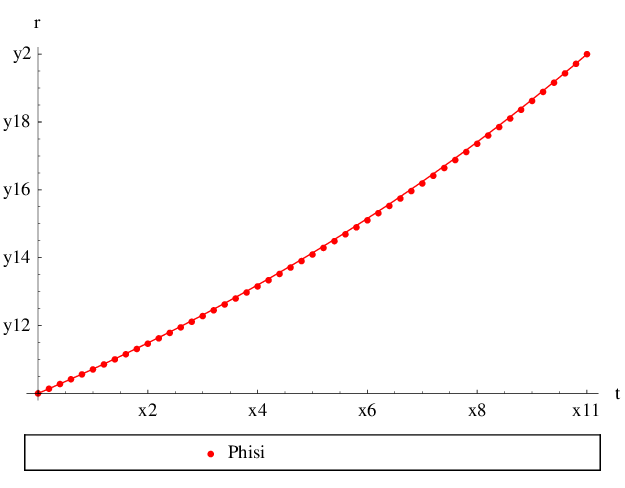}
\end{psfrags}
\caption{Geodesics between concentric spheres for the scale-invariant metric.}
\label{fig:numerics:concentric:si}
\end{center}
\end{figure}

For mean curvature weighted metrics, the differential equation for the radius is
\begin{align*}
&\Ph=1+A \on{Tr}(L)^{2k}: && r_{tt}=-r_t^2 \Big(\frac1r - \frac{2k A 2^{2k-1}}{r^{2k+1}+A 2^{2k} r}\Big).
\end{align*}

The numerics for these metrics are shown in figure~\ref{fig:numerics:concentric:ga_varyingk} and figure~\ref{fig:numerics:concentric:ga_varyingA}.
Note that we got convergence to a path consisting of concentric spheres even for the $G^0$-metric ($A=0$), even though
we know from the theory that this is not the shortest path. In fact, there are no shortest paths for the $G^0$--metric since 
it has vanishing geodesic distance \cite{Michor102}.

For the scale-invariant metric, the differential equation is given by
\begin{align*}
&\Ph=\Vol^{-2}+A \frac{\Tr(L)^2}{\Vol}: && r_{tt}= \frac{r_t^2}{r}.
\intertext{This equation has an explicit analytical solution}
&\Ph=\Vol^{-2}+A \frac{\Tr(L)^2}{\Vol}: && r=C_1 e^{C_2 t}.
\end{align*}
Note that this equation, and therefore its solution, is independent of $A$. 
Again, this is confirmed by the numerics: see Figure \ref{fig:numerics:concentric:si}. 

\subsection{Translation of a sphere}
In this section we will study geodesics between a sphere and a translated sphere for various almost 
local metrics 
of the type $\Ph=\Vol^k$, $\Ph=e^{\Vol}$, and $\Ph=1+A\Tr(L)^{2k}$.

Depending on the distance (relative to the radius) of the two translated spheres,  different behaviors can be observed. 

\noindent
{\bf High distance:} 
\begin{itemize}
\item {\it Shrink and grow:} For some metrics it is possible to shrink a sphere in finite time to zero. For these metrics long translation goes via a shrinking part and growing part. Metrics with this behavior are $\Ph=\Vol^k$, $\Ph=e^{\Vol}$ and $\Ph=1+A\Tr(L)^{2}$.  This phenomenon is studied in more detail in section~\ref{shrinkandblow}; see also Figure~\ref{fig:numerics:move:longdistance}.
\item {\it Moving an optimal middle shape:} For some of the metrics translation of a sphere with a certain optimal radius is a geodesic. For these metrics geodesics for long translations scale the sphere to the optimal radius and translate the sphere with the optimal radius.  Metrics with this behavior are $\Ph=1+A\Tr(L)^{2k}$ for $k>1$. This behavior is studied in section~\ref{optimalshape}.
\end{itemize}

\begin{figure}[ht]
\centering
\includegraphics[width=\textwidth]{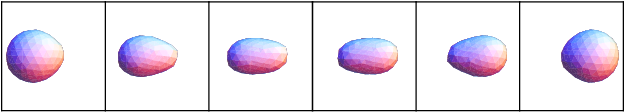}
\caption{Geodesic between two unit spheres translated by distance $1.5$ for $\Ph=\Vol$. 20 timesteps and a triangulation with 500 triangles were used. Time progresses from left to right. 
Boundary shapes $t=0$ and $t=1$ are not included. }
\label{fig:numerics:move:middledistance}
\end{figure}

\noindent
{\bf Low distance:} 
\begin{itemize}
\item Geodesics of pure translation ($\Ph=1+A\Tr(L)^{2k}$ for $k>1$; c.f. Figure~\ref{fig:numerics:move:trace4andtrace6}).
\item Geodesics that pass through an ellipsoid, where the longer principal axis is in the direction of the translation (conformal metrics, c.f. Figure~\ref{fig:numerics:move:middledistance}).
\item Geodesics that pass through an  ellipsoid, where the principal axis in the direction of the translation is shorter ($\Ph=1+A\Tr(L)^{2k}$ for $k>1$, c.f. Figure~\ref{fig:numerics:move:trace4andtrace6}).
\item Geodesics that pass through a cigar--shaped figure ($\Ph=1+A\Tr(L)^{2},$ c.f. Figure~\ref{fig:numerics:move:ga_cigar}).
\end{itemize}

\subsection{Shrink and grow}\label{shrinkandblow}

\begin{figure}[ht]
\centering
\includegraphics[width=\textwidth]{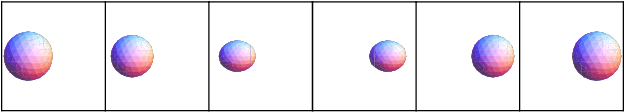}
\caption{Geodesic between two unit spheres translated by distance $2$ for $\Ph=e^{\Vol}$. 20 timesteps and a triangulation with 500 triangles were used. Time progresses from left to right. 
Boundary shapes $t=0$ and $t=1$ are not included. }
\label{fig:numerics:move:longdistance}
\end{figure}

In section~\ref{con_spheres} we showed that it is possible to shrink a sphere to zero in finite time for some of the metrics, namely, conformal metrics with $\Ph=\Vol^k$ or $\Ph=e^{\Vol}$ and for the $G^A$--metric. For these metrics geodesics of long translation will go via a shrinking part and growing part, and almost all of the translation will be done with the shrunken version of the shape. An example of such a geodesic can be seen in figure~\ref{fig:numerics:move:longdistance}.  

We could not determine numerically whether a collapse of the sphere to a point occurs or not. But the more
time steps were used, the smaller the ellipsoid in the middle turned out. Also, the energy of the geodesic path comes very close to the 
energy needed to shrink the sphere to a point and blow it up again. It is remarkable that almost all of the translation is concentrated at 
a single time step, independently of the number of timesteps that were used. The reason for this behavior is that high volumes are penalized so much: 
In the case of figure~\ref{fig:numerics:move:longdistance}, 
$e^{\Vol}$ is more than 1000 times  smaller in the middle than at the boundary shapes.

We now want find out under what conditions on the distance and radius of the boundary spheres of the geodesic this behavior can occur.
To do this, we compare the energy needed for a pure translation with the energy needed to first shrink the sphere to almost zero, 
then move it, and then blow it up again. 

The energy needed for a pure translation of a sphere with radius $r$ by distance $\ell$ in the direction of a unit vector $e_1$ is given by
\begin{align*}
E &=  \int_0^1\int_{S^2} \Ph(\Vol,\Tr(L)) \g( \ell.e_1, \nu )^2 \vol(g) dt \\&= 
\Ph(4r^2\pi,-\frac{2}{r}) \int_0^\pi \int_0^{2\pi} \g\left( \ell.e_1, 
\left(\begin{array}{c}\cos \varphi \sin \theta \\ \sin \varphi \sin \theta \\ \cos \theta \end{array}\right) 
\right)^2 r^2 \sin \theta\ d\varphi d\theta \\&= 
\Ph(4r^2\pi,-\frac{2}{r}) \int_0^\pi \int_0^{2\pi} \ell^2. (\cos \varphi \sin \theta)^2  r^2 \sin \theta d\varphi d\theta =\Ph(4r^2\pi,-\frac{2}{r}) . \frac{4\pi}{3} \ell^2.  r^2 . 
\end{align*}
Any other unit vector can be chosen instead of $e_1$, yielding the same result.

\begin{figure}[ht]
\centering
\begin{psfrags}
\def\PFGstripminus-#1{#1}%
\def\PFGshift(#1,#2)#3{\raisebox{#2}[\height][\depth]{\hbox{%
  \ifdim#1<0pt\kern#1 #3\kern\PFGstripminus#1\else\kern#1 #3\kern-#1\fi}}}%
\providecommand{\PFGstyle}{}%
%
\psfrag{l}[bc][bc]{\PFGstyle $l$}%
\psfrag{Phi1TrL2}[cl][cl]{\PFGstyle $G^1$}%
\psfrag{PhieVol}[cl][cl]{\PFGstyle $\Ph=e^{\Vol}$}%
\psfrag{PhiVol2}[cl][cl]{\PFGstyle $\Ph=\Vol^2$}%
\psfrag{PhiVol3}[cl][cl]{\PFGstyle $\Ph=\Vol^3$}%
\psfrag{PhiVol}[cl][cl]{\PFGstyle $\Ph=\Vol$}%
\psfrag{rA}[cl][cl]{\PFGstyle $r$}%
\psfrag{r}[bc][bc]{\PFGstyle $r$}%
\psfrag{t}[cl][cl]{\PFGstyle $t$}%
\psfrag{x00}[tc][tc]{\PFGstyle $0.0$}%
\psfrag{x05}[tc][tc]{\PFGstyle $0.5$}%
\psfrag{x0}[tc][tc]{\PFGstyle $0$}%
\psfrag{x10}[tc][tc]{\PFGstyle $1.0$}%
\psfrag{x11}[tc][tc]{\PFGstyle $1$}%
\psfrag{x1254}[tc][tc]{\PFGstyle $1250$}%
\psfrag{x12}[tc][tc]{\PFGstyle $10$}%
\psfrag{x14}[tc][tc]{\PFGstyle $1000$}%
\psfrag{x154}[tc][tc]{\PFGstyle $1500$}%
\psfrag{x15}[tc][tc]{\PFGstyle $1.5$}%
\psfrag{x1754}[tc][tc]{\PFGstyle $1750$}%
\psfrag{x21}[tc][tc]{\PFGstyle $2$}%
\psfrag{x253}[tc][tc]{\PFGstyle $250$}%
\psfrag{x2}[tc][tc]{\PFGstyle $0.2$}%
\psfrag{x41}[tc][tc]{\PFGstyle $4$}%
\psfrag{x4}[tc][tc]{\PFGstyle $0.4$}%
\psfrag{x53}[tc][tc]{\PFGstyle $500$}%
\psfrag{x61}[tc][tc]{\PFGstyle $6$}%
\psfrag{x6}[tc][tc]{\PFGstyle $0.6$}%
\psfrag{x753}[tc][tc]{\PFGstyle $750$}%
\psfrag{x81}[tc][tc]{\PFGstyle $8$}%
\psfrag{x8}[tc][tc]{\PFGstyle $0.8$}%
\psfrag{y00}[cr][cr]{\PFGstyle $0.0$}%
\psfrag{y05}[cr][cr]{\PFGstyle $0.5$}%
\psfrag{y0}[cr][cr]{\PFGstyle $0$}%
\psfrag{y10}[cr][cr]{\PFGstyle $1.0$}%
\psfrag{y11}[cr][cr]{\PFGstyle $1$}%
\psfrag{y151}[cr][cr]{\PFGstyle $1.5$}%
\psfrag{y1}[cr][cr]{\PFGstyle $0.1$}%
\psfrag{y21}[cr][cr]{\PFGstyle $2$}%
\psfrag{y2}[cr][cr]{\PFGstyle $0.2$}%
\psfrag{y3}[cr][cr]{\PFGstyle $0.3$}%
\psfrag{y4}[cr][cr]{\PFGstyle $0.4$}%
\psfrag{y5}[cr][cr]{\PFGstyle $0.5$}%
\psfrag{y6}[cr][cr]{\PFGstyle $0.6$}%
\psfrag{y8}[cr][cr]{\PFGstyle $0.8$}%
\psfrag{ym11}[cr][cr]{\PFGstyle $-1$}%
\psfrag{ym123}[cr][cr]{\PFGstyle $-120$}%
\psfrag{ym13}[cr][cr]{\PFGstyle $-100$}%
\psfrag{ym151}[cr][cr]{\PFGstyle $-1.5$}%
\psfrag{ym21}[cr][cr]{\PFGstyle $-2$}%
\psfrag{ym22}[cr][cr]{\PFGstyle $-20$}%
\psfrag{ym42}[cr][cr]{\PFGstyle $-40$}%
\psfrag{ym5}[cr][cr]{\PFGstyle $-0.5$}%
\psfrag{ym62}[cr][cr]{\PFGstyle $-60$}%
\psfrag{ym82}[cr][cr]{\PFGstyle $-80$}%
\includegraphics[width=\textwidth]{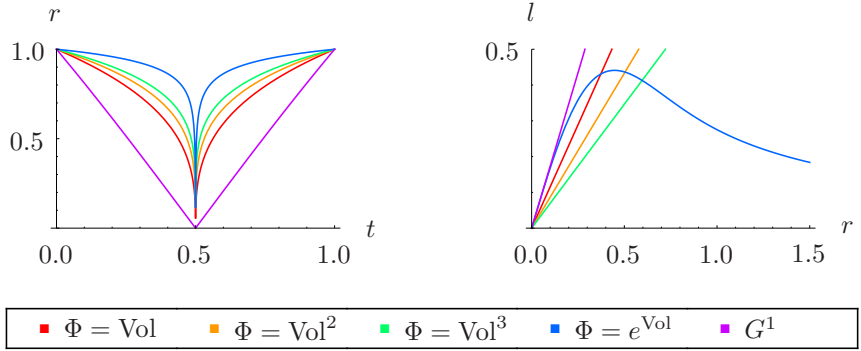}
\end{psfrags}
\caption{Left: Shrinking a sphere to zero along a geodesic path and blowing it up again. 
Right: Pairs of $\ell$ and $r$ such that translating a sphere of radius $r$ by distance $\ell$ needs as much energy as shrinking it to zero and 
blowing it up again. $G^1$ stands for the $G^A$ metric with $A=1$.}
\label{fig:numerics:move:shrinknblow}
\end{figure}

We will now calculate the energy needed for shrinking the sphere, moving it, and blowing it up again. 
The energy needed for translating a sphere of radius almost zero can be neglected. 
Shrinking and blowing up are done using the solutions to the geodesic equation for the radius from the last section, where 
one has to adapt the constants to the boundary conditions. For the shrinking part we have $r(0)=r$ and $r(\tfrac12)=0$, and for the growing part
we have $r(\tfrac12)=0,r(1)=r$; see Figure~\ref{fig:numerics:move:shrinknblow}~(left).

The energy of the path is
\begin{align*}
& \Ph=\Vol^k: && E= \int_0^1 \Vol^k \int_{S^2} r_t^2 \vol(g) dt = \frac{4^{k+2} \pi^{k+1}}{(k+2)^2} r^{2k+4}, \\
& \Ph=e^{\Vol}: && E=\int_0^1 e^{\Vol} \int_{S^2} r_t^2 \vol(g) dt = \frac{1}{\pi}(e^{2\pi r^2}-1)^2.
\end{align*}
The energy of the two different paths are the same when 
\begin{align*}
& \Ph=\Vol^k: &&   \ell   = \frac{2 \sqrt{3}r }{k+2}, \\
& \Ph=e^{\Vol}: &&   \ell   = \frac{\sqrt{3}(1-e^{-2\pi r^2})}{2 r\pi } .
\end{align*}
These curves are shown in Figure~~\ref{fig:numerics:move:shrinknblow}~(right). We did not derive an analytic solution for the $G^A$--metric, but for $A=1$  one can see the solution curves in figure~\ref{fig:numerics:move:shrinknblow}.

\subsection{Moving an optimal shape}\label{optimalshape}

In the following we want to determine whether pure translation of a sphere is a geodesic.
Therefore, let $f_t=f_0+b(t) \cdot e_1$, where $f_0$ is a sphere of radius $r$ and where $b(t)$ is constant on $M$. Plugging this into the geodesic equation from 
section~\ref{geodesic_equation_imm:geodesic_equation} yields an ODE for $b(t)$ and a part which has to vanish identically. The latter is given by
\begin{equation}\label{radius}
(\p_\v \Ph) \frac{2}{r} 4r^2\pi +
(\p_\m\Ph) \frac{2}{r^2} 
+\Ph  \frac{2}{r}=0 
\end{equation}
For conformal metrics this equation is  satisfied only if $\Ph=\Vol\i$. Since this metric induces vanishing geodesic distance (see section~\ref{geodesic_distance}) we are not interested in this case.
For curvature weighted metrics the above equation reads as
\begin{align*}
\Ph&=1+A\Tr(L)^{2k}: \qquad \frac{4^k A (k-1) }{r^{4 k}}=1
\end{align*}
Solutions to these equations are given by
\begin{align*}
\Ph&=1+A\Tr(L)^{2k}: \qquad r=2\sqrt[2k]{A(k-1)}, \quad k\geq 1.
\end{align*}

\begin{figure}[ht]
\centering
\includegraphics[width=.72\textwidth]{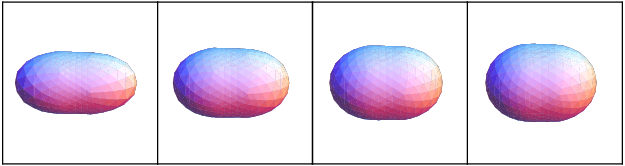}%
\caption{Middle figure of a geodesic between two unit spheres translated by distance $3$ for $\Ph=1+A \Tr(L)^2$. From left to right: $A=0.2$, $A=0.4$, $A=0.6$, $A=0.8$. In each of the simulations 20 timesteps and a triangulation with 720 triangles were used. }
\label{fig:numerics:move:ga_cigar}
\end{figure}

For the most prominent example, the $G^A$--metric, this yields $r=0$, and therefore translation can never be a geodesic for this type of metric. The numerics have shown that the $G^A$--metric  yields geodesics that resemble the geodesics of the $G^A$--metric for planar curves from \cite[section~5.2]{Michor98}. Namely, when the two spheres are sufficiently far apart, the geodesic passes through a cigar-like middle shape, 
see figure~\ref{fig:numerics:move:ga_cigar}. As predicted by the theory (see section~\ref{fr:fr2}),  geodesics for very high distances tend to have  behavior similar to that  of $\Vol^k$ metrics; i.e.,
the geodesic first shrinks the sphere, then moves it, and then blows it up again 
(cf. section~\ref{shrinkandblow}). 

\begin{figure}[ht]
\centering
\includegraphics[width=.9\textwidth]{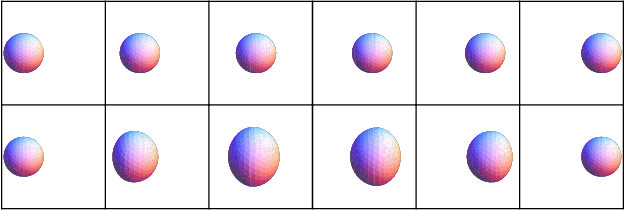}
\caption{Geodesic between two unit spheres translated by distance $3$ for $\Ph=1+\frac{1}{16}\Tr(L)^4$ (first row) and
$\Ph=1+\Tr(L)^6$ (second row). In each of the experiments 20 timesteps and a triangulation with 720 triangles were used. Time progresses from left to right. 
Boundary shapes $t=0$ and $t=1$ are not included.}
\label{fig:numerics:move:trace4andtrace6}
\end{figure}

For metrics weighted by higher factors of mean curvature, 
the above equation for the radius has a positive solution.
For these metrics, geodesics for translations tend to
scale the sphere until it has reached the optimal radius and then
translate it. If the radius is already optimal, the resulting geodesic is a pure 
translation (see figure~\ref{fig:numerics:move:trace4andtrace6}). 

\begin{figure}[ht]
\centering
\includegraphics[width=.8\textwidth]{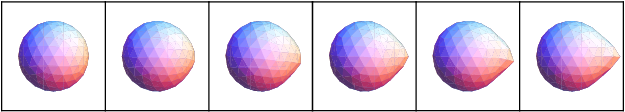}
\caption{Geodesic between a sphere and a sphere with a small bump for $\Ph=\Vol$. 20 timesteps and a triangulation with 500 triangles were used. Time progresses from left to right.}
\label{fig:numerics:deform:minibump_vol}
\end{figure}

If the distance is not high enough a scaling towards the optimal size still occurs, but the middle 
figure is not a perfect sphere anymore. Instead it is an ellipsoid as in 
figure~\ref{fig:numerics:move:trace4andtrace6}.

\subsection{Deformation of a shape}

We will calculate numerically the geodesic between a shape and a deformation of the shape for various almost local metrics. 
Small deformations are handled well by all metrics, and they all yield similar results. 
An example of a geodesic 
resulting in a small deformation can be seen in figure~\ref{fig:numerics:deform:minibump_vol}, where a small bump is 
grown out of a sphere. The energy needed for this deformation is reasonable compared to the energy needed for a pure translation. 
Taking the metric with $\Ph=\Vol$ as an example, growing a bump of size 0.4 as in figure~\ref{fig:numerics:deform:minibump_vol} costs 
about a third of a translation of the sphere by 0.4. 

\begin{figure}[ht]
\centering
\includegraphics[width=\textwidth]{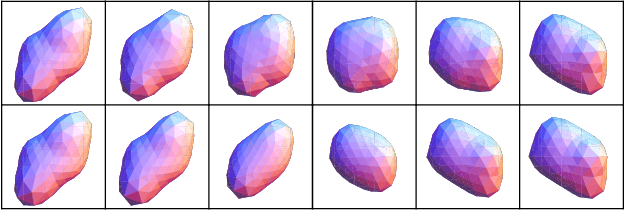}
\caption{Large deformation of a shape for $\Ph=\Vol$ and $\Ph=e^{\Vol}$. 20 timesteps and a triangulation with 500 triangles were used. Time progresses from left to right.}
\label{fig:numerics:deform:large_deformation_conformal}
\end{figure}

\begin{figure}[ht]
\includegraphics[width=\textwidth]{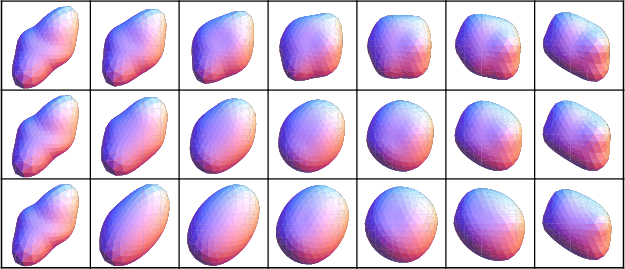}
\caption{Large deformation of a shape for  $\Ph=1+0.1\Tr(L)^2$ (top), $\Ph=1+10\Tr(L)^2$ (middle),
and $\Ph=1+\Tr(L)^6$ (bottom). 20 timesteps and a triangulation with 720 triangles were used. Time progresses from left to right.}
\label{fig:numerics:deform:trace}
\end{figure}

Bigger deformations work well with $\Vol^k$-metrics and curvature weighted metrics, but not with the $e^{\Vol}$-metric, 
which tends to shrink the object and to concentrate almost all of the deformation at a single time step. 
In figure~\ref{fig:numerics:deform:large_deformation_conformal}, a large deformation can be seen for the case of $\Ph=\Vol$ and $\Ph=e^{\Vol}$. 
Clearly one can see that the $e^{\Vol}$-metric concentrates almost all of the deformation in a single time step. 
We have met this misbehavior of the $e^{\Vol}$-metric already with translations. Again, the reason is that $e^{\Vol}$ is so sensitive to changes in volume.
In figure~\ref{fig:numerics:deform:trace} one sees that geodesics are smoothed further by higher curvature weights.

\bibliographystyle{plain}

\begin{thebibliography}{10}

\bibitem{Michor120}
M.~Bauer, P.~Harms, and P.~W. Michor.
\newblock Curvature weighted metrics on shape space of hypersurfaces in n-space.
\newblock To appear in: Differential Geometry and its Applications.
\url{arXiv:1102.0678}.

\bibitem{Michor119}
M.~Bauer, P.~Harms, and P.~W. Michor.
\newblock Sobolev metrics on shape space of surfaces in n-space.
\newblock To appear in: Journal of Geometric Mechanics. \url{arXiv:1009.3616}.

\bibitem{Bauer2010}
Martin Bauer.
\newblock {\em Almost local metrics on shape space of surfaces}.
\newblock PhD thesis, University of Vienna, 2010.

\bibitem{Besse2008}
Arthur~L. Besse.
\newblock {\em Einstein manifolds}.
\newblock Classics in Mathematics. Springer-Verlag, Berlin, 2008.

\bibitem{Michor40}
V.~Cervera, F.~Mascar{\'o}, and P.~W. Michor.
\newblock The action of the diffeomorphism group on the space of immersions.
\newblock {\em Differential Geom. Appl.}, 1(4):391--401, 1991.

\bibitem{Fourer2002}
R.~Fourer and B.~W. Kernighan.
\newblock {\em AMPL: A Modeling Language for Mathematical Programming}.
\newblock Duxbury Press, 2002.

\bibitem{Fulton1997}
William Fulton.
\newblock {\em Young tableaux}, volume~35 of {\em London Mathematical Society
  Student Texts}.
\newblock Cambridge University Press, 1997.

\bibitem{Harms2010}
Philipp Harms.
\newblock {\em Sobolev metrics on shape space of surfaces}.
\newblock PhD thesis, University of Vienna, 2010.

\bibitem{Kobayashi1996a}
Shoshichi Kobayashi and Katsumi Nomizu.
\newblock {\em Foundations of differential geometry. {V}ol. {I}}.
\newblock Wiley Classics Library. John Wiley \& Sons Inc., New York, 1996.

\bibitem{MichorF}
I.~Kol{\'a}{\v{r}}, P.~W. Michor, and J.~Slov{\'a}k.
\newblock {\em Natural operations in differential geometry}.
\newblock Springer-Verlag, Berlin, 1993.

\bibitem{MichorG}
Andreas Kriegl and Peter~W. Michor.
\newblock {\em The convenient setting of global analysis}, volume~53 of {\em
  Mathematical Surveys and Monographs}.
\newblock American Mathematical Society, Providence, RI, 1997.

\bibitem{MennucciYezzi2008}
A.~Mennucci, A.~Yezzi, and G.~Sundaramoorthi.
\newblock Properties of {S}obolev-type metrics in the space of curves.
\newblock {\em Interfaces Free Bound.}, 10(4):423--445, 2008.

\bibitem{MichorH}
Peter~W. Michor.
\newblock {\em Topics in differential geometry}, volume~93 of {\em Graduate
  Studies in Mathematics}.
\newblock American Mathematical Society, Providence, RI, 2008.

\bibitem{Michor102}
Peter~W. Michor and David Mumford.
\newblock Vanishing geodesic distance on spaces of submanifolds and
  diffeomorphisms.
\newblock {\em Doc. Math.}, 10:217--245 (electronic), 2005.

\bibitem{Michor98}
Peter~W. Michor and David Mumford.
\newblock Riemannian geometries on spaces of plane curves.
\newblock {\em J. Eur. Math. Soc. (JEMS) 8 (2006), 1-48}, 2006.

\bibitem{Michor107}
Peter~W. Michor and David Mumford.
\newblock An overview of the {R}iemannian metrics on spaces of curves using the
  {H}amiltonian approach.
\newblock {\em Appl. Comput. Harmon. Anal.}, 23(1):74--113, 2007.

\bibitem{Salvai2009}
Marcos Salvai.
\newblock Geodesic paths of circles in the plane.
\newblock {\em Revista Matem{\'a}tica Complutense}, 2009.

\bibitem{Shah2008}
Jayant Shah.
\newblock {$H\sp 0$}-type {R}iemannian metrics on the space of planar curves.
\newblock {\em Quart. Appl. Math.}, 66(1):123--137, 2008.

\bibitem{Sullivan2008}
John~M. Sullivan.
\newblock Curvatures of smooth and discrete surfaces.
\newblock In {\em Discrete differential geometry}, volume~38 of {\em
  Oberwolfach Semin.}, pages 175--188. Birkh\"auser, Basel, 2008.

\bibitem{Thompson1942}
D'Arcy Thompson.
\newblock {\em On Growth and Form}.
\newblock Cambridge University Press, 1942.

\bibitem{Verpoort2008}
Steven Verpoort.
\newblock {\em The geometry of the second fundamental form: Curvature
  properties and variational aspects}.
\newblock PhD thesis, Katholieke Universiteit Leuven, 2008.

\bibitem{Waechter2002}
A.~W\"achter.
\newblock {\em An Interior Point Algorithm for Large-Scale Nonlinear
  Optimization with Applications in Process Engineering}.
\newblock PhD thesis, Carnegie Mellon University, 2002.

\bibitem{YezziMennucci2004}
A.~Yezzi and A.~Mennucci.
\newblock Conformal riemannian metrics in space of curves.
\newblock EUSIPCO, 2004.

\bibitem{YezziMennucci2004a}
A.~{Yezzi} and A.~{Mennucci}.
\newblock {Metrics in the space of curves}.
\newblock \url{arXiv:math/0412454}, December 2004.

\bibitem{YezziMennucci2005}
Anthony Yezzi and Andrea Mennucci.
\newblock Conformal metrics and true "gradient flows" for curves.
\newblock In {\em Proceedings of the Tenth IEEE International Conference on
  Computer Vision}, volume~1, pages 913--919, Washington, 2005. IEEE Computer
  Society.

\end{thebibliography}

\begin{landscape}
\appendix
\section*{The AMPL model file}\label{ampl}
\lstset{emph={param,sum,var,sum,if,else,minimize,then,sqrt,abs,acos,set},emphstyle=\textbf, extendedchars=true,
  backgroundcolor=\color[gray]{0.9},
  numbers=left, numberstyle=\scriptsize, stepnumber=2, numbersep=5pt,
  xleftmargin=12pt, xrightmargin=12pt,
  commentstyle=\small,columns=flexible,
  showstringspaces=false}

\begin{lstlisting}[basicstyle=\small,caption= AMPL model file]
param A default 1;
param k default 1;
param B default 1;
param l default 1;
param TimestepsN > 1 integer;
param VerticesN integer;
param PenaltyFactor default 1;
param PenaltyExponent default 2;
set VerticesI := 1..VerticesN;
set VerticesOfEdgesI within {VerticesI,VerticesI};
set VerticesOfFacesI within {VerticesI,VerticesI,VerticesI};
set FacesOfVerticesI {v in VerticesI} within VerticesOfFacesI;
set LinkOfVerticesI {VerticesI} within {VerticesOfFacesI,VerticesOfEdgesI,{-1,1}}; 
set AdjacentEdgesOfVerticesI {VerticesI} within {VerticesOfEdgesI,{1,-1},VerticesOfEdgesI,{1,-1}};
set EdgesOfFacesI {VerticesOfFacesI} within VerticesOfEdgesI;
set EdgesOfVerticesI {v in VerticesI} := setof {(f1,f2,f3,e1,e2,o) in LinkOfVerticesI[v]}(e1,e2);

param Pi default 3.141592653589793;
param PerfectAngle {v in VerticesI} default cos(2*Pi/card(FacesOfVerticesI[v]));
param InitialVertices {VerticesI,1..3};
param FinalVertices {VerticesI,1..3};

var MiddleVertices {2..TimestepsN,VerticesI,1..3};

var Vertices {t in 1..TimestepsN+1,v in VerticesI,i in 1..3} = 
  (if t=1 then InitialVertices[v,i] 
   else if t=TimestepsN+1 then FinalVertices[v,i] 
   else MiddleVertices[t,v,i]);

var VectorOfEdges {t in 1..TimestepsN+1, (v1,v2) in VerticesOfEdgesI,i in 1..3} =
	 Vertices[t,v2,i] - Vertices[t,v1,i];

var LengthOfEdges {t in 1..TimestepsN+1, (v1,v2) in VerticesOfEdgesI} =
	  sqrt(VectorOfEdges[t,v1,v2,1]^2+VectorOfEdges[t,v1,v2,2]^2+VectorOfEdges[t,v1,v2,3]^2);

var CrossOfFaces {t in 1..TimestepsN+1,(v1,v2,v3) in VerticesOfFacesI,i in 1..3} =
  if i=1 then (Vertices[t,v2,2]-Vertices[t,v1,2])*(Vertices[t,v3,3]-Vertices[t,v1,3]) - 
              (Vertices[t,v2,3]-Vertices[t,v1,3])*(Vertices[t,v3,2]-Vertices[t,v1,2]) 
  else if i=2 then -(Vertices[t,v2,1]-Vertices[t,v1,1])*(Vertices[t,v3,3]-Vertices[t,v1,3]) +
                    (Vertices[t,v2,3]-Vertices[t,v1,3])*(Vertices[t,v3,1]-Vertices[t,v1,1]) 
  else (Vertices[t,v2,1]-Vertices[t,v1,1])*(Vertices[t,v3,2]-Vertices[t,v1,2]) - 
       (Vertices[t,v2,2]-Vertices[t,v1,2])*(Vertices[t,v3,1]-Vertices[t,v1,1]) ;

var NormCrossOfFaces {t in 1..TimestepsN+1,(v1,v2,v3) in VerticesOfFacesI} =
  sqrt(CrossOfFaces[t,v1,v2,v3,1]^2 + CrossOfFaces[t,v1,v2,v3,2]^2 + CrossOfFaces[t,v1,v2,v3,3]^2);

var NuOfFaces {t in 1..TimestepsN+1,(v1,v2,v3) in VerticesOfFacesI,i in 1..3} =
  CrossOfFaces[t,v1,v2,v3,i]/NormCrossOfFaces[t,v1,v2,v3];

var AreaOfFaces {t in 1..TimestepsN+1,(v1,v2,v3) in VerticesOfFacesI} =
  NormCrossOfFaces[t,v1,v2,v3]/2;

var AreaOfVertices {t in 1..TimestepsN+1, v in VerticesI} =
  (sum {(f1,f2,f3) in FacesOfVerticesI[v]} AreaOfFaces[t,f1,f2,f3])/3;

var VectorAreaOfVertices {t in 1..TimestepsN+1, v in VerticesI, i in 1..3} = 
  (sum {(v1,v2,v3) in FacesOfVerticesI[v]} CrossOfFaces[t,v1,v2,v3,i])/6; 

var SquareOfNormOfVectorAreaOfVertices {t in 1..TimestepsN+1, v in VerticesI} = 
  VectorAreaOfVertices[t,v,1]^2+VectorAreaOfVertices[t,v,2]^2+VectorAreaOfVertices[t,v,3]^2;

var NormOfVectorAreaOfVertices {t in 1..TimestepsN+1, v in VerticesI} = 
  sqrt( SquareOfNormOfVectorAreaOfVertices[t,v]);

var Volume {t in 1..TimestepsN+1} =
  sum{(v1,v2,v3) in VerticesOfFacesI} AreaOfFaces[t,v1,v2,v3];

var VectorMeanCurvatureOfVertices {t in 1..TimestepsN+1, v in VerticesI, i in 1..3} = 
  if i=1 then 
    sum {(f1,f2,f3,e1,e2,o) in LinkOfVerticesI[v]} o*
      ( VectorOfEdges[t,e1,e2,2]*NuOfFaces[t,f1,f2,f3,3] - 
        VectorOfEdges[t,e1,e2,3]*NuOfFaces[t,f1,f2,f3,2] ) 
  else if i=2 then 
    sum {(f1,f2,f3,e1,e2,o) in LinkOfVerticesI[v]} o*
      (-VectorOfEdges[t,e1,e2,1]*NuOfFaces[t,f1,f2,f3,3] + 
        VectorOfEdges[t,e1,e2,3]*NuOfFaces[t,f1,f2,f3,1] ) 
  else 
    sum {(f1,f2,f3,e1,e2,o) in LinkOfVerticesI[v]} o*
      ( VectorOfEdges[t,e1,e2,1]*NuOfFaces[t,f1,f2,f3,2] - 
        VectorOfEdges[t,e1,e2,2]*NuOfFaces[t,f1,f2,f3,1] ) ;

var SquareOfScalarMeanCurvatureOfVertices {t in 1..TimestepsN+1, v in VerticesI} = 
  (VectorMeanCurvatureOfVertices[t,v,1]^2+VectorMeanCurvatureOfVertices[t,v,2]^2
    +VectorMeanCurvatureOfVertices[t,v,3]^2)/SquareOfNormOfVectorAreaOfVertices[t,v];

var PhiOfVertices {t in 1..TimestepsN+1,v in VerticesI} =
  1+ A*(SquareOfScalarMeanCurvatureOfVertices[t,v])^k +B*(Volume[t])^l;

var IncrementsOfVertices {t in 1..TimestepsN,v in VerticesI,i in 1..3} =
  TimestepsN*(Vertices[t+1,v,i] - Vertices[t,v,i]);

var Energy = 1/ 12 / TimestepsN *(
  sum {t in 1..TimestepsN,v in VerticesI} 
    PhiOfVertices[t,v] * sum {(w1,w2,w3) in FacesOfVerticesI[v]}
        ( IncrementsOfVertices[t,v,1]*CrossOfFaces[t,w1,w2,w3,1] +
          IncrementsOfVertices[t,v,2]*CrossOfFaces[t,w1,w2,w3,2] +
          IncrementsOfVertices[t,v,3]*CrossOfFaces[t,w1,w2,w3,3] )^2 /
        NormCrossOfFaces[t,w1,w2,w3] +
  sum {t in 1..TimestepsN,v in VerticesI} 
    PhiOfVertices[t+1,v] * sum {(w1,w2,w3) in FacesOfVerticesI[v]}
        ( IncrementsOfVertices[t,v,1]*CrossOfFaces[t+1,w1,w2,w3,1] +
          IncrementsOfVertices[t,v,2]*CrossOfFaces[t+1,w1,w2,w3,2] +
          IncrementsOfVertices[t,v,3]*CrossOfFaces[t+1,w1,w2,w3,3] )^2 /
        NormCrossOfFaces[t+1,w1,w2,w3] );

var Penalty = 
  sum {t in 1..TimestepsN+1, v in VerticesI,(v1,w1,o1,v2,w2,o2) in AdjacentEdgesOfVerticesI[v]} 
    abs(
      ( VectorOfEdges[t,v1,w1,1]*VectorOfEdges[t,v2,w2,1] +
        VectorOfEdges[t,v1,w1,2]*VectorOfEdges[t,v2,w2,2] +
        VectorOfEdges[t,v1,w1,3]*VectorOfEdges[t,v2,w2,3] ) * o1 * o2 
      / LengthOfEdges[t,v1,w1] / LengthOfEdges[t,v2,w2]
      - PerfectAngle[v]
    )^PenaltyExponent;

minimize f: 
  Energy+Penalty*PenaltyFactor;

\end{lstlisting}
\end{landscape}

\end{document}